\newif\ifPreprint\Preprinttrue
\ifPreprint
  \documentclass[a4paper]{amsart}
  \usepackage[a4paper,margin=22mm]{geometry}
\else
  \PassOptionsToPackage{leqno}{amsmath}
  \documentclass[acrhead]{ws-book9x6}
  \usepackage{ws-book-van}
  
  \def\theoremstyle#1{}
  \def\swapnumbers{}
\fi
\usepackage{microtype}
\usepackage[svgnames]{xcolor}
\ifPreprint
  \usepackage{amsmath}
  \usepackage{amssymb}
\else
\fi

\usepackage{enumitem}
\setlist[enumerate]{topsep=0mm}

\usepackage{todonotes}

\usepackage{mathtools}
\usepackage{mathrsfs}
\usepackage{etoolbox} 
\usepackage{stmaryrd}
\usepackage{booktabs} 
\DeclareMathAlphabet{\mathpzc}{OT1}{pzc}{m}{it}
\DeclareMathAlphabet{\mathcal}{OMS}{cmsy}{m}{n}
\DeclareSymbolFontAlphabet{\mathbb}{AMSb}

\usepackage{hyperref}

\hypersetup{pdftitle={Cyclotomic quiver Hecke algebras of type A},
  pdfauthor={Andrew Mathas},
  colorlinks = true,
  linkcolor =blue,
  anchorcolor = red,
  citecolor = DarkGreen,
  urlcolor = MediumBlue
}

\makeatletter
\newcommand\Autoref[1]{\@first@ref#1,@}
\def\@throw@dot#1.#2@{#1}
\def\@set@refname#1{
    \edef\@tmp{\getrefbykeydefault{#1}{anchor}{}}%
    \def\@refname{\@nameuse{\expandafter\@throw@dot\@tmp.@autorefname}s}%
}
\def\@first@ref#1,#2{%
  \ifx#2@\autoref{#1}\let\@nextref\@gobble
  \else%
    \@set@refname{#1}
    \@refname~\ref{#1}\let\@nextref\@next@ref%
  \fi%
  \@nextref#2%
}
\def\@next@ref#1,#2{%
   \ifx#2@ and~\ref{#1}\let\@nextref\@gobble
   \else, \ref{#1}
   \fi%
   \@nextref#2%
}
\makeatother

\usepackage{cite}

\newcommand\Comment[2][AM]{\space\par\medskip\noindent%
   \fbox{\begin{minipage}{\textwidth}\textbf{Comment\ifx\relax#1\else---#1\fi}\newline%
        #2\end{minipage}}\medskip
}
\usepackage{tikz}
\usetikzlibrary{matrix,calc,intersections,through,backgrounds,arrows}
\usetikzlibrary{patterns,decorations.markings}
\usepgflibrary{arrows}

\makeatletter
\newdimen\LineSpace
\newdimen\LineWidth
\tikzset{line space/.code={\LineSpace=#1}, line space=3pt}
\tikzset{pattern line width/.code={\LineWidth=#1}, pattern line width=.4pt}
\pgfdeclarepatternformonly[\LineSpace,\LineWidth]{my north east lines}
  {\pgfqpoint{-\LineWidth}{-\LineWidth}}{\pgfqpoint{\LineSpace}{\LineSpace}}
  {\pgfqpoint{\LineSpace}{\LineSpace}}%
  {
  \pgfsetcolor{\tikz@pattern@color}
  \pgfsetlinewidth{\LineWidth}
  \pgfpathmoveto{\pgfqpoint{-\LineWidth}{-\LineWidth}}
  \pgfpathlineto{\pgfqpoint{\LineSpace + 0.1pt}{\LineSpace + 0.1pt}}
  \pgfusepath{stroke}
  }
\pgfdeclarepatternformonly[\LineSpace,\LineWidth]{my north west lines}
  {\pgfqpoint{-\LineWidth}{-\LineWidth}}{\pgfqpoint{\LineSpace}{\LineSpace}}
  {\pgfqpoint{\LineSpace}{\LineSpace}}%
  {
  \pgfsetcolor{\tikz@pattern@color}
  \pgfsetlinewidth{\LineWidth}
  \pgfpathmoveto{\pgfqpoint{-\LineWidth}{\LineSpace}}
  \pgfpathlineto{\pgfqpoint{\LineSpace + 0.1pt}{-\LineWidth}}
  \pgfusepath{stroke}
  }
\makeatother

\newcount\tableauRow
\newcount\tableauCol
\def\tableauDim{0.4}
\def\TikTableau[#1][#2](#3){%
  \begin{tikzpicture}[scale=\tableauDim,draw/.append style={thick,black},baseline=#1]
    \tableauRow=0
    \ifx\relax#2\else
      \foreach \r/\c in {#2} {
        \filldraw[blue!55,opacity=0.5](\r,\c)+(-.5,-.5)rectangle++(.5,.5);
      }
    \fi
    \foreach \Row in {#3} {
       \tableauCol=1
       \foreach\k in \Row {
         \draw[thin](\the\tableauCol,\the\tableauRow)+(-.5,-.5)rectangle++(.5,.5);
         \draw[thin](\the\tableauCol,\the\tableauRow)node{$\k$};
         \global\advance\tableauCol by 1
       }
       \global\advance\tableauRow by -1
    }
  \end{tikzpicture}
}
\def\Tableau(#1){\TikTableau[-3.6mm][{}]({#1})}

\def\multitab#1|#2{%
  \,\TikTableau[-1mm][{}]({{#1}})\ifx:#2\,\else\,\big|\,\multitab#2\fi}
\def\onetab(#1){{\def\tableauDim{0.38}\big(\multitab#1|:\big)}}
\def\otab(#1){{\def\tableauDim{0.38}%
  \tikzset{every node/.append style={font=\fontsize{7}{7}\selectfont}}%
  \TikTableau[-1mm][{}]({{#1}})}%
}
\def\ootab(#1){{\def\tableauDim{0.34}%
  \tikzset{every node/.append style={font=\fontsize{8}{8}\selectfont}}%
  \TikTableau[-1mm][{}]({{#1}})}
}
\def\Otab[#1](#2){\TikTableau[-1mm][{#1}]({{#2}})}
\def\tableau(#1){{\def\tableauDim{0.38}%
  \tikzset{every node/.append style={font=\fontsize{8}{8}\selectfont}}%
  \TikTableau[-2.8mm][{}]({#1})%
}}

\def\Dots{\raisebox{-1.5mm}{\scriptsize${\cdot}{\cdot}{\cdot}$}}

\newcommand\YoungDiagram[2][\relax]{
  \begin{tikzpicture}[scale=0.43,draw/.append style={thick,black},baseline=-3mm]
    \ifx\relax#1\relax%
    \else 
    \foreach\box in {#1} {
      \filldraw[blue!30]\box rectangle ++(1,1);
    }
    \fi
    \newcount\tableauRow
    \tableauRow=0
    \foreach \diagramCol in {#2} {
      \draw[thin](1,\the\tableauRow)grid ++(\diagramCol,1);
      \global\advance\tableauRow by -1
    }
  \end{tikzpicture}
}

\newcommand\Bigger[2][7]{\left#2\rule{0mm}{#1truemm}\right.}

\def\TriDiagram(#1|#2|#3){\Bigger(%
  \YoungDiagram{#1}\Bigger|\
  \YoungDiagram{#2}\Bigger|\YoungDiagram{#3}
  \Bigger)
}

\def\Tritab(#1|#2|#3){\Bigger(
             \TikTableau[-5mm][{}]({#1})\,\Bigger|\,%
             \TikTableau[-5mm][{}]({#2})\,\Bigger|\,%
             \TikTableau[-5mm][{}]({#3})
             \Bigger)%
}

\colorlet{darkgreen}{green!50!black}
\tikzset{dots/.style={very thick,loosely dotted},
         greendot/.style={fill,circle,color=darkgreen,inner sep=2pt,outer sep=0}
}
\def\greendot(#1,#2){\node[greendot] at(#1,#2){}}

\newenvironment{braid}{
  \begin{tikzpicture}[baseline=7mm,blue,line width=1.5pt, scale=0.35,
                      draw/.append style={rounded corners},
                      every node/.append style={font=\fontsize{5}{5}\selectfont}]%
  }{\end{tikzpicture}
}

\def\BraidHeight{4}  
\newcommand\YCrossing[5]{%
  \coordinate (A) at (#2,\BraidHeight);
  \coordinate (B) at (#3,\BraidHeight);
  \coordinate (C) at (#4,0);
  \coordinate (D) at (#5,0);
  \draw[blue](A)node[above]{$#1$};
  \draw[blue](B)node[above]{$#1$};
  \draw[blue](C);
  \draw[blue](D);
  \path[name path=AD](A)--(D);
  \path[name path=BC](B)--(C);
  \path[name intersections={of=AD and BC,by=I}];
  \draw[blue](A)--($ (A)!0.8!(I) $) -- ($ (I)!0.2!(C) $)--(C);
  \draw[blue](B)--($ (B)!0.8!(I) $) -- ($ (I)!0.2!(D) $)--(D);
}

\let\Sect=\S
\let\<=\langle
\let\>=\rangle
\def\({\big(}
\def\){\big)}
\def\bijection{\overset{\sim}{\longrightarrow}}
\let\surjection\twoheadrightarrow
\def\Sum{\displaystyle\sum}

\newcommand{\C}{\mathbb{C}}
\def\D{\mathscr D}
\newcommand{\F}{\mathbb{F}}
\newcommand{\Fp}[1][p]{\mathbb{F}_{#1}}
\renewcommand\L[1][n]{\mathscr{L}_{#1}}
\newcommand{\N}{\mathbb{N}}
\def\O{\mathcal{O}}
\def\Zcal{\mathcal{Z}}
\def\p{\mathfrak{p}}
\newcommand{\Q}{\mathbb{Q}}
\newcommand{\Sym}{\mathfrak S}
\newcommand{\Z}{\mathbb{Z}}
\newcommand{\A}{\mathcal{A}}
\renewcommand{\AA}{\mathbb{A}}

\newcommand\Fun{\textsf{F}_n}

\newcommand\adj[1][F]{\mathbf{a}^{#1}_q}
\newcommand\car{\mathbf{c}_q}
\newcommand\dec{\mathbf{d}_q}
\newcommand\deck{\dec^{\mathcal{K}}}
\newcommand\dinv{\mathbf{e}^{\mathcal{K}}_q}
\newcommand\Mull{\mathbf{m}}

\def\map#1#2{\,{:}\,#1\!\longrightarrow\!#2}
{\catcode`\|=\active
  \gdef\set#1{\mathinner{\lbrace\,{\mathcode`\|"8000%
                                   \let|\midvert #1}\,\rbrace}}
}
{\catcode`\|=\active
  \gdef\gen#1{\mathinner{\langle\,{\mathcode`\|"8000%
                                   \let|\midvert #1}\,\rangle}}
}
\def\midvert{\egroup\,\mid\,\bgroup}

\def\pmod#1{\text{ }(\textrm{mod } #1)\,}

\def\Hcal{\mathscr{H}}
\renewcommand\H[1][n]{\Hcal^\Lambda_{#1}}
\newcommand\UnH[1][n]{\underline{\Hcal}^\Lambda_{#1}}
\newcommand\Hlam[1][\blam]{\Hcal^{\gdom#1}_n}
\newcommand\Hlamp[1][\blam]{\Hcal^{\prime\ldom#1}_n}
\newcommand\HK{\Hcal^K_n}
\newcommand\HO{\Hcal^\O_n}
\def\Rcal{\mathscr{R}}
\newcommand\Rn[1][n]{\Rcal_{#1}}
\newcommand\R[1][n]{\Rcal^\Lambda_{#1}}
\newcommand\Rp[1][\beta']{\Rcal^{\Lambda'}_{#1}}

\newcommand\Poly{\mathscr{C}_n}
\newcommand\Fock[1][\Q(q)]{\mathscr{F}^\Lambda_{#1}}
\newcommand\Cry{\mathscr{L}^\Lambda_0}
\newcommand\fock[1][\blam]{{\mid}#1{\>}}

\def\Uq{\cramped{U_q(\widehat{\mathfrak{sl}}_e)}}
\def\Usl{U(\widehat{\mathfrak{sl}}_e)}
\def\UA{U_\A(\widehat{\mathfrak{sl}}_e)}

\def\dyo[#1]{y^{\<#1\>}}

\newcommand\yo{y^{\O}}
\def\psio{\psi^\O}
\def\eps{1+\rho_r(\bi)}
\def\spe{1-\rho_r(\bi)}
\newcommand\fo[1][\bi]{f_{#1}^\O}

\def\Pcal{\mathcal{P}}
\newcommand{\Klesh}[1][n]{\mathcal{K}^{\Lambda}_{#1}}
\newcommand{\Parts}[1][n]{\Pcal^\Lambda_{#1}}
\newcommand\blam{{\boldsymbol\lambda}}

\newcommand\balpha{{\boldsymbol\alpha}}
\newcommand\bbeta{{\boldsymbol\beta}}
\newcommand\bgamma{{\boldsymbol\gamma}}
\newcommand\charge{{\boldsymbol\kappa}}
\newcommand\bmu{{\boldsymbol\mu}}

\newcommand\bnu{{\boldsymbol\nu}}

\newcommand\bsig{{\boldsymbol\sigma}}
\newcommand\btau{{\boldsymbol\tau}}
\def\Zero{\mathbf{0}_\ell}
\def\zero{v_\Lambda}

\def\diag#1{\llbracket#1\rrbracket}
\newcommand\Belt{\mathbf{B}_A}

\def\i{\hat\imath}

\newcommand\bi{\mathbf{i}}
\newcommand\bj{\mathbf{j}}
\newcommand\bx{\mathbf{x}}
\newcommand\Cont[1][n]{\mathscr{C}_{#1}}

\def\ilam{\bi^\blam}
\def\imu{\bi^\bmu}

\newcommand\iarrow[1][i]{\mathrel{\raisebox{-0.5mm}{%
   $\xrightarrow{\raisebox{-0.5mm}[0mm][0mm]{\scriptsize\ $#1$\ }}$}%
}}
\newcommand\goodarrow{\iarrow[\text{good}]}

\let\uold=\u
\def\a{{\mathsf a}}
\def\b{{\mathsf b}}
\def\c{{\mathsf c}}
\def\d{{\mathsf d}}
\def\s{{\mathsf s}}
\def\t{{\mathsf t}}
\def\u{{\mathsf u}}
\def\v{{\mathsf v}}
\def\T{{\mathsf T}}
\def\V{{\mathsf V}}

\def\tlam{\t^\blam}
\def\tllam{\t_\blam}
\def\tmu{\t^\bmu}

\def\Std{\mathop{\rm Std}\nolimits}

\def\SStd{\mathop{\rm Std}^2\nolimits}

\let\gedom=\trianglerighteq
\let\gdom=\vartriangleright
\let\Gdom=\blacktriangleright

\makeatletter
\newcommand{\superimpose}[2]{%
  {\ooalign{$#1\@firstoftwo#2$\cr\hfil$#1\@secondoftwo#2$\hfil\cr}}}
\makeatother

\newcommand{\GDom}{\mathop{\mathpalette\superimpose{{\gdom}{\cdot\,}}}}
\newcommand{\notGDom}{\mathop{\mathpalette\superimpose{{\not\gdom}{\cdot\,}}}}

\def\UnC{\underline{C}}
\def\UnD{\underline{D}}
\def\UnG{\underline{G}}
\def\UnP{\underline{P}}
\def\UnS{\underline{S}}

\newcommand{\DeclareMyOperator}[1]{%
  \expandafter\DeclareMathOperator\csname #1\endcsname{#1}
}
\newcommand{\DeclareOperators}{\forcsvlist{\DeclareMyOperator}}

\DeclareOperators{codeg,comp,Deg,End,Hom,Ind,Mat,Res,rad,Shape,res,Add,Rem}

\DeclareMathOperator\noedge{\:\rlap{\hspace*{0.25em}/}\text{---}\:}

\def\gdim{\mathop{\mathpzc{dim}}\nolimits_q}
\DeclareMathOperator\ZHom{\mathpzc{Hom}}

\newcommand\rep[1][]{\mathop{\rm Rep}\nolimits_{#1}}
\newcommand\repH[1][n]{[\rep(\H[#1])]}
\newcommand\proj[1][]{\mathop{\rm Proj}\nolimits_{#1}}
\newcommand\projH[1][n]{[\proj(\H[#1])]}
\newcommand\Rep[1][\Q(q)]{[\mathop{\rm Rep}\nolimits^\Lambda_{#1}]}
\newcommand\UnRep[1][\Q]{\underline{\mathop{\rm Rep}\nolimits}^\Lambda_{#1}}
\newcommand\Proj[1][\Q(q)]{[\mathop{\rm Proj}\nolimits^\Lambda_{#1}]}
\newcommand\UnProj[1][\Q]{\underline{\mathop{\rm Proj}\nolimits}^\Lambda_{#1}}

\newcommand{\qbinom}[3][q]{\genfrac{\llbracket}{\rrbracket}{0pt}{}{#2}{#3}_{#1}}

\DeclareMathOperator\Gedom{\,{\underline{\kern-.1ex{\blacktriangleright}\kern-0.1ex}}\,}
\DeclareMathOperator\defect{def}
\DeclareMathOperator\wt{wt}
\DeclareMathOperator{\Ch}{Ch\text{${}_q$}}
\DeclareMathOperator{\ch}{ch\text{${}_q$}}
\DeclareMathOperator{\iInd}{\text{$i$}-Ind}
\DeclareMathOperator{\iRes}{\text{$i$}-Res}

\def\sgn{\mathtt{sgn}}

\let\ledom=\trianglelefteq
\let\ldom=\vartriangleleft

\def\ceil#1/#2{\lceil\tfrac{#1}{#2}\rceil}

\usepackage{aliascnt}
\def\NewTheorem#1{%
  \newaliascnt{#1}{equation}
  \newtheorem{#1}[#1]{#1}
  \aliascntresetthe{#1}
  \expandafter\def\csname #1autorefname\endcsname{#1}
}
\def\equationautorefname~#1\null{(#1)\null}
\makeatletter
\def\subsectionautorefname{\Sect\@gobble}
\makeatother

\swapnumbers
\numberwithin{equation}{subsection}

\NewTheorem{Proposition}
\NewTheorem{Theorem}
\NewTheorem{Corollary}
\NewTheorem{Conjecture}
\NewTheorem{Lemma}
\NewTheorem{Definition}
\theoremstyle{definition}
\theoremstyle{remark}
\NewTheorem{Remark}
\NewTheorem{Remarks}

\newaliascnt{Example}{equation}
\aliascntresetthe{Example}

\newenvironment{Example}[1][\relax]%
  {\refstepcounter{Example}\trivlist
   \ifx#1\relax
     \item[\hskip\labelsep\theequation.~\textbf{Example}\space]
    \else
     \item[\hskip\labelsep\theequation.~\textbf{Example }(#1)\space]
   \fi
   \ignorespaces
  }{\unskip\nobreak\hfil%
    \penalty50\hskip2em\hbox{}\nobreak\hfil$\Diamond$%
    \parfillskip=0pt\finalhyphendemerits=0\penalty-100\endtrivlist%
}
\newaliascnt{Examples}{equation}
\aliascntresetthe{Examples}

  {\refstepcounter{Examples}\trivlist
   \item[\hskip\labelsep\theequation.~\textbf{Examples}\space]
   \ignorespaces #1\enumerate
   }{\endenumerate\unskip\nobreak\hfil%
    \penalty50\hskip2em\hbox{}\nobreak\hfil$\Diamond$%
    \parfillskip=0pt\finalhyphendemerits=0\penalty-100\endtrivlist%
}

\newcounter{MyCase}
\numberwithin{MyCase}{equation}
\newcommand\Case[1]{\refstepcounter{MyCase}
  \medskip\noindent\textbf{Case \arabic{MyCase}.}
  \space\textit{#1.}\newline
}

\def\Step#1.{\medskip\noindent\textit{Step #1.}\space}

\def\Email#1{\email{\href{mailto:#1}{#1}}}

\ifPreprint
  \title{Cyclotomic quiver Hecke algebras of type~$A$}
  \author{Andrew Mathas}
  \address{School of Mathematics and Statistics F07,
  University of Sydney, NSW 2006, Australia}
  \Email{andrew.mathas@sydney.edu.au}
\fi

\begin{document}

%

\ifPreprint
  \maketitle
  \tableofcontents
\else
  \chapter[Mathas --- Cyclotomic quiver Hecke algebras of type~$A$]%
  {Cyclotomic quiver Hecke algebras\\ of type~$A$}
  \centerline{\textsc{Andrew Mathas}}
\fi

\section*{Introduction}

  The \textit{cyclotomic Hecke algebras} of type~$A$ are a much studied
  class of algebras that include, as special cases, the group algebras
  of the symmetric groups and the Iwahori-Hecke algebras of types~$A$
  and~$B$. They have a rich representation theory that can be approached
  using algebraic combinatorics, standard tools of representation
  theory, or via the theory of Lie and algebraic groups, which brings
  deep methods from geometry into play.

  In 2008 Khovanov and Lauda~\cite{KhovLaud:diagI,KhovLaud:diagII} and
  Rouquier~\cite{Rouq:2KM} introduced the \textit{quiver Hecke
  algebras}, or \textit{KLR algebras}.  These are a remarkable family
  $\set{\Rn(\Gamma)|n\ge0}$ of
  $\Z$-graded algebras defined by generators and relations depending on
  a quiver~$\Gamma$. The motivation for defining
  and studying these algebras came, at least in part, from questions in
  geometry and 2-representation theory. The algebras $\Rn(\Gamma)$ categorify
  the negative part of the associated quantum group
  $U_q(\mathfrak{g}_\Gamma)$~\cite{Rouquier:QuiverHecke2Lie,
                                   VaragnoloVasserot:CatAffineKLR}.
  That is, there are natural isomorphisms
  $$  U^-_q(\mathfrak{g}_\Gamma) \cong \bigoplus_{n\ge0}[\proj(\Rn(\Gamma))],$$
  where $[\proj(\Rn(\Gamma))]$ is the Grothendieck group of the category of
  finitely generated graded projective $\Rn(\Gamma)$-modules. For each dominant
  weight~$\Lambda$ the quiver Hecke algebra $\Rn(\Gamma)$ has a cyclotomic
  quotient $\R(\Gamma)$ that categorifies the highest weight module $L(\Lambda)$
  \cite{Rouquier:QuiverHecke2Lie,KangKashiwara:CatCycKLR,Webster:HigherRep}.
  These results can be thought of as far reaching generalizations of Ariki's
  Categorification Theorem in type~$A$~\cite{Ariki:can}.

  Spectacularly, Brundan and Kleshchev~\cite{BK:GradedKL,Rouq:2KM}
  proved that each cyclotomic Hecke algebra of type~$A$ is isomorphic
  to a cyclotomic quotient of a quiver Hecke algebra of
  type~$A$. Thus, the KLR algebras give a new window for understanding
  the cyclotomic Hecke algebras of type~$A$.
  This chapter is an attempt to open this window and show how the
  ``classical'' ungraded representation theory and the emerging graded
  representation theory of the cyclotomic quiver Hecke algebras
  interact.

  With the advent of the KLR algebras the cyclotomic Hecke
  algebras can now be studied from many different perspectives
  including:
  \begin{enumerate}
    \item As ungraded cyclotomic Hecke algebras.
    \item As graded cyclotomic quiver Hecke algebras or KLR algebras.
    \item Geometrically as the ext-algebras of Lusztig
    sheaves~\cite{VaragnoloVasserot:CatAffineKLR,RouquierShanVaragnoloVasserot,%
    Maksimau:CanonicalKLR}.
    \item Through the lens of 2-representation theory using Rouquier's theory of
    2-Kac Moody
    algebras~\cite{Rouq:2KM,KangKashiwara:CatCycKLR,Webster:HigherRep}.
  \end{enumerate}
  Here we focus on (a) and (b) taking an unashamedly combinatorial approach,
  although we will see shadows of geometry and $2$-representation theory.

  For every quiver there is a corresponding family of KLR algebras,
  however, the quiver Hecke algebras attached to the quivers of type $A$
  are special because these are the only quiver Hecke algebras that
  existed in the literature prior to~\cite{KhovLaud:diagI,Rouq:2KM} ---
  all of the other quiver Hecke algebras are ``new'' algebras. In
  type~$A$, when we are working over a field, the quiver Hecke algebras
  are isomorphic to affine Hecke algebras of
  type~$A$~\cite{Rouquier:QuiverHecke2Lie} and the cyclotomic quiver
  Hecke algebras are isomorphic to the cyclotomic Hecke algebras of
  type~$A$~\cite{BK:GradedKL}.  The cyclotomic Hecke algebras of
  type~$A$ have a uniform description but, historically, they have been
  studied either as \textit{Ariki-Koike} algebras $(v\ne1$), or as
  degenerate Ariki-Koike algebras $(v=1$). The existence of gradings on
  Hecke algebras, at least in the ``abelian defect case'', was predicted
  by Rouquier~\cite[Remark~3.11]{Rouquier:DerivICM} and
  Turner~\cite{Turner:Rock}.

  The cyclotomic quiver Hecke algebras of type~$A$ are better understood than
  other types because we already know a lot about the isomorphic, but ungraded,
  cyclotomic Hecke algebras~\cite{M:cyclosurv}. For example, by piggybacking on
  the ungraded theory, homogeneous bases have been constructed for the
  cyclotomic quiver Hecke algebras of type~$A$\cite{HuMathas:GradedCellular} but
  such bases are not yet known in other types. Many of the major results for
  general quiver Hecke algebras were first proved in type~$A$ and then
  generalized to other types. In fact, the type~$A$ algebras, through Ariki's
  theorem and Chuang and Rouquier's seminal work on
  $\mathfrak{sl}_2$-categorifications~\cite{ChuangRouq:sl2}, has motivated many
  of these developments.

  The first section starts by giving a uniform description of the
  degenerate and non-degenerate cyclotomic Hecke algebras, recalling
  some structural results from the ungraded representation theory of
  these algebras. Everything mentioned in this section is applied
  later in the graded setting.

  The second section introduces the cyclotomic KLR algebras as abstract algebras
  given by generators and relations. We use the relations in a series of
  extended examples to give the reader a feel for these algebras. In
  particular, using just the relations we show that the semisimple cyclotomic
  quiver Hecke algebras of type~$A$ are direct sums of matrix rings.  From
  this we deduce Brundan and Kleshchev's Graded Isomorphism
  Theorem in the semisimple case.

  The third section starts with Brundan and Kleshchev's Graded Isomorphism
  Theorem~\cite{BK:GradedKL}. We develop the representation theory of the
  cyclotomic quiver Hecke algebras as graded cellular algebras, focusing on the
  graded Specht modules.  The highlight of this section is a self-contained
  proof of Brundan and Kleshchev's Graded Categorification
  Theorem~\cite{BK:GradedDecomp}, starting from the graded branching rules for
  the graded Specht modules and then using Ariki's Categorification
  Theorem~\cite{Ariki:can} to make the link with canonical bases. We also give a
  new treatment of graded adjustment matrices using a cellular algebra approach.

  In the final section we sketch one way of proving Brundan and Kleshchev's
  Graded Isomorphism Theorem using the classical theory of seminormal forms. As an
  application we describe how to construct a new graded cellular basis for the
  cyclotomic quiver Hecke algebras that appears to have remarkable properties.
  We end with a conjecture for the $q$-characters of the graded simple modules.

  Although the experts will find some new results, most of the novelty
  is in our approach and the arguments that we use. We include many
  examples and a comprehensive survey of the literature. For a different
  perspective we recommend Kleshchev's survey
  article~\cite{KLeshchev:KLRSurvey} on the applications of quiver Hecke
  algebras to symmetric groups.

  \subsection*{Acknowledgements}
  This chapter grew out of a series of lectures at
  the IMS at the National University of Singapore. I thank the
  organizers for inviting me to give these lectures and to write this
  chapter.  The direction taken in these notes, and the conjecture
  formulated in \autoref{S:Conjecture}, is partly motivated by the
  author's joint work with Jun Hu and I thank him for his implicit
  contributions. I thank Susumu Ariki, Anton Evseev, Matthew Fayers, Jun
  Hu, Kai Meng Tan and the referee for their comments on earlier
  versions of this manuscript and Jon Brundan for some helpful
  discussions about crystal bases. This chapter was written while
  visiting Universit\"at Stuttgart and Charles University, Prague.

\section{Cyclotomic Hecke algebras of type~$A$}\label{Sec:CycHeckeAlgebras}
  This section surveys the representation theory of the cyclotomic Hecke
  algebras of type~$A$ and, at the same time, introduces the results and the
  combinatorics that we need later.

  \subsection{Cyclotomic Hecke algebras and Ariki-Koike
  algebras}\label{S:CycHeckeAlgebras}
  Hecke algebras of the complex reflections groups
  $G_{\ell,n}=\Z/\ell\Z\wr\Sym_n$ of type $G(\ell,1,n)$ were introduced by
  Ariki-Koike~\cite{AK}, motivated by the Iwahori-Hecke algebras of Coxeter
  groups~\cite{IwahoriMatsumoto}. Soon afterwards, Brou\'e and
  Malle~\cite{BM:cyc} defined Hecke algebras for arbitrary complex reflection
  groups. The following refinement of the definition of these algebras unifies
  the treatment of the degenerate and non-degenerate cyclotomic Hecke algebras of type
  $G(\ell,1,n)$.

  Let $\Zcal$ be a commutative domain with one.

  \begin{Definition}[%
    \protect{Hu-Mathas\cite[Definition~2.2]{HuMathas:SeminormalQuiver}}]
    \label{D:HeckeAlgebras}
    Fix integers $n\ge0$ and $\ell\ge1$. The \textbf{cyclotomic Hecke
    algebra of type~$A$}, with \textbf{Hecke parameter} $v\in\Zcal^\times$ and
    \textbf{cyclotomic parameters} $Q_1,\dots,Q_\ell\in\Zcal$, is the unital
    associative $\Zcal$-algebra $\Hcal_n=\Hcal_n(\Zcal,v,Q_1,\dots,Q_\ell)$
    with generators $L_1,\dots,L_n$, $T_1,\dots,T_{n-1}$ and relations
  {\setlength{\abovedisplayskip}{2pt}
   \setlength{\belowdisplayskip}{1pt}
  \begin{align*}
    \textstyle\prod_{l=1}^\ell(L_1-Q_l)&=0,  &
    (T_r+v^{-1})(T_r-v)&=0,                  &
    L_{r+1}&=T_rL_rT_r+T_r,
  \end{align*}
  \begin{align*}
      L_rL_t&=L_tL_r, &
    T_rT_s&=T_sT_r\text{ if }|r-s|>1,\\
    T_sT_{s+1}T_s&=T_{s+1}T_sT_{s+1}, &
    T_rL_t&=L_tT_r\text{ if }t\ne r,r+1,
  \end{align*}
  }
  where $1\le r<n$, $1\le s<n-1$ and $1\le t\le n$.
\end{Definition}

By definition, $\Hcal_n$ is generated by $L_1,T_1,\dots,T_{n-1}$ but we prefer
including $L_2,\dots,L_n$ in the generating set.

Let $\Sym_n$ be the \textbf{symmetric group} on~$n$ letters. For $1\le r<n$ let
$s_r=(r,r+1)$ be the corresponding simple transposition. Then
$\{s_1,\dots,s_{n-1}\}$ is the standard set of Coxeter generators for~$\Sym_n$.
A \textbf{reduced expression} for $w\in\Sym_n$ is a word
$w=s_{r_1}\dots s_{r_k}$ with $k$ minimal and $1\le r_j<n$ for $1\le j\le k$.
If $w=s_{r_1}\dots s_{r_k}$ is reduced then set $T_w=T_{r_1}\dots T_{r_k}$. Then
$T_w$ is independent of the choice of reduced expression by Matsumoto's Monoid
Lemma~\cite{Matsumoto} since the braid relations hold in~$\Hcal_n$; see, for
example, \cite[Theorem~1.8]{M:Ulect}.  Arguing as in~\cite[Theorem~3.3]{AK}, it
follows that~$\Hcal_n$ is free as a~$\Zcal$-module with basis
\begin{equation}\label{E:AKBasis}
  \set{L_1^{a_1}\dots L_n^{a_n} T_w|0\le a_1,\dots,a_n<\ell\text{ and }w\in\Sym_n}.
\end{equation}
Consequently, $\Hcal_n$ is free as a $\Zcal$-module of rank $\ell^n n!$, which is the
order of the complex reflection group $G_{\ell,n}=\Z/\ell\Z\wr\Sym_n$ of type
$G(\ell,1,n)$.

\autoref{D:HeckeAlgebras} is different from Ariki and Koike's~\cite{AK}
definition of the cyclotomic Hecke algebras of type $G(\ell,1,n)$ because we
have changed the commutation relation for $T_r$ and $L_r$.  Ariki and
Koike~\cite{AK} defined their algebra to be the unital associative algebra
generated by $T_0,T_1,\dots,T_{n-1}$ subject to the relations
\begin{gather*}
  \begin{align*}
    \textstyle\prod_{l=1}^\ell(T_0-Q'_l)&=0,  &
    (T_r+v^{-1})(T_r-v)&=0,  \\
    T_0T_1T_0T_1&=T_1T_0T_1T_0 &
    T_sT_{s+1}T_s&=T_{s+1}T_sT_{s+1},
  \end{align*}\\
    T_rT_s=T_sT_r \text{ if }|r-s|>1
\end{gather*}
We have renormalised the quadratic relation for the $T_r$, for $1\le r<n$, so
that $q=v^2$ in the notation of~\cite{AK}. Ariki and Koike then defined $L'_1=T_0$ and
set $L'_{r+1}=T_rL'_rT_r$ for $1\le r<n$. In fact, if $v-v^{-1}$ is invertible
in $\Zcal$ then $\Hcal_n$ is (isomorphic to) the Ariki-Koike algebra with
parameters $Q'_l=1+(v-v^{-1})Q_l$ for $1\le l\le\ell$. To see this set
$L_r'=1+(v-v^{-1})L_r$ in $\Hcal_n$, for $1\le r\le n$. Then
$T_rL'_rT_r=(v-v^{-1})T_rL_rT_r+T_r^2=L'_{r+1}$, which implies our claim.
Therefore, over a field, $\Hcal_n$ is an Ariki-Koike algebra whenever $v^2\ne1$.
On the other hand, if $v^2=1$ then $\Hcal_n$ is a \textit{degenerate} cyclotomic
Hecke algebra~\cite{Klesh:book,AMR}.

We note that the Ariki-Koike algebras with $v^2=1$ include as a special the
group algebras $\Zcal G_{\ell,n}$ of the complex reflection groups $G_{\ell,n}$,
for $n\ge0$.  One consequence of the last paragraph is that $\Zcal G_{\ell,n}$
is \textit{not} a specialization of~$\Hcal_n$. This said, if $F$ is a field such
that $\Hcal_n$ and $FG_{\ell,n}$ are both split semisimple then $\Hcal_n\cong
FG_{\ell,n}$.  On the other hand, the algebras~$\Hcal_n$ always fit into the
\textit{spetses} framework of Brou\'e, Malle and Michel~\cite{BMM:spetses}.

The algebras $\Hcal_n$ with $v^2=1$ are the \textit{degenerate}
cyclotomic Hecke algebras of type $G(\ell,1,n)$ whereas if $v^2\ne1$
then~$\Hcal_n$ is an Ariki-Koike algebra in the sense of~\cite{AK}. Our
definition of~$\Hcal_n$ is more natural in the sense that many features
of the algebras~$\Hcal_n$ have a uniform description in both the
degenerate and non-degenerate cases:
\begin{itemize}
  \item The centre of $\Hcal_n$ is the set of symmetric polynomials in $L_1,\dots,L_n$ (Brundan~\cite{Brundan:degenCentre} in the degenerate case when $v^2=1$ and
  announced when $v^2\ne1$ by Graham and Francis building on~\cite{FrancisGraham}).
  \item The blocks of $\Hcal_n$ are indexed by the same combinatorial data
  (Lyle and Mathas~\cite{LM:AKblocks} when $v^2\ne1$ and  Brundan~\cite{Brundan:degenCentre} when $v^2=1$).
  \item The irreducible $\Hcal_n$-modules are indexed by the crystal graph of the
  integral highest weight module $L(\Lambda)$ for $\Uq$
  (Ariki~\cite{Ariki:can} when $v^2\ne1$ and Brundan and Kleshchev~\cite{BK:DegenAriki} when $v^2=1$).
  \item The algebras $\Hcal_n$ categorify $L(\Lambda)$. Moreover, in
  characteristic zero the projective indecomposable $\Hcal_n$-modules correspond to
  the canonical basis of~$L(\Lambda)$.
  (Ariki~\cite{Ariki:can} when $v^2\ne1$ and Brundan and Kleshchev~\cite{BK:DegenAriki} when $v^2=1$).
  \item The algebra $\Hcal_n$ is isomorphic to a cyclotomic quiver Hecke algebras of
  type~$A$ (Brundan and Kleshchev~\cite{BK:GradedKL}).
\end{itemize}
In contrast, the Ariki-Koike algebras with $v^2=1$ do not share any of these
properties: their center can be larger than the set of symmetric polynomials in
$L_1,\dots,L_n$ (Ariki~\cite{Ariki:can}); if $\ell>1$ then they have only one
block (Lyle and Mathas~\cite{LM:AKblocks}); their irreducible modules are
indexed by a different set (Mathas~\cite{Mathas:AKSimples}); they do not
categorify~$L(\Lambda)$ and no non-trivial grading on these algebras is known.
In this sense, the definition of the Ariki-Koike algebras from~\cite{AK} gives
the wrong algebras when $v^2=1$. \autoref{D:HeckeAlgebras} corrects for this.

Historically, many results for the cyclotomic Hecke algebras $\Hcal_n$
were proved separately in the degenerate ($v^2=1$) and non-degenerate
cases $(v^2\ne1)$.  Using \autoref{D:HeckeAlgebras} it should now be
possible to give uniform proofs of all of these results. In fact, in the
cases that we have checked uniforms arguments can now be given for the
degenerate and non-degenerate cases.

\subsection{Quivers of type $A$ and integral parameters}\label{S:Quivers}
  Rather than work with arbitrary cyclotomic parameters $Q_1,\dots,Q_\ell$, as in
  \autoref{D:HeckeAlgebras}, we now specialize to the \textit{integral case} using
  the Morita equivalence results of Dipper and the author~\cite{DM:Morita} (when
  $v^2\ne1)$ and Brundan and Kleshchev~\cite{BK:HigherSchurWeyl} (when $v^2=1$).
  First, however, we need to introduce quivers and quantum integers.

  Fix $e\in\{1,2,3,4,\dots\}\cup\{\infty\}$ and let $\Gamma_e$ be the
  quiver with vertex set $I_e=\Z/e\Z$ and edges $i\longrightarrow i+1$, for
  $i\in I_e$, where we adopt the convention that $e\Z=\{0\}$ when $e=\infty$. If
  $i,j\in I_e$ and~$i$ and~$j$ are not connected by an edge in~$\Gamma_e$ then
  we write $i\noedge j$. When $e$ is fixed we write $\Gamma=\Gamma_e$ and
  $I=I_e$.  Hence, we are considering either the linear quiver~$\Z$ ($e=\infty$)
  or a cyclic quiver ($e<\infty$):

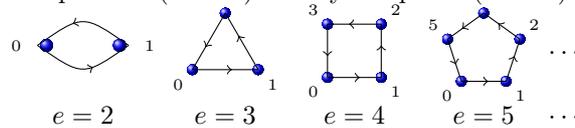
\begin{center}
\vskip-2mm
\begin{tabular}{*5{@{}c}}
  \begin{tikzpicture}[scale=0.6,decoration={curveto, markings,
            mark=at position 0.6 with {\arrow{>}}
    }]
    \useasboundingbox (-1.7,-0.7) rectangle (1.7,0.7);
    \foreach \x in {0, 180} {
      \shade[ball color=blue] (\x:0.8cm) circle(4pt);
    }
    \tikzstyle{every node}=[font=\tiny]
    \draw[postaction={decorate}] (0:1cm) .. controls (90:7mm) .. (180:1cm)
           node[left,xshift=-0.9mm]{$0$};
    \draw[postaction={decorate}] (180:1cm) .. controls (270:7mm) .. (0:1cm)
           node[right,xshift=0.9mm]{$1$};
  \end{tikzpicture}
& 
  \begin{tikzpicture}[scale=0.5,decoration={ markings,
            mark=at position 0.6 with {\arrow{>}}
    }]
    \useasboundingbox (-1.7,-0.7) rectangle (1.7,1.4);
    \foreach \x in {90,210,330} {
      \shade[ball color=blue] (\x:1cm) circle(4pt);
    }
    \tikzstyle{every node}=[font=\tiny]
    \draw[postaction={decorate}]( 90:1cm)--(210:1cm) node[below left]{$0$};
    \draw[postaction={decorate}](210:1cm)--(330:1cm) node[below right]{$1$};
    \draw[postaction={decorate}](330:1cm)--( 90:1cm)
           node[above,yshift=.5mm]{$2$};
  \end{tikzpicture}
& 
  \begin{tikzpicture}[scale=0.5,decoration={ markings,
            mark=at position 0.6 with {\arrow{>}}
    }]
    \useasboundingbox (-1.7,-0.7) rectangle (1.7,0.7);
    \foreach \x in {45, 135, 225, 315} {
      \shade[ball color=blue] (\x:1cm) circle(4pt);
    }
    \tikzstyle{every node}=[font=\tiny]
    \draw[postaction={decorate}] (135:1cm) -- (225:1cm) node[below left]{$0$};
    \draw[postaction={decorate}] (225:1cm) -- (315:1cm) node[below right]{$1$};
    \draw[postaction={decorate}] (315:1cm) -- (45:1cm) node[above right]{$2$};
    \draw[postaction={decorate}] (45:1cm) -- (135:1cm) node[above left]{$3$};
  \end{tikzpicture}
& 
  \begin{tikzpicture}[scale=0.5,decoration={ markings,
            mark=at position 0.6 with {\arrow{>}}
    }]
    \useasboundingbox (-1.7,-0.7) rectangle (1.7,0.7);
    \foreach \x in {18,90,162,234,306} {
      \shade[ball color=blue] (\x:1cm) circle(4pt);
    }
    \tikzstyle{every node}=[font=\tiny]
    \draw[postaction={decorate}] (162:1cm) -- (234:1cm)node[below left]{$0$};
    \draw[postaction={decorate}] (234:1cm) -- (306:1cm) node[below right]{$1$};
    \draw[postaction={decorate}] (306:1cm) -- (18:1cm) node[above right]{$2$};
    \draw[postaction={decorate}] (18:1cm) -- (90:1cm)
           node[above,yshift=.5mm]{$4$};
    \draw[postaction={decorate}] (90:1cm) -- (162:1cm) node[above left]{$5$};
  \end{tikzpicture}
&\raisebox{3mm}{$\dots$}
\\[2mm]
  $e=2$&$e=3$&$e=4$&$e=5$&$\cdots$
\end{tabular}
\end{center}

  In the literature the case $e=\infty$ is often written as $e=0$, however, we
  prefer $e=\infty$ because then $e=|I_e|$. There are also several results that
  hold when $e>n$ --- using the ``$e=0$ convention'' this condition must be
  written as~$e>n$ or~$e=0$. We write $e\ge2$ to mean
  $e\in\{2,3,4,5,\dots\}\cup\{\infty\}$.

  To the quiver $\Gamma_e$ we attach the symmetric Cartan matrix
  $(c_{ij})_{i,j\in I}$, where
  $$c_{ij}=\begin{cases} 2,&\text{if } i=j,\\
    -1,&\text{if $i\rightarrow j$ or $i\leftarrow j$},\\
    -2,&\text{if }i\leftrightarrows j,\\
    0,&\text{otherwise},
  \end{cases}$$
  Following\cite[Chapter~1]{Kac}, let $\widehat{\mathfrak{sl}}_e$ be the
  Kac-Moody algebra of~$\Gamma_e$~\cite{Kac} with simple roots
  $\set{\alpha_i|i\in I}$, fundamental weights $\set{\Lambda_i|i\in I}$,
  positive weight lattice $P^+=\bigoplus_{i\in I}\N\Lambda_i$ and
  positive root lattice $Q^+=\bigoplus_{i\in I}\N\alpha_i$. Let
  $(\cdot,\cdot)$ be the usual invariant form associated with this data,
  normalised so that $(\alpha_i,\alpha_j)=c_{ij}$ and
  $(\Lambda_i,\alpha_j)=\delta_{ij}$, for $i,j\in I$.

  Fix a sequence $\charge=(\kappa_1,\dots,\kappa_\ell)\in\Z^\ell$, the
  \textbf{multicharge}, and define
  $\Lambda=\Lambda(\charge) =
      \Lambda_{\overline{\kappa}_1}+\dots+\Lambda_{\overline{\kappa}_\ell},$
  where $\overline{a}=a+e\Z\in I$ for $a\in\Z$. Then $\Lambda\in P^+$ is
  dominant weight of \textbf{level}~$\ell$.  The integral cyclotomic Hecke
  algebras defined below depend only on~$\Lambda$, however, the bases and our
  combinatorics often depends upon the choice of multicharge~$\charge$.

  Recall that $\Zcal$ is an integral domain. For $t\in\Zcal^\times$ and $k\in\Z$
  define the \textbf{$t$-quantum integer} $[k]_t$ by
  $$[k]_t=\begin{cases}
    t+t^3+\dots+t^{2k-1},&\text{if }k\ge0,\\
    -(t^{-1}+t^{-3}+\dots+t^{2k+1}),&\text{if }k<0.
  \end{cases}$$
  When $t$ is understood we simply write $[k]=[k]_t$. Unpacking the definition,
  if $t^2\ne1$ then $[k]=(t^{2k}-1)/(t-t^{-1})$ whereas $[k]=\pm k$ if
  $t=\pm1$.

  The \textbf{quantum characteristic} of $v$ is the smallest element of
  $e\in\{2,3,4,5,\dots\}\cup\{\infty\}$ such that $[e]_v=0$, where we set
  $e=\infty$ if $[k]_v\ne0$ for all $k>0$.

  \begin{Definition}\label{D:integral}
    Suppose that $\Lambda=\Lambda(\charge)\in P^+$, for $\charge\in\Z^\ell$, and
    that $v\in\Zcal$ has quantum characteristic~$e$.
    The \textbf{integral} cyclotomic Hecke algebra of type $A$ of weight
    $\Lambda$ is the cyclotomic Hecke algebra
    $\H=\Hcal_n(\Zcal,v,Q_1,\dots,Q_r)$ with Hecke parameter~$v$ and cyclotomic
    parameters $Q_r=[\kappa_r]_v$, for $1\le r\le\ell$.
  \end{Definition}

  When $v^2\ne1$ the parameter $Q_r$ corresponds to the Ariki-Koike
  parameters $Q_r'=v^{2\kappa_r}$, for $1\le r\le\ell$, where we use the
  notation of \autoref{S:CycHeckeAlgebras}.

  As observed in~\cite[\S2.2]{HuMathas:SeminormalQuiver}, translating the Morita
  equivalence theorems of~\cite[Theorem~1.1]{DM:Morita}
  and~\cite[Theorem~5.19]{BK:HigherSchurWeyl} into the current setting explains the
  significance of the integral cyclotomic Hecke algebras.

  \begin{Theorem}[\protect{Dipper-Mathas~\cite{DM:Morita}, %
         Brundan-Kleshchev~\cite{BK:HigherSchurWeyl} }]\label{T:Morita}
    Every cyclotomic quiver Hecke algebra $\Hcal_n$ is Morita equivalent to a
    direct sum of tensor products of integral cyclotomic Hecke algebras.
  \end{Theorem}

  Brundan and Kleshchev treated the degenerate case when $v^2=1$ using very
  different arguments than those in~\cite{DM:Morita}. With the benefit of
  \autoref{D:HeckeAlgebras} the argument of~\cite{DM:Morita} now applies
  uniformly to both the degenerate and non-degenerate cases. The Morita
  equivalences in~\cite{DM:Morita,BK:HigherSchurWeyl} are described explicitly,
  with the equivalence being determined by orbits of the cyclotomic parameters.
  See~\cite{DM:Morita,BK:HigherSchurWeyl} for more details.

  In view of \autoref{T:Morita}, it is enough to consider the integral
  cyclotomic Hecke algebras~$\H$ where $v\in\Zcal^\times$ has quantum
  characteristic~$e$ and $\Lambda\in P^+$. This said, for most of
  \autoref{Sec:CycHeckeAlgebras} we consider the general case of a not
  necessarily integral cyclotomic Hecke algebra because we will need this
  generality in~\autoref{S:OKLR}.

  \subsection{Cellular algebras}\label{S:CellularAlgebras} For convenience we
  recall Graham and Lehrer's cellular algebra framework~\cite{GL}. This will
  allow us to define Specht modules for $\Hcal_n$ as cell modules.
  Significantly, the cellular algebra machinery endows the Specht modules with
  an associative bilinear form. Here is the definition.

\begin{Definition}[Graham and Lehrer~\cite{GL}]
  \label{graded cellular def}\label{D:cellular}
  Suppose that $A$ is a $\Zcal$-algebra that is $\Zcal$-free and of finite rank
  as a $\Zcal$-module. A \textbf{cell datum} for $A$ is an ordered triple
  $(\Pcal,T,C)$, where $(\Pcal,\gdom)$ is the \textbf{weight poset}, $T(\lambda)$
  is a finite set for $\lambda\in\Pcal$, and
  $$C\map{\coprod_{\lambda\in\Pcal}T(\lambda)\times T(\lambda)}A;
     (\s,\t)\mapsto c_{\s\t},
  $$
  is an injective map of sets such that:
  \begin{enumerate}
    \item[(GC$_1$)]
    $\set{c_{\s\t}|\s,\t\in T(\lambda) \text{ for } \lambda\in\Pcal}$ is a
      $\Zcal$-basis of $A$.
    \item[(GC$_2$)] If $\s,\t\in T(\lambda)$, for some $\lambda\in\Pcal$, and
    $a\in A$ then there exist scalars $r_{\t\v}(a)$, which do not depend on
    $\s$, such that
      $$c_{\s\t} a=\sum_{\v\in T(\lambda)}r_{\t\v}(a)c_{\s\v}\pmod
      {A^{\gdom\lambda}},$$
      where $A^{\gdom\lambda}=\gen{c_{\a\b}|\mu\gdom\lambda\text{ and }\a,\b\in T(\mu)}_\Zcal$.
    \item[(GC$_3$)] The $\Zcal$-linear map $*\map AA$ determined by
      $c_{\s\t}^*=c_{\t\s}$, for all $\lambda\in\Pcal$ and
      all $\s,\t\in T(\lambda)$, is an anti-isomorphism of $A$.
  \end{enumerate}
  A \textbf{cellular algebra} is an algebra that has a cell datum. If $A$ is a
  cellular algebra with cell datum $(\Pcal,T,C)$ then the basis
  $\set{c_{\s\t}|\lambda\in\Pcal\text{ and } \s,\t\in T(\lambda)}$ is a
  \textbf{cellular basis} of~$A$ with cellular algebra
  anti-isomorphism~$*$.
\end{Definition}

K\"onig and Xi~\cite{KX:CellAlg} have given an equivalent definition of
cellular algebras that does not depend upon a choice of basis. Goodman
and Graber~\cite{GoodmanGraber:JonesBasic} have shown that (GC$_3$) can
be relaxed to the requirement that $(c_{\s\t})^*\equiv
c_{\t\s}\pmod{A^{\gdom\lambda}}$ for some anti-isomorphism~$*$ of~$A$.

The prototypical example of a cellular algebra is a matrix algebra with its
basis of matrix units, which we call a \textit{Wedderburn basis}. As any split
semisimple algebra is isomorphic to a direct sum of matrix algebras it follows
that every split semisimple algebra is cellular.  The cellular algebra framework
is, however, most useful in studying non-semisimple algebras that are not
isomorphic to a direct sum of matrix rings. In general, a cellular basis can be
thought of as an approximation, or weakening, of a basis of matrix units. (This
idea is made more explicit in~\cite{M:seminormal}.)

The cellular basis axioms determines a filtration of the cellular algebra, via
the ideals $A^{\gdom\lambda}$. As we will see, this leads to a quick
construction of its irreducible representations.

For $\lambda\in\Pcal$, let
$A^{\gedom\lambda}=\gen{c_{\a\b}|\mu\gedom\lambda\text{ and }\a,\b\in T(\mu)}_\Zcal$.
Then it follows from \autoref{D:cellular} that $A^{\gedom\blam}$ is a two-sided ideal of~$A$.

Fix $\lambda\in\Pcal$. The \textbf{cell module} $\UnC^\lambda$ is the (right)
$A$-module with basis $\set{c_\t|\t\in T(\lambda)}$ and where $a\in A$ acts
on~$\UnC^\lambda$ by:
$$c_\t a=\sum_{\v\in T(\lambda)}r_{\t\v}(a)c_\v,\qquad\text{for }\t\in T(\lambda),$$
where the scalars $r_{\t\v}(a)\in\Zcal$ are those appearing in (GC$_2$). It follows
immediately from \autoref{D:cellular} that $\UnC^\lambda$ is an $A$-module.
Indeed, if $\s\in T(\lambda)$ then $\UnC^\lambda$ is isomorphic to the submodule
$(c_{\s\t}+A^{\gdom\lambda})A$ of $A/A^\lambda$ via the map
$c_\t\mapsto c_{\s\t}+A^\lambda$, for $\t\in T(\lambda)$.
The cell module $\UnC^\lambda$ comes with a symmetric
bilinear form $\<\ ,\ \>_\lambda$ that is uniquely determined by
\begin{equation}\label{E:CellularBilinearForm}
  \<c_\t,c_\v\>_\lambda c_{\a\b}\equiv c_{\a\t}c_{\v\b}\pmod{A^{\gdom\lambda}},
\end{equation}
for $\a,\b,\t,\v\in T(\lambda)$. By (GC$_2$) of \autoref{D:cellular},  the inner
product $\<c_\t,c_\v\>_\lambda$ depends only on $\t$ and $\v$, and not on the
choices of $\a$ and $\b$. In addition, $\<xa,y\>_\lambda=\<x,ya^*\>_\lambda$, for
all $x,y\in \UnC^\lambda$ and $a\in A$. Therefore,
\begin{equation}\label{E:radical}
  \rad\UnC^\lambda=\set{x\in\UnC^\lambda|
            \<x,y\>_\lambda=0\text{ for all }y\in \UnC^\lambda}
\end{equation}
is an $A$-submodule of~$\UnC^\lambda$. Set
$\UnD^\lambda=\UnC^\lambda/\rad \UnC^\lambda$. Then $\UnD^\lambda$ is an
$A$-module.

The following theorem summarizes some of the main properties of a cellular
algebra. The proof is surprisingly easy given the strength of the result. In
applications the main difficulty is in showing that a given algebra is cellular.

If $M$ is an $A$-module and $D$ is an irreducible $A$-module, let
$[M:D]$ be the decomposition multiplicity of~$D$ in~$M$.

\begin{Theorem}[Graham and Lehrer~\cite{GL}]\label{T:CellularSimples}
  Suppose that $\Zcal=F$ is a field. Then:
  \begin{enumerate}
    \item If $\mu\in\Pcal$ then $\UnD^\mu$ is either zero or
    absolutely irreducible.
    \item Let $\mathcal{K}=\set{\mu\in\Pcal|\UnD^\mu\ne0}$. Then
    $\set{\UnD^\mu|\mu\in\mathcal{K}}$ is a complete set of pairwise non-isomorphic
    irreducible $A$-modules.
    \item If $\lambda\in\Pcal$ and $\mu\in\mathcal{K}$ then
    $[\UnC^\lambda{:}\UnD^\mu]\ne0$
    only if $\lambda\gedom\mu$. Moreover, $[\UnC^\mu{:}\UnD^\mu]=1$.
\end{enumerate}
\end{Theorem}

If $\mu\in\mathcal{K}$ let $\UnP^\mu$ be the projective cover
of~$\UnD^\mu$. It follows from \autoref{D:cellular} that $\UnP^\mu$ has
a filtration in which the quotients are cell modules such that
$\UnC^\lambda$ appears with multiplicity $[\UnC^\lambda{:}\UnD^\mu]$.
Consequently, an analogue of Brauer-Humphreys reciprocity holds for~$A$.
In particular, the Cartan matrix of~$A$ is symmetric.

  \subsection{Multipartitions and tableaux}\label{S:Tableaux} A
  \textbf{partition} of~$m$ is a weakly decreasing sequence
  $\lambda=(\lambda_1,\lambda_2,\dots)$ of non-negative integers such that
  $|\lambda|=\lambda_1+\lambda_2+\dots=m$. An ($\ell$-)\textbf{multipartition}
  of~$n$ is an $\ell$-tuple $\blam=(\lambda^{(1)},\dots,\lambda^{(\ell)})$ of
  partitions such that $|\lambda^{(1)}|+\dots+|\lambda^{(\ell)}|=n$. We identify
  the multipartition $\blam$ with its \textbf{diagram}, which is the set of
  \textbf{nodes} $\diag\blam=\set{(l,r,c)|1\le c\le \lambda^{(l)}_r\text{ for
  }1\le l\le\ell}$.  In this way, we think of $\blam$ as an ordered $\ell$-tuple
  of arrays of boxes in the plane and we talk of the \textbf{components}
  of~$\blam$. Similarly, by the \textbf{rows} and \textbf{columns} of~$\blam$ we
  will mean the rows and columns in each component. For example, if
  $\blam=(3,1^2|2,1|3,2)$ then
  $$\blam=\diag\blam=\TriDiagram(3,1,1|2,1|3,2).$$

  A node $A$ is an \textbf{addable node} of~$\blam$ if $A\notin\blam$ and
  $\blam\cup\{A\}$ is the (diagram of) a multipartition of~$n+1$.  Similarly, a
  node $B$ is a \textbf{removable node} of~$\blam$ if $B\in\blam$ and
  $\blam\setminus\{B\}$ is a multipartition of~$n-1$. If $A$ is an addable node
  of $\blam$ let $\blam+A$ be the multipartition $\blam\cup\{A\}$ and, similarly,
  if $B$ is a removable node let $\blam-B=\blam\setminus\{B\}$. Order the nodes
  lexicographically by~$\le$.

  The set of multipartitions of~$n$ becomes a poset under
  \text{dominance} where~$\blam$ \textbf{dominates}~$\bmu$,
  written as $\blam\gedom\bmu$, if
  $$\sum_{k=1}^{l-1}|\lambda^{(k)}|+\sum_{j=1}^i\lambda^{(l)}_j
       \ge\sum_{k=1}^{l-1}|\mu^{(k)}|+\sum_{j=1}^i\mu^{(l)}_j,$$
  for $1\le l\le\ell$ and $i\ge1$. If $\blam\gedom\bmu$ and $\blam\ne\bmu$
  then write $\blam\gdom\bmu$. Let $\Parts=\Parts[\ell,n]$ be the set of
  multipartitions of~$n$. We consider $\Parts$ as a poset ordered by dominance.

  Fix $\blam\in\Parts$. A \textbf{$\blam$-tableau} is a bijective map
  $\t\map{\diag\blam}\{1,2,\dots,n\}$, which we identify with a labelling of
  (the diagram of) $\blam$ by $\{1,2,\dots,n\}$. For example,
  $$\Tritab({1,2,3},{4},{5}|{6,7},{8}|{9,10,11},{12,13})
    \hspace*{2mm}\text{and}\hspace*{2mm}
    \Tritab({9,12,13},{10},{11}|{6,8},{7}|{1,3,5},{2,4})$$
  are both $\blam$-tableaux when $\blam=(3,1^2|2,1|3,2)$.

  A $\blam$-tableau is \textbf{standard} if its entries increase along rows and
  down columns in each component. For example, the two tableaux above are
  standard. Let $\Std(\blam)$ be the set of standard $\blam$-tableaux. If
  $\Pcal$ is any set of multipartitions let
  $\Std(\Pcal)=\bigcup_{\blam\in\Pcal}\Std(\blam)$. Similarly set
  $\Std^2(\Pcal)=\set{(\s,\t)|\s,\t\in\Std(\blam) \text{ for }\blam\in\Pcal}$.

  If $\t$ is a $\blam$-tableau set $\Shape(\t)=\blam$ and let $\t_{\downarrow m}$
  be the subtableau of~$\t$ that contains the numbers $\{1,2,\dots,m\}$. If~$\t$ is
  a standard $\blam$-tableau then $\Shape(\t_{\downarrow m})$ is a multipartition
  for all $m\ge0$. We extend the dominance ordering to $\Std(\Parts)$, the set
  of all standard tableaux, by defining $\s\gedom\t$ if
  $\Shape(\s_{\downarrow m})\gedom\Shape(\t_{\downarrow m}),$
  for $1\le m\le n$. As before, write $\s\gdom\t$ if $\s\gedom\t$ and
  $\s\ne\t$. Finally, define the \textbf{strong dominance ordering} on
  $\Std^2(\Parts)$ by $(\s,\t)\Gedom(\u,\v)$ if $\s\gedom\u$ and $\t\gedom\v$.
  Similarly, $(\s,\t)\Gdom(\u,\v)$ if $(\s,\t)\Gedom(\u,\v)$ and
  $(\s,\t)\ne(\u,\v)$

  It is easy to see that there are unique standard $\blam$-tableaux $\tlam$ and
  $\tllam$ such that $\tlam\gedom\t\gedom\tllam$, for all $\t\in\Std(\blam)$.
  The tableau $\tlam$ has the numbers $1,2,\dots,n$ entered in order from left
  to right along the rows of $\t^{\lambda^{(1)}}$, and then
  $\t^{\lambda^{(2)}},\dots,\t^{\lambda^{(\ell)}}$. Similarly,~$\tllam$ is
  the tableau with the numbers $1,\dots,n$ entered in order down the columns of
  $\t^{\lambda^{(\ell)}},\dots,\t^{\lambda^{(2)}},\t^{\lambda^{(1)}}$. If
  $\blam=(3,1^2|2,1|3,2)$ then the two $\blam$-tableaux displayed above are
  $\tlam$ and $\tllam$, respectively.

  Given a standard $\blam$-tableau $\t$ define permutations
  $d(\t),d'(\t)\in\Sym_n$ by $\tlam d(\t)=\t=\t_\blam d'(\t)$. Then
  $d(\t)d'(\t)^{-1}=d(\t_{\blam})$ with
  $\ell(d(\t))+\ell(d'(\t))=\ell(d(\t_\blam))$, for all $\t\in\Std(\blam)$. Let
  $\le$ be the Bruhat order on~$\Sym_n$ with the convention that $1\le w$ for
  all $w\in\Sym_n$. Independently, Ehresmann and James~\cite{James} showed that if
  $\s,\t\in\Std(\blam)$ then $\s\gedom\t$ if and only if $d(\s)\le d(\t)$ and if
  and only if $d'(\t)\le d'(\s)$. A proof can be found, for example,
  in~\cite[Theorem~3.8]{M:Ulect}.

  Finally, we will need to know how to conjugate multipartitions and tableaux.
  The conjugate of a partition $\lambda$ is the partition
  $\lambda'=(\lambda'_1,\lambda'_2,\dots)$ where
  $\lambda'_r=\#\set{s\ge1|\lambda_s\ge r}$. That is, we swap the rows and columns
  of~$\lambda$.  The \textbf{conjugate} of a multipartition
  $\blam=(\lambda^{(1)}|\dots|\lambda^{(\ell)})$ is the multipartition
  $\blam'=(\lambda^{(\ell)}{}'|\dots|\lambda^{(1)}{}')$. Similarly, the conjugate
  of a $\blam$-tableau $\t=(\t^{(1)}|\dots|\t^{(\ell)})$ is the
  $\blam'$-tableau $\t'=(\t^{(\ell)}{}'|\dots|\t^{(1)}{}')$ where $\t^{(k)}{}'$
  is the tableau obtained by swapping the rows and columns of $\t^{(k)}$, for
  $1\le k\le\ell$. Then $\blam\gedom\bmu$ if and only if
  $\bmu'\gedom\blam'$, and $\s\gedom\t$ if and only if $\t'\gedom\s'$.

  \subsection{The Murphy basis of $\H$}\label{S:Murphy} Graham and
  Lehrer~\cite{GL} showed that the cyclotomic Hecke algebras (when $v^2\ne1$) are
  cellular algebras. In this section we recall another cellular basis for these
  algebras that was constructed in~\cite{DJM:cyc} when $v^2\ne1$ and
  in~\cite{AMR} when $v^2=1$. When $\ell=1$ these results are due to
  Murphy~\cite{Murphy:basis}.

  First observe that \autoref{D:HeckeAlgebras} implies that there is a unique
  anti-isomorphism~$*$ on $\Hcal_n$ that fixes each of the generators
  $T_1,\dots,T_{n-1},L_1,\dots,L_n$ of $\Hcal_n$. It is easy to see that
  $T_w^*=T_{w^{-1}}$, for $w\in\Sym_n$

  Fix a multipartition $\blam\in\Parts$. Following
  \cite[Definition~3.14]{DJM:cyc} and~\cite[\Sect6]{AMR}, if $\s,\t\in\Std(\blam)$
  define $m_{\s\t}=T_{d(\s)^{-1}}m_\blam T_{d(\t)}$, where
  $m_\blam= u_\blam x_\blam$,
  $$u_\blam=\prod_{1\le l<\ell}\prod_{r=1}^{|\lambda^{(1)}|+\dots+|\lambda^{(l)}|}
  \frac1{Q'_{l+1}}(L_r-[\kappa_{l+1}])
  \quad\text{and}\quad x_\blam=\sum_{w\in\Sym_\blam}v^{\ell(w)}T_w,$$
  where $Q_l'=1+(v-v^{-1})Q_l$ as in \autoref{S:CycHeckeAlgebras}. The
  renormalization of~$u_\blam$ by $1/Q'_{l+1}$ is not strictly necessary. When
  $Q'_{l+1}=0$ this factor can be omitted from the definition of~$u_\blam$, at
  the expense of some aesthetics in some of the formulas that follow. In the
  integral case, which is what we care most about, this problem does not arise
  because $Q'_l=v^{\kappa_l}\ne0$ since $Q_l=[\kappa_l]$, for $1\le l\le\ell$.

  Using the relations in $\H$ it is not hard to show that $u_\blam$ and $x_\blam$
  commute. Consequently, $m_{\s\t}^*=m_{\t\s}$, for all $(\s,\t)\in\Std^2(\Parts)$.

  \begin{Theorem}[\protect{%
      \cite[Theorem~3.26]{DJM:cyc} and \cite[Theorem~6.3]{AMR}}]
    \label{T:MurphyBasis}
    The cyclotomic Hecke algebra $\H$ is free as a $\Zcal$-module with cellular basis
    $\set{m_{\s\t}|\s, \t\in\Std(\blam)\text{ for }\blam\in\Parts}$
    with respect to the poset $(\Parts,\gedom)$.
  \end{Theorem}

  Consequently, $\H$ is a cellular algebra so all of the theory in
  \autoref{S:CellularAlgebras} applies. In particular, for each
  $\blam\in\Parts$ there exists a Specht module $\UnS^\blam$ with basis
  $\set{m_\t|\t\in\Std(\blam)}$. Concretely, we could take
  $m_\t=m_{\tlam\t}+\Hlam$, for $\t\in\Std(\blam)$.

  Let $\UnD^\blam=\UnS^\blam/\rad\UnS^\blam$ be the quotient of $\UnS^\blam$ by
  the radical of its bilinear form. Set
  $\Klesh=\set{\bmu\in\Parts|\UnD^\bmu\ne0}$. By
  \autoref{T:CellularSimples} we obtain:

  \begin{Corollary}[\cite{GL,DJM:cyc,AMR}]\label{C:UngradedSimples}
    Suppose that $\Zcal=F$ is a field. Then $\set{\UnD^\bmu|\bmu\in\Klesh}$ is a
    complete set of pairwise non-isomorphic irreducible $\H$-modules.
  \end{Corollary}

  The set of multipartitions $\Klesh$ has been determined by
  Ariki~\cite{Ariki:class}; see also
  \cite{BK:GradedDecomp,ArikiJaconLecouvey:ModularBranching}. We describe
  and recover his classification of the irreducible $\H$-modules in
  \autoref{C:SimpleClass} below. When $\ell\ge3$ the only known
  descriptions of $\Klesh$ are recursive.
  See~\cite{DJ:reps,ArikiKriemanShunsuke} for $\ell\le2$
  and~\cite{Mathas:AKSimples} when $e=2$.

  \subsection{Semisimple cyclotomic Hecke algebras of type~$A$}
  \label{S:Semisimple}

  We now explicitly describe the semisimple representation theory of~$\H$ using
  the \textit{seminormal coefficient systems} introduced
  in~\cite{HuMathas:SeminormalQuiver}. As we are ultimately interested in the
  cyclotomic quiver Hecke algebras, which are intrinsically non-semisimple
  algebras, it is a little surprising that we are interested in these results.
  We will see, however, that the semisimple representation theory of~$\H$ and
  the KLR grading are closely intertwined.

  The \textbf{Gelfand-Zetlin subalgebra} of $\Hcal_n$ is the subalgebra
  $\L=\L(\Zcal)=\<L_1,L_2,\dots,L_n\>$. We believe that understanding this
  subalgebra is crucial to understanding the representation theory of~$\Hcal_n$.
  To explain how~$\L$ acts on~$\H$ define two \textbf{content} functions
  for $\t\in\Std(\Parts)$ by
  \begin{equation}\label{E:content}
  c^\Zcal_r(\t)=v^{2(c-b)}Q_l+[c-b]_v\in\Zcal    \qquad\text{and}\qquad
             c^\Z_r(\t)=\kappa_l+c-b\in\Z,
  \end{equation}
  where $\t(l,b,c)=r$ and $1\le r\le n$. In the special case of the integral
  parameters, where $Q_l=[\kappa_l]_v$ for $1\le l\le\ell$, the reader can check
  that $c^\Zcal_r(\t)=[c^\Z_r(\t)]_v$, for $1\le r\le n$.

  The next result is well-known and extremely useful.

  \begin{Lemma}[\protect{James-Mathas~\cite[Proposition~3.7]{JM:cyc-Schaper}}]
    \label{L:JucysMurphyAction}
    Suppose that $1\le r\le n$ and that $\s,\t\in\Std(\blam)$, for
    $\blam\in\Parts$. Then
    $$m_{\s\t}L_r \equiv c^\Zcal_r(\t)m_{\s\t}
    +\sum_{\substack{\v\gdom\t\\\v\in\Std(\blam)}}a_\v m_{\s\v}\pmod\Hlam,$$
    for some $a_\v\in\Zcal$.
  \end{Lemma}

  \begin{proof}Let $(l,b,c)=\t^{-1}(r)$. Using our notation,
    \cite[Proposition~3.7]{JM:cyc-Schaper} says that
    $m_{\s\t}L'_r=Q'_lv^{2(c-b)}m_{\s\t}$ plus linear combination of more dominant
    terms, where $Q_l'=1+(v-v^{-1})Q_l$. As $L_r=1+(v-v^{-1})L'_r$ this easily
    implies the result when $v^2\ne1$. The case when $v^2=1$ now follows by
    specialization --- or, see \cite[Lemma~6.6]{AMR}.
  \end{proof}

  In the integral case,
  $m_{\s\t}L_r\equiv [c^\Z_r(\t)]m_{\s\t}+\sum_{\v\gdom\t}a_\v m_{\s\t}\pmod\Hlam$.
  This agrees with~\cite[Lemma~2.9]{HuMathas:SeminormalQuiver}.

  The Hecke algebra $\Hcal_n$ is \textbf{content separated} if whenever
  $\s,\t\in\Std(\Parts)$ are standard tableaux, not necessarily of the same
  shape, then $\s=\t$ if and only if~$c^\Zcal_r(\s)=c^\Zcal_r(\t)$, for $1\le r\le n$.
  The following is an immediate corollary of
  \autoref{L:JucysMurphyAction} using the theory of JM-elements developed in
   \cite[Theorem~3.7]{M:seminormal}.

  \begin{Corollary}[\protect{\cite[Proposition~3.4]{HuMathas:SeminormalQuiver}}]
    \label{C:LBimod}
    Suppose that $\Zcal=F$ is a field and that $\Hcal_n$ is content separated. Then,
    as an $(\L,\L)$-bimodule,
    $$\Hcal_n=\bigoplus_{(\s,\t)\in\Std^2(\Parts)}H_{\s\t},$$
    where $H_{\s\t}=\set{h\in\Hcal_n|L_rh=c^\Zcal_r(\s)h\text{ and }
            hL_r=c^\Zcal_r(\t)h,\text{ for }1\le r\le n}$.
  \end{Corollary}

  For the rest of~\autoref{S:Semisimple} we assume that $\Hcal_n$ is content
  separated.  \autoref{C:LBimod} motivates the following definition.

  \begin{Definition}[%
     \protect{Hu-Mathas~\cite[Definition~3.7]{HuMathas:SeminormalQuiver}}]
     \label{D:SeminormalBasis}
     Suppose that $\Zcal=K$ is a field. A $*$-seminormal basis of $\Hcal_n$ is
     a basis of the form
     $$\set{f_{\s\t}|0\ne f_{\s\t}\in H_{\s\t} \text{ and } f_{\s\t}^*=f_{\t\s},
                    \text{ for }(\s,\t)\in\Std^2(\Parts)}.$$
  \end{Definition}

  There is a vast literature on seminormal bases. This story started with
  Young's seminormal forms for the symmetric groups~\cite{QSAI} and has now been
  extended to Hecke algebras and many other diagram algebras including the
  Brauer, BMW and partition algebras; see, for
  example,~\cite{Ram:seminormal,Nazarov:brauer,M:tilting}.

  Suppose that $\{f_{\s\t}\}$ is a $*$-seminormal basis and that
  $(\s,\t),(\u,\v)\in\Std^2(\Parts)$. Let
  $\Cont=\set{c^\Zcal_r(\s)|\s\in\Std(\Parts)\text{ for }1\le r\le n}$
  be the set of all possible contents for the tableaux in $\Std(\Parts)$.
  Following Murphy~\cite{M:Nak,M:seminormal},
  for a standard tableau $\s\in\Std(\Parts)$ define
  $$F_\s=\prod_{r=1}^n\prod_{\substack{c\in\Cont\\c\ne c^\Zcal_r(\s)}}
            \frac{L_r-c}{c^\Zcal_r(\s)-c}.$$
  By \autoref{D:SeminormalBasis}, if $(\s,\t),(\u,\v)\in\Std^2(\Parts)$ then
  $\delta_{\s\u}\delta_{\t\v}f_{\s\t}=F_\u f_{\s\t}F_\v$.
  In particular, $F_\s$ is a non-zero element of~$\Hcal_n$. It follows that
  $F_\s$ is a scalar multiple of $f_{\s\s}$, which implies that
  $\set{F_\s|\s\in\Std(\Parts)}$ is a complete set of pairwise orthogonal idempotents
  in~$\Hcal_n$. (In fact, in \cite{M:seminormal} these properties are used to
  establish \autoref{C:LBimod}.) Consequently, there exists a non-zero
  scalar~$\gamma_\s\in F$ such that $F_\s=\frac1{\gamma_s}f_{\s\s}$.
  If $(\s,\t),(\u,\v)\in\Std^2(\Parts)$ then
  \begin{equation}\label{E:SeminormalProduct}
    f_{\s\t}f_{\u\v} = f_{\s\t}F_\t F_\u f_{\u\v}=\delta_{\t\u}\gamma_\t f_{\s\v},
  \end{equation}

  The next definition allows us to classify all seminormal bases and to
  describe how $\H$ acts on them.

  \begin{Definition}[\protect{Hu-Mathas~\cite[\S3]{HuMathas:SeminormalQuiver}}]
    \label{D:SNCS}
    A \textbf{$*$-seminormal coefficient system} is a collection of scalars
    $$\balpha=\set{\alpha_r(\t)|\t\in\Std(\Parts)\text{ and }1\le r\le n}$$
    such that $\alpha_r(\t)=0$ if $\v=\t(r,r+1)$ is not standard, if
    $\v\in\Std(\Parts)$ then
    $$\alpha_r(\v)\alpha_r(\t)
        =\frac{\(1-v^{-1}c^\Zcal_r(\t)+vc^\Zcal_r(\v)\)%
               \(1+vc^\Zcal_r(\t)-v^{-1}c^\Zcal_r(\v)\)}
               {\(c^\Zcal_r(\t)-c^\Zcal_r(\v)\)%
               \(c^\Zcal_r(\v)-c^\Zcal_r(\t)\)},
    $$
    and
    $\alpha_r(\t)\alpha_{r+1}(\t s_r)\alpha_r(\t s_rs_{r+1})
          =\alpha_{r+1}(\t)\alpha_r(\t s_{r+1})\alpha_{r+1}(\t s_{r+1}s_r)$,
     and if~$|r-r'|>1$ then
    $\alpha_r(\t)\alpha_{r'}(\t s_r)=\alpha_{r'}(\t)\alpha_r(\t
    s_{r'})$, for $1\le r,r'<n$.
  \end{Definition}

  As the reader might guess, the conditions on the scalars $\alpha_r(\t)$ in
  \autoref{D:SNCS} correspond to the quadratic relations $(T_r-v)(T_r+v^{-1})=0$
  and the braid relations for $T_1,\dots,T_{n-1}$. The simplest example of a seminormal
  coefficient system is
  $$ \alpha_r(\t)=\frac{\(1-v^{-1}c^\Zcal_{r+1}(\t)+vc^\Zcal_r(\t)\)}%
                       {\(c^\Zcal_{r+1}(\t)-c^\Zcal_r(\t)\)},
  $$
  whenever $1\le r<n$ and $\t, \t(r,r+1)\in\Std(\Parts)$. Another seminormal
  coefficient system is given in \autoref{E:MurphySNCS} below.

  Seminormal coefficient systems arise because they describe the action
  of~$\Hcal_n$ on a seminormal basis. More precisely, we have the following:

  \begin{Theorem}[Hu-Mathas~\cite{HuMathas:SeminormalQuiver}]\label{T:SeminormalBasis}
    Suppose that $\Zcal=K$ is a field and that $\Hcal_n$ is content separated
    and that
    $\set{f_{\s\t}|(\s,\t)\in\Std^2(\Parts)}$ is a seminormal basis
    of~$\Hcal_n$. Then $\{f_{\s\t}\}$ is a cellular basis of~$\Hcal_n$ and there
    exists a unique seminormal coefficient system $\balpha$ such that
    $$ f_{\s\t}T_r = \alpha_r(\t)f_{\s\v}
    +\frac{1+(v-v^{-1})c^\Zcal_{r+1}(\t)}{c^\Zcal_{r+1}(\t)-c^\Zcal_r(\t)}f_{\s\t},$$
    where $\v=\t(r,r+1)$. Moreover, if $\s\in\Std(\blam)$ then
    $F_\s=\frac1{\gamma_\s}f_{\s\s}$ is a primitive idempotent and $\UnS^\blam\cong
    F_\s\Hcal_n$ is irreducible for all $\blam\in\Parts$.
  \end{Theorem}

  \begin{proof}[Sketch of proof]
    By definition, $\{f_{\s\t}\}$ is a basis of $\Hcal_n$ such that
    $f_{\s\t}^*=f_{\t\s}$ for all $(\s,\t)\in\Std^2(\Parts)$. Therefore, it follows
    from \autoref{E:SeminormalProduct} that $\{f_{\s\t}\}$ is a cellular basis
    of~$\Hcal_n$ with cellular automorphism~$*$.

    It is an amusing application of the relations in \autoref{D:HeckeAlgebras} to
    show that there exists a seminormal coefficient system that describes the
    action of~$T_r$ on the seminormal basis. See
    \cite[Lemma~3.13]{HuMathas:SeminormalQuiver} for details. The uniqueness
    of~$\balpha$ is clear.

    We have already observed in \autoref{E:SeminormalProduct} that
    $F_\s=\frac1{\gamma_\s}f_{\s\s}$, for
    $\s\in\Std(\blam)$, so it remains to show
    that $F_\s$ is primitive and that $\UnS^\blam\cong F_\s\Hcal_n$. By what we have
    already shown, $F_\s\Hcal_n$ is contained in the span of
    $\set{f_{\s\t}|\t\in\Std(\blam)}$. On the other hand, if
    $f=\sum_\t r_\t f_{\s\t}\in F_\s\Hcal_n$ and $r_\v\ne0$ then
    $r_\v f_{\s\v}=fF_\v\in F_\s\Hcal_n$.  It follows that $F_\s\Hcal_n=\sum_\t
    Kf_{\s\t}$, as a vector space. Consequently, $F_\s\Hcal_n$ is irreducible
    and $F_\s$ is a primitive idempotent in~$\Hcal_n$. Finally, $\UnS^\blam\cong
    F_\s\Hcal_n$ by \autoref{L:JucysMurphyAction} since $\Hcal_n$ is content
    separated.
  \end{proof}

  \begin{Corollary}[\protect{\cite[Corollary~3.17]{HuMathas:SeminormalQuiver}}]
  Suppose that $\balpha$ is a seminormal coefficient system and that
  $\s\gdom\t=\s(r,r+1)$, for tableaux $\s,\t\in\Std(\Parts)$ and where
  $1\le r<n$. Then $\alpha_r(\s)\gamma_\t=\alpha_r(\t)\gamma_\s$.
  \end{Corollary}

  Consequently, if the seminormal coefficient system $\balpha$ is known then
  fixing $\gamma_\t$, for some $\t\in\Std(\blam)$, determines $\gamma_\s$ for
  all $\s\in\Std(\blam)$. Conversely, these scalars, together with~$\balpha$,
  determines the seminormal basis.

  \begin{Corollary}[\protect{Classifying seminormal bases
       \cite[Theorem~3.14]{HuMathas:SeminormalQuiver}}]\label{C:SeminormalClass}
    There is a one-to-one correspondence between the $*$-seminormal
    bases of~$\Hcal_n$ and the pairs $(\balpha,\bgamma)$ where
    $\balpha=\set{\alpha_r(\s)|1\le r<n\text{ and }\s\in\Std(\Parts)}$ is a
    seminormal coefficient system and
    $\bgamma=\set{\gamma_{\tlam}|\blam\in\Parts}$.
  \end{Corollary}

  Finally, the seminormal basis machinery in this section can be used to classify
  the semisimple cyclotomic Hecke algebras~$\Hcal_n$, thus re-proving Ariki's
  semisimplicity criterion \cite{Ariki:ss}, when $v^2\ne1$, and
  \cite[Theorem~6.11]{AMR}, when $v^2=1$.

  \begin{Theorem}[\protect{Ariki~\cite{Ariki:ss} and \cite[Theorem~6.11]{AMR}}]
    \label{T:SSimple}
    Suppose that $F$ is a field. The following are equivalent:
    \begin{enumerate}
      \item $\Hcal_n=\Hcal_n(F,v,Q_1,\dots,Q_\ell)$ is semisimple.
      \item $\Hcal_n$ is content separated.
      \item $\displaystyle[1]_v[2]_v\dots[n]_v\prod_{1\le r<s\le\ell}\prod_{-n<d<n}
    (v^{2d}Q_r+[d]_v-Q_s)\ne0$.
    \end{enumerate}
  \end{Theorem}

  We want to rephrase the semisimplicity criterion of \autoref{T:SSimple} for
  the integral cyclotomic Hecke algebras $\H$, for $\Lambda\in P^+$. For
  each $i\in I$ define the \textbf{$i$-string of length~$n$} to be
  $\alpha_{i,n}=\alpha_i+\alpha_{i+1}+\dots+\alpha_{i+n-1}$.  Then
  $\alpha_{i,n}\in Q^+$.

  \begin{Corollary}\label{C:IntegralSSimple}
    Suppose that $\Lambda\in P^+$ and that $\Zcal=F$ is a field. Then~$\H$ is
    semisimple if and only if $e>n$ and~$(\Lambda,\alpha_{i,n})\le 1$, for
    all $i\in I$.
  \end{Corollary}

  \begin{proof} As $Q_r=[\kappa_r]$, for $1\le r\le\ell$, the statement of
    \autoref{T:SSimple}(c) simplifies because
    $v^{2d}Q_r+[d]_v-Q_s=v^{-2\kappa_s}[d+\kappa_r-\kappa_s]_v$.  Therefore,
    $\H$ is semisimple if and only if
    $$[1]_v[2]_v\dots[n]_v\prod_{1\le r<s\le\ell}\prod_{-n<d<n}
    [d+\kappa_r-\kappa_s]_v\ne0.$$
    On the other hand, $[1]_v[2]_v\dots[n]_v\ne0$ if and only if $e>n$.
    Furthermore, $(\Lambda,\alpha_{i,n})\le1$, for all $i\in I$, if and only
    if $\kappa_r+d\ne\kappa_s$, for $1\le r<s\le\ell$ and all~$-n<d<n$.
    The result follows.
  \end{proof}

  \subsection{Gram determinants and the Jantzen sum formula} For
  future use, we now recall the closed formula for the Gram determinants of the
  Specht modules $\UnS^\blam$ and the connection between these formulas and
  Jantzen filtrations. Throughout this section we assume that $\Hcal_n$ is
  content separated over the field $K=\Zcal$.

  For $\blam\in\Parts$ let $\UnG^\blam=\(\<m_\s,m_\t\>\)_{\s,\t\in\Std(\blam)}$
  be the \textbf{Gram matrix} of the Specht module $\UnS^\blam$, where we fix an
  arbitrary ordering of the rows and columns of~$\UnG^\blam$.

  For $(\s,\t)\in\Std^2(\Parts)$ set $f_{\s\t}=F_\s m_{\s\t}F_\t$.
  By \autoref{L:JucysMurphyAction} and \autoref{E:SeminormalProduct},
  $$f_{\s\t}=m_{\s\t} + \sum_{(\u,\v)\gdom(\s,\t)}r_{\u\v}m_{\u\v},$$
  for some $r_{\u\v}\in K$. By construction, $\{f_{\s\t}\}$ is a seminormal
  basis of $\Hcal_n$. By \cite[Proposition~3.18]{HuMathas:SeminormalQuiver} this
  basis corresponds to the seminormal coefficient system given by
  \begin{equation}\label{E:MurphySNCS}
      \alpha_r(\t)=\begin{cases}
              1,&\text{if }\t\gdom\t(r,r+1),\\
              \frac{\(1-v^{-1}c_r(\t)+vc_r(\v)\)\(1+vc_r(\t)-v^{-1}c_r(\v)\)}
        {\(c_r(\t)-c_r(\v)\)\(c_r(\v)-c_r(\t)\)}, &\text{otherwise,}
      \end{cases}
   \end{equation}
   for $\t\in\Std(\Parts)$ and $1\le r<n$ such that $\t s_r$ is standard.
   The $\gamma$-coefficients $\{\gamma_\t\}$ for this basis are explicitly
   known  by \cite[Corollary~3.29]{JM:cyc-Schaper}. Moreover,
   \begin{equation}\label{E:GammaProduct}
           \det\UnG^\blam=\prod_{\t\in\Std(\blam)}\gamma_\t
   \end{equation}
   By explicitly computing the scalars $\gamma_\t$, and using an intricate
   inductive argument based on the semisimple branching rules for the Specht
   modules, James and the author proved the following:

   \begin{Theorem}[\protect{James-Mathas~\cite[Corollary~3.38]{JM:cyc-Schaper}}]
     \label{T:SchaperDet}
      Suppose that $\Hcal_n$ is content separated. Then there exist explicitly
      known scalars $g_{\blam\bmu}$ and signs $\varepsilon_{\blam\bmu}=\pm1$ such
      that
      $$\det\UnG^\blam = \prod_{\substack{\bmu\in\Parts\\\blam\gdom\bmu}}
             g_{\blam\bmu}^{\varepsilon_{\blam\bmu}\dim\UnS^\bmu}.$$
   \end{Theorem}

   The scalars $g_{\blam\bmu}$ are described combinatorially as the quotient of
   at most two \textit{hook lengths}.
   The sign $\varepsilon_{\blam\bmu}$ is the parity of the sum of the leg
   lengths of these hooks.

  \autoref{T:SchaperDet} gives a very pretty closed formula for the Gram
  determinant~$\UnG^\blam$, generalizing a classical result of James and
  Murphy~\cite{JM:det}. One problem with this formula is that $\det\UnG^\blam$
  is a \textit{polynomial} in $v,v^{-1},Q_1,\dots,Q_\ell$ whereas
  \autoref{T:SchaperDet} computes this determinant as a
  \textit{rational function} in~$v,Q_1,\dots,Q_\ell$. On the other hand, as we
  now recall, \autoref{T:SchaperDet} has an impressive module theoretic
  application in the \textit{Jantzen sum formula}.

  Fix a modular system $(K,\Zcal,F)$, where $\Zcal$ is a discrete valuation ring
  with maximal ideal~$\p$ and such that $\Zcal$ contains
  $v,v^{-1},Q_1,\dots,Q_\ell$, Let $K$ be the field of fractions of~$\Zcal$ and
  let $F=\Zcal/\p$ be the residue field of~$\Zcal$. Let $\Hcal^\Zcal_n$,
  $\Hcal^K_n\cong\Hcal^\Zcal_n\otimes_\Zcal K$ and
  $\Hcal^F_n=\Hcal^\Zcal_n\otimes_\Zcal F$ be the corresponding Hecke algebras.
  Therefore, $\Hcal^F_n$ has Hecke parameter $v+\p$ and cyclotomic parameters
  $Q_l+\p$, for $1\le l\le\ell$.

  Let $\blam\in\Parts$ and let $\UnS^\blam_\Zcal$ and
  $\UnS^\blam_F\cong\UnS^\blam_\Zcal\otimes_\Zcal F$ be the corresponding Specht
  modules for $\Hcal^\Zcal_n$ and $\Hcal^F_n$, respectively. Define a filtration
  of the Specht module~$\UnS^\blam_\Zcal$ by
  $J_k(\UnS^\blam_\Zcal) = \set{x\in\UnS^\blam_\Zcal|\<x,y\>_\blam\in\p^k
                \text{ for all } y\in\UnS^\blam_\Zcal},$
  for $k\ge0$. The \textbf{Jantzen filtration} of $\UnS_F^\blam$ is the filtration
  $$\UnS_F^\blam=J_0(\UnS_F^\blam)\supseteq J_1(\UnS_F^\blam)\supseteq\dots\supseteq
       J_z(\UnS_F^\blam)=0,$$
  where
  $J_k(\UnS_F^\blam)=\(J_k(\UnS_\Zcal^\blam)+\p\UnS^\blam_\Zcal\)/\p\UnS^\blam_\Zcal$
  for $k\ge0$.  (As $\UnS^\blam_F$ is finite dimensional, $J_z(\UnS^\blam_F)=0$ for $z\gg0$.)

  Let $\rep(\Hcal_n)$ be the category of finitely generated $\Hcal_n$-modules
  and let $[\rep(\Hcal_n)]$ be its Grothendieck group.  Let~$[M]$ be the image
  of the $\Hcal_n$-module $M$ in $[\rep(\Hcal_n)]$. Let $\nu_\p$ be the
  $\p$-adic valuation map on~$\Zcal^\times$.

  \begin{Theorem}[\protect{James-Mathas~\cite[Theorem~4.6]{JM:cyc-Schaper}}]
    \label{T:Jantzen}
    Let $(K,\Zcal,F)$ be a modular system and suppose that
    $\blam\in\Parts$. Then, in $[\rep(\Hcal^F_n)]$,
    $$\sum_{k>0}\,[J_k(\UnS^\blam_F)] = \sum_{\blam\gdom\bmu}
              \varepsilon_{\blam\bmu}\nu_\p(g_{\blam\bmu})[\UnS_F^\bmu].$$
  \end{Theorem}

  Intuitively, the proof of \autoref{T:Jantzen} amounts to taking the $\p$-adic
  valuation of the formula in \autoref{T:SchaperDet}. In fact, this is exactly
  how \autoref{T:Jantzen} is proved except that you need the corresponding
  formulas for the Gram determinants of the weight spaces of the Weyl modules of
  the cyclotomic Schur algebras of~\cite{DJM:cyc}. This is enough because the
  dimensions of the weight spaces of a module uniquely determine its image in
  the Grothendieck group of the Schur algebra. The proof given in \cite{JM:cyc-Schaper}
  is stated only for the non-degenerate case~$v^2\ne1$, however, the arguments
  apply equally well for the degenerate case when $v^2=1$.

  The main point that we want to emphasize in this section is that the
  \textit{rational} formula for $\det\UnG^\blam_K$ in \autoref{T:SchaperDet}
  corresponds to writing the left-hand side of the Jantzen sum formula sum as a
  $\Z$-linear combination of Specht modules. Therefore, when the right-hand side
  of the sum formula is written as a linear combination of simple modules some
  of the terms must cancel. We give a cancellation free sum formula in
  \autoref{S:GramDetsII}.

  \autoref{T:Jantzen} is a useful inductive tool because it gives an upper bound
  on the decomposition numbers of~$\UnS^\blam_F$. Let
  $j_{\blam\bmu}= \varepsilon_{\blam\bmu}\nu_\p(g_{\blam\bmu})$, for
  $\blam,\bmu\in\Parts$ and set $d^F_{\blam\bmu}=[\UnS^\blam_F:\UnD^\bmu_F]$.
  Using \autoref{T:Jantzen} to compute the multiplicity of $\UnD^\bmu_F$ in
  $\bigoplus_{k>0}J_k(\UnS^\blam_F)$ yields the following.

  \begin{Corollary}
    Suppose that $\blam,\bmu\in\Parts$. Then
    $0\le d^F_{\blam\bmu}\le
        \Sum_{\substack{\bnu\in\Parts\\\blam\gdom\bnu\gedom\bmu}}
        j_{\blam\bnu}d^F_{\bnu\bmu}.$
  \end{Corollary}

  As a second application, \autoref{T:Jantzen} classifies the irreducible Specht
  modules $\UnS^\blam_F$, for $\blam\in\Klesh$.

  \begin{Corollary}[\protect{James-Mathas~\cite[Theorem~4.7]{JM:cyc-Schaper}}]
    \label{C:IrredSpechts}
    Suppose that $F$ is a field and $\blam\in\Klesh$. Then the Specht
    module $\UnS^\blam_F$ is irreducible if and only if~$j_{\blam\bmu}=0$
    for all $\bmu\gdom\blam$.
  \end{Corollary}

  \subsection{The blocks of $\Hcal^F_n$}\label{S:blocks} The most
  important application of the Jantzen Sum Formula (\autoref{T:Jantzen})
  is to the classification of the blocks of~$\Hcal^F_n$. The algebra
  $\Hcal^F_n$, and in fact any finite-dimensional algebra over a field,
  can be written as a direct sum of indecomposable two-sided ideals:
  $\Hcal^F_n=B_1\oplus\dots\oplus B_d$. The subalgebras $B_1,\dots,B_z$,
  which are the \textbf{blocks} of~$\Hcal^F_n$, are uniquely determined up
  to permutation. Any $\Hcal^F_n$-module $M$ splits into a direct sum of
  block components $M=MB_1\oplus\dots\oplus MB_d$, where we allow some
  of the summands to be zero. The module $M$ \textbf{belongs} to the
  block~$B_r$ if $M=MB_r$.  It is a standard fact that two simple
  modules $\UnD^\blam$ and $\UnD^\bmu$ belong to the same block if and
  only if they are in the same \textbf{linkage class}. That is, there
  exists a sequence of indecomposable modules
  $\underline{M}_0,\dots,\underline{M}_z$ and a sequence of
  multipartitions $\bnu_0=\blam$, $\bnu_1,\dots,\bnu_z=\bmu$ such that
  $[\underline{M}_r:\UnD^{\bnu_r}]\ne0$ and
  $[\underline{M}_r:\UnD^{\bnu_{r+1}}]\ne0$, for $0\le r<z$. In fact, we
  can assume that the $\underline{M}_r$ are Specht modules, even though
  the Specht modules are not necessarily indecomposable.

  We want an explicit combinatorial description of the blocks of~$\Hcal^F_n$.
  Define two equivalence relations $\sim_C$ and $\sim_J$ on $\Parts$ as follows.
  First, $\blam\sim_C\bmu$ if there is an equality of \textit{multisets}
  $\set{c^F_{\tlam}(r)|1\le r\le n}=\set{c^F_{\tmu}(r)|1\le r\le n}$. The
  second relation, \textbf{Jantzen equivalence}, is more involved:
  $\blam\sim_J\bmu$ if there exists a sequence
  $\bnu_0=\blam,\bnu_1,\dots,\bnu_z=\bmu$ of multipartitions in~$\Parts$ such
  that $j_{\bnu_r\bnu_{r+1}}\ne0$ or $j_{\bnu_{r+1}\bnu_r}\ne0$, for $0\le r<z$.

  \begin{samepage}
  \begin{Theorem}[\protect{Lyle-Mathas~\cite{LM:AKblocks},
      Brundan~\cite{Brundan:degenCentre}}]\label{T:Blocks}\leavevmode\newline
      Suppose that $F$ is a field and that $\blam,\bmu\in\Parts$. Then
      the following are equivalent:
      \begin{enumerate}
        \item $\UnD^\blam$ and $\UnD^\bmu$ are in the same $\Hcal^F_n$-block.
        \item $\UnS^\blam$ and $\UnS^\bmu$ are in the same $\Hcal^F_n$-block.
        \item $\blam\sim_J\bmu$.
        \item $\blam\sim_C\bmu$.
    \end{enumerate}
  \end{Theorem}
  \end{samepage}

  Parts~(a) and (b) are equivalent by the general theory of cellular
  algebras~\cite{GL} whereas the equivalence of parts (b) and (c) is a general
  property of Jantzen filtrations from~\cite{LM:AKblocks}. (In fact, part~(c) is
  general property of the standard modules of a quasi-hereditary algebra.) In
  practice, part~(d) is the most useful because it is easy to compute.

  The hard part in proving \autoref{T:Blocks} is in showing that parts (c)
  and (d) are equivalent. The argument is purely combinatorial with work of
  Fayers~\cite{Fayers:AKweight,Fayers:multicore} playing an important role.

  In the integral case, when $\Hcal^F_n=\H$ for some $\Lambda\in P^+$, there is
  a nice reformulation of \autoref{T:Blocks}.
  The \textbf{residue sequence} of a standard tableau~$\t$ is
  $\bi^\t=(i^\t_1,\dots,i^\t_n)\in I^n$ where $i^\t_r=c^\Z_r(\t)+e\Z$. If
  $\t\in\Std(\blam)$, for $\blam\in\Parts$, define
  $$\beta^\blam = \sum_{r=1}^n\alpha_{i^\t_r}=\sum_{r=1}^n\alpha_{i^\blam_r}
                \in Q^+.$$
  By definition, $\beta^\blam\in Q^+$ depends only on~$\blam$, and not
  on the choice of~$\t$. Moreover, $\blam\sim_C\bmu$ if and only if
  $\beta^\blam=\beta^\bmu$. Hence, we have the following:

  \begin{Corollary}\label{C:IntegralBlocks}
      Suppose that $\Lambda\in P^+$ and $\blam,\bmu\in\Parts$. Then $\UnS^\blam$
      and $\UnS^\bmu$ are in the same $\H$-block if and only
      if~$\beta^\blam=\beta^\bmu$.
  \end{Corollary}

\section{Cyclotomic quiver Hecke algebras of type~$A$}

  This section introduces the quiver Hecke algebras, and their cyclotomic
  quotients. We use the relations to reveal some of the properties of these
  algebras. The main aim of this section is to give the reader an appreciation
  of, and some familiarity with, the KLR relations without appealing to any
  general theory.

  \subsection{Graded algebras}\label{S:GradedCellular}
  In this section we quickly review the theory of graded (cellular) algebras.
  For more details the reader is referred to
  \cite{BGS:Koszul,NV:graded,HuMathas:GradedCellular}. Throughout,~$\Zcal$ is a
  commutative integral domain. Unless otherwise stated, all modules and algebras
  will be free and of finite rank as $\Zcal$-modules.

  In this chapter a \textbf{graded module} will always mean a $\Z$-graded
  module. That is, a $\Zcal$-module $M$ that has a decomposition
  $M=\bigoplus_{d\in\Z}M_d$ as a $\Zcal$-module. A \textbf{positively graded}
  module is a graded module $M=\bigoplus_d M_d$ such that $M_d=0$ if $d<0$.

  A \textbf{graded algebra} is a  unital associative
  $\Zcal$-algebra $A=\bigoplus_{d\in\Z}A_d$ that is a $\Z$-graded
  $\Zcal$-module such that $A_dA_e\subseteq A_{d+e}$, for all
  $d,e\in\Z$. It follows that $1\in A_0$ and that $A_0$ is a $\Z$-graded
  subalgebra of $A$.  A graded (right) $A$-module is a graded
  $\Zcal$-module $M$ such that $\underline{M}$ is an
  $\underline{A}$-module and $M_dA_e\subseteq M_{d+e}$, for all
  $d,e\in\Z$. Here $\underline{M}$ is the (ungraded) module, and
  $\underline{A}$ is the (ungraded) algebra, obtained by forgetting the
  $\Z$-gradings on $M$ and $A$ respectively. Graded submodules, graded
  left $A$-modules and so on are all defined in the obvious way.

  Suppose that $M$ is a graded $A$-module. If  $m\in M_d$, for $d\in\Z$, then
  $m$ is \textbf{homogeneous} of \textbf{degree} $d$ and we set $\deg m=d$.
  Every element $m\in M$ can be written uniquely as a linear combination
  $m=\sum_d m_d$ of its \textbf{homogeneous components}, where $\deg m_d=d$ and
  $m_d\in M$.

  A homomorphism of graded $A$-modules $M$ and $N$ is an $\underline{A}$-module
  homomorphism $f\map{\underline{M}}{\underline{N}}$ such that $\deg f(m)=\deg m$,
  whenever $m\in M$ is homogeneous. That is, $f$ is a degree preserving
  $\underline{A}$-module homomorphism.

  Let $\rep(A)$ be the category of finitely generated graded $A$-modules
  together with degree preserving homomorphisms. Similarly, $\proj(A)$ is the
  category of finitely generated projective $A$-modules with degree preserving
  maps. A \textbf{graded functor} between such categories is any functor
  that commutes with the grading shift functor that sends $M$ to
  $M\<1\>$.

  If $M$ is a graded $\Zcal$-module and $s\in\Z$ let $M\<s\>$ be the graded
  $\Zcal$-module obtained by shifting the grading on~$M$ up by~$s$; that is,
  $M\<s\>_d=M_{d-s}$, for $d\in\Z$. If $M\ne0$ then $M\cong M\<s\>$ as
  $A$-modules if and only if $s=0$. In contrast,
  $\underline{M}\cong\underline{M\<s\>}$ as $\underline{A}$-modules, for
  all $s\in\Z$.

  Let $\Hom_A(M,N)$ be the space of (degree preserving) $A$-module homomorphisms
  and set
  $$\ZHom_A(M,N)=\bigoplus_{d\in\Z}\ZHom_A(M,N\<d\>)
                \cong \bigoplus_{d\in\Z}\ZHom_A(M\<-d\>,N).$$
  The reader may check that
  $\ZHom_A(M,N)\cong\Hom_{\underline{A}}(\underline{M},\underline{N})$ as
  $\Zcal$-modules.

  Suppose that~$q$ is an indeterminate and that $M$ is a graded module.
  The \textbf{graded dimension} of~$M$ is the Laurent polynomial $$\gdim
  M=\sum_{d\in\Z} (\dim M_d)q^d\in\N[q,q^{-1}].$$ If~$M$ is a graded
  $A$-module, and $D$ is an irreducible graded $A$-module, then the
  \textbf{graded decomposition number} is the Laurent polynomial
  $$[M:D]_q = \sum_{s\in\Z} [M:D\<s\>]\,q^s\in\N[q,q^{-1}].$$
  By definition, the (ungraded) decomposition multiplicity
  $[\underline{M}:\underline{D}]$ is given by evaluating $[M:D]_q$ at $ q=1$,

  Suppose that  $A$ is a graded algebra and that $\underline{m}$ is an
  (ungraded) $\underline{A}$-module. A graded \textbf{lift} of $\underline{m}$
  is any graded $A$-module $M$ such that $\underline{M}\cong\underline{m}$ as
  $\underline{A}$-modules. If $M$ is a graded lift of $\underline{m}$ then so is
  $M\<s\>$, for any $s\in\Z$, so graded lifts are not unique. By Fitting's Lemma,  if
  $\underline{m}$ is indecomposable then its graded lift, if it exists,
  is unique up to grading shift~\cite[Lemma~2.5.3]{BGS:Koszul}.

  Following~\cite{HuMathas:GradedCellular}, the theory of cellular algebras from
  \autoref{S:CellularAlgebras} extends to the graded setting in a natural way.

  \begin{Definition}[\protect{\cite[\S2]{HuMathas:GradedCellular}}]
    Suppose that $A$ is $\Z$-graded $\Zcal$-algebra that is free of finite rank
  over $\Zcal$. A \textbf{graded cell datum} for $A$ is a cell datum
  $(\Pcal,T,C)$ together with a \textit{degree function}
    $$\deg\map{\coprod_{\lambda\in\Pcal}T(\lambda)}\Z$$
    such that
    \begin{enumerate}
      \item[(GC$_d$)]the element $c_{\s\t}$ is homogeneous of degree $\deg
      c_{\s\t}=\deg(\s)+\deg(\t)$, for all $\lambda\in\Pcal$ and
      $\s,\t\in T(\lambda)$.
    \end{enumerate}
    Then, $A$ is a \textbf{graded cellular algebra} with
    \textbf{graded cellular basis} $\{c_{\s\t}\}$.
  \end{Definition}

  We use $\star$ for the \textit{homogeneous} cellular algebra involution of~$A$
  that is determined by $c_{\s\t}^\star=c_{\t\s}$, for $\s,\t\in T(\lambda)$.

  \begin{Example}(Toy example)
    The most basic example of a graded algebra is the truncated polynomial ring
    $A=F[x]/(x^{n+1})$, for some integer $n>0$, where $\deg x=2$. As an ungraded
    algebra, $\underline{A}$ has exactly one simple module, namely the field~$F$
    with $x$ acting as multiplication by zero. This algebra is a graded cellular
    algebra with $\Pcal=\{0,1,\dots,n\}$, with its natural order, and
    $T(d)=\{d\}$ and $c_{dd}=x^d$. The irreducible graded
    $A$-modules are $F\<d\>$, for $d\in\Z$, and $\gdim A=1+q^2+\dots+q^{2n}$.
  \end{Example}

  \begin{Example}
    Let $A=\text{Mat}_n(\Zcal)$ be the $\Zcal$-algebra of $n\times n$-matrices.
    The basis of matrix units $\set{e_{st}|1\le s,t\le n}$ is a cellular basis
    for~$A$, where $\Pcal=\{\heartsuit\}$ and $T(\heartsuit)=\{1,2,\dots,n\}$.
    We want to put a non-trivial grading on~$A$. Let
    $\{d_1,\dots,d_n\}\subset\Z$ be a set of integers such that
    $d_s+d_{n-s+1}=0$, for $1\le s\le n$. Set $c_{st}=e_{s(n-t+1)}$ and define a
    degree function $\deg\map{T(\heartsuit)}\Z$ by $\deg s=d_s$. Then
    $\set{c_{st}|1\le s,t\le n}$ is a graded cellular basis of~$A$ and
    $\gdim A=\sum_{s=1}^n q^{d_s}$.
    Consequently, semisimple algebras can have non-trivial gradings.
  \end{Example}

  Exactly as in \autoref{S:CellularAlgebras}, for each $\lambda\in\Pcal$ we
  obtain a \textit{graded} cell module $C^\lambda$ with homogeneous basis
  $\set{c_\t|\t\in T(\lambda)}$ and $\deg c_\t=\deg\t$. Generalizing
  \autoref{E:CellularBilinearForm}, the graded cell module $C^\lambda$ comes
  equipped with a \textit{homogeneous} symmetric bilinear form $\<\ ,\
  \>_\lambda$ of degree zero. Therefore, if $x$ and $y$ are homogeneous
  elements of $C^\lambda$ then $\<x,y\>_\lambda\ne0$ only if $\deg x+\deg y=0$.
  Moreover, $\<xa,y\>_\lambda=\<x,ya^\star\>_\lambda$, for all $x,y\in C^\lambda$
  and all~$a\in A$. Consequently,
  $$\rad C^\lambda=\set{x\in C^\lambda|\<x,y\>_\lambda=0
        \text{ for all }y\in C^\lambda}$$
  is a graded submodule of~$C^\lambda$ so that
  $D^\lambda=C^\lambda/\rad C^\lambda$ is a graded $A$-module.

  If $M$ is an $A$-module then its (graded) \textbf{dual} is the
  $A$-module
  \begin{equation}\label{E:dual}
    M^\circledast=\ZHom_{\Zcal}(M,\Zcal),
  \end{equation}
  with $A$-action $(f\cdot a)(m)=f(ma^\star)$, for $f\in M^\circledast$,
  $a\in A$ and $m\in M$.

  \begin{Theorem}[\protect{\cite[Theorem~2.10]{HuMathas:GradedCellular}}]
    \label{T:GradedSimples}
    Suppose that $\Zcal$ is a field and that~$A$ is a graded cellular algebra. Then:
    \begin{enumerate}
      \item If $\lambda\in\Pcal$ then $D^\lambda$ is either~$0$ or an
      absolutely irreducible graded $A$-module. If $D^\lambda\ne0$ then
      $(D^\lambda)^\circledast\cong D^\lambda$.
      \item Let $\mathcal{K}=\set{\mu\in\Pcal|D^\mu\ne0}$. Then
      $\set{D^\lambda\<s\>|\lambda\in\mathcal{K} \text{ and }s\in\Z}$
      is a complete set of pairwise non-isomorphic irreducible (graded) $A$-modules.
      \item If $\lambda\in\Pcal$ and $\mu\in\mathcal{K}$ then
      $[C^\lambda{:}D^\mu]_q\ne0$
      only if $\lambda\gedom\mu$. Moreover, $[C^\mu{:}D^\mu]_q=1$.
    \end{enumerate}
  \end{Theorem}

  Forgetting the grading, the basis $\{c_{\s\t}\}$ is still a cellular basis
  of~$\underline{A}$. Comparing \autoref{T:CellularSimples} and
  \autoref{T:GradedSimples} it follows that every (ungraded) irreducible
  $\underline{A}$-module has a graded lift that is unique up to shift.
  Conversely, if $D$ is an irreducible graded $A$-module then $\underline{D}$ is
  an irreducible $\underline{A}$-module. (This holds more generally whenever a
  grading is put on a finite dimensional algebra;
  see~\cite[Theorem~4.4.4]{NV:graded}.) It is an instructive exercise to prove
  that if $A$ is a finite dimensional graded algebra then every simple
  $\underline{A}$-module has a graded lift and, up to shift, every graded simple
  $A$-module is of this form.

  By \cite[Theorems 3.2 and 3.3]{GordonGreen:GradedArtin} every projective
  indecomposable $\H$-module has a graded lift. More generally, as shown in
  \cite[\S4]{NV:graded}, if $M$ is a finitely generated graded $A$-module then
  the Jacobson radical of $\underline{M}$ has a graded lift.

  The matrix
  $\mathbf{D}_A(q)=([C^\lambda:D^\mu]_q)_{\lambda\in\Pcal,\mu\in\mathcal{K}}$
  is the \textbf{graded decomposition matrix} of $A$. For each $\mu\in\mathcal{K}$ let
  $P^\mu$ be the projective cover of~$D^\mu$ in~$\rep(A)$. The matrix
  $\mathbf{C}_A(q)=([P^\lambda:D^\mu]_q)_{\lambda,\mu\in\mathcal{K}}$ is the
  \textbf{graded Cartan matrix} of~$A$.

  An $A$-module $M$ has a \textbf{cell filtration} if there exists a filtration
  $M=M_0\supset M_1\supset\dots\supset M_z\supset 0$ such
  that each subquotient $M_r/M_{r+1}$ is isomorphic, up to shift, to some
  graded cell module. Fixing isomorphisms $M_r/M_{r+1}\cong C^{\lambda_r}\<d_r\>$,
  for some $\lambda_r\in\Pcal$ and $d_r\in\Z$, define
  $(M:C^\lambda)_q=\sum_d m_dq^d$, where
        $m_d=\#\set{1\le r\le z|\lambda_r=\lambda \text{ and }d_r=d}$.
  In general, the multiplicities $(M:C^\lambda)_q$ depend upon the
  choice of filtration \textit{and} the labelling of the isomorphisms
  $M_r/M_{r+1}\cong C^{\lambda_r}\<d_r\>$ because the cell modules are not
  guaranteed to be pairwise non-isomorphic, even up to shift.

  \begin{Corollary}[\protect{\cite[Theorem 2.17]{HuMathas:GradedCellular}}]%
    \label{C:CartanSymmetric}%
    Suppose that $\Zcal=F$ is a field. If $\mu\in\mathcal{K}$ then $P^\mu$ has a
    cell filtration such that $(P^\mu:C^\lambda)_q=[C^\lambda:D^\bmu]_q$, for all
    $\lambda\in\Pcal$. Consequently,
    $\mathbf{C}_A(q)=\mathbf{D}_A(q)^{tr}\mathbf{D}_A(q)$
    is a symmetric matrix.
  \end{Corollary}

  \subsection{Cyclotomic quiver Hecke algebras}\label{S:CycQuiverHecke}
  We are now ready to define cyclotomic quiver Hecke algebras. We start by
  defining the affine versions of these algebras and then pass to the cyclotomic
  quotients. Throughout this section we will make extensive use of the Lie theoretic
  data that is attached to the quiver $\Gamma_e$ in \autoref{S:Quivers}.

  If $\beta\in Q^+$ let $I^\beta=\set{\bi\in I^n|\beta=\alpha_{i_1}+\dots+\alpha_{i_n}}$.

  \begin{Definition}[Khovanov and Lauda~\cite{KhovLaud:diagI,KhovLaud:diagII}
    and Rouquier~\cite{Rouq:2KM}] \label{D:QuiverRelations}
    Suppose that $n\ge0$, $e\ge1$, and $\beta\in Q^+$. The
    \textbf{quiver Hecke algebra}, or \textbf{Khovanov-Lauda--Rouquier algebra},
    $\Rn[\beta]=\Rn[\beta](\Zcal)$ of type~$\Gamma_e$ is the unital associative
    $\Zcal$-algebra with generators
    $$\{\psi_1,\dots,\psi_{n-1}\} \cup \{ y_1,\dots,y_n \}
                                 \cup \set{e(\bi)|\bi\in I^\beta}$$
    and relations
    \bgroup
      \setlength{\abovedisplayskip}{1pt}
      \setlength{\belowdisplayskip}{1pt}
      \begin{align*}
        e(\bi) e(\bj) &= \delta_{\bi\bj} e(\bi),
          &{\textstyle\sum_{\bi \in I^\beta}} e(\bi)&= 1,\\
        y_r e(\bi) &= e(\bi) y_r,
        &\psi_r e(\bi)&= e(s_r{\cdot}\bi) \psi_r,
        &y_r y_s &= y_s y_r,
      \end{align*}
      \begin{align*}
        \psi_r \psi_s &= \psi_s \psi_r,&&\text{if }|r-s|>1,\\
        \psi_r y_s  &= y_s \psi_r,&&\text{if }s \neq r,r+1,
    \end{align*}
    \begin{equation}\label{E:ypsi}
      \begin{aligned}
        \psi_r y_{r+1} e(\bi)&=(y_r\psi_r+\delta_{i_ri_{r+1}})e(\bi),\\
        y_{r+1}\psi_re(\bi)&=(\psi_r y_r+\delta_{i_ri_{r+1}})e(\bi),
      \end{aligned}
    \end{equation}
    \begin{equation}
      \psi_r^2e(\bi) = \begin{cases}\label{E:quadratic}
           (y_{r+1} - y_{r})(y_{r}-y_{r+1}) e(\bi),
                 &\text{if }i_r\rightleftarrows i_{r+1},\\
           (y_{r}-y_{r+1})e(\bi),&\text{if  }i_r\rightarrow i_{r+1},\\
           (y_{r+1} - y_{r})e(\bi),&\text{if }i_r\leftarrow i_{r+1},\\
           0,&\text{if }i_r = i_{r+1},\\
           e(\bi),&\text{otherwise},\\
        \end{cases}
    \end{equation}
    and $(\psi_{r}\psi_{r+1}\psi_{r}-\psi_{r+1}\psi_{r}\psi_{r+1})e(\bi)$ is
    equal to
    \begin{equation}
       \begin{cases}\label{E:braid}
      (y_r+y_{r+2}-2y_{r+1})e(\bi),&\text{if }i_{r+2}=i_r \rightleftarrows i_{r+1},\\
        -e(\bi),\hspace*{30mm} &\text{if }i_{r+2}=i_r\rightarrow i_{r+1} ,\\
         e(\bi), &\text{if }i_{r+2}=i_r\leftarrow i_{r+1},\\
      0,&\text{otherwise,}
    \end{cases}
    \end{equation}
  \egroup
  for $\bi,\bj\in I^\beta$ and all admissible $r$ and $s$.
  \end{Definition}

  Part of the point of these definitions is that $\Rn[\beta]$ is a
  $\Z$-graded algebra with degree function determined by
  $$\deg e(\bi)=0,\qquad \deg y_r=2\qquad\text{and}\qquad \deg
    \psi_s e(\bi)=-c_{i_s,i_{s+1}},$$
  for $1\le r\le n$, $1\le s<n$ and $\bi\in I^n$.

  Suppose that $n\ge0$. Then
  $I^n=\bigsqcup_\beta I^\beta$ is the decomposition of $I^n$ into a disjoint
  union  of $\Sym_n$-orbits. Define
  \begin{equation}\label{E:RnBlocks}
    \Rn=\bigoplus_{\beta\in Q^+}\Rn[\beta].
  \end{equation}
  Set $\beta = \sum_{\bi\in I^\beta}e(\bi)$, for $\beta\in Q^+$.
  Then $\Rn[\beta]=e_\beta\Rn e_\beta$ is a two-sided ideal of $\Rn$ and
  \autoref{E:RnBlocks} is the decomposition of~$\Rn$ into
  blocks. That is, $\Rn[\beta]$ is indecomposable for all $\beta\in Q^+$.

  Khovanov and Lauda~\cite{KhovLaud:diagI,KhovLaud:diagII} and
  Rouquier~\cite{Rouq:2KM} define quiver Hecke algebras for quivers of arbitrary
  type. In the short time since their inception a lot has been discovered about
  these algebras. The first important result is that these algebras categorify
  the negative part of the corresponding quantum
  group~\cite{KhovLaud:diagI,Rouquier:QuiverHecke2Lie,
  BK:GradedDecomp,VaragnoloVasserot:CatAffineKLR}.

  \begin{Remark}\label{R:QuiverHecke}
    We have defined only a special case of the quiver Hecke algebras defined in
    \cite{KhovLaud:diagI,Rouq:2KM}. In addition to allowing arbitrary quivers,
    Khovanov and Lauda allow a more general choice of signs.  Rouquier's
    definition, which is the most general, defines the quiver Hecke algebras in
    terms of a matrix $Q=(Q_{ij})_{i,j\in I}$ with entries in a polynomial ring
    $\Zcal[u,v]$ with the properties that $Q_{ii}=0$, $Q_{ij}$ is not a zero
    divisor in $\Zcal[u,v]$ for $i\ne j$ and $Q_{ij}(u,v)=Q_{ji}(v,u)$, for
    $i,j\in I$. For an arbitrary quiver~$\Gamma$,
    Rouquier~\cite[Definition~3.2.1]{Rouq:2KM} defines $\Rn[\beta](\Gamma)$ to be the
    algebra generated by $\psi_r,y_s,e(\bi)$ subject to the relations above
    except that the quadratic and braid relations are replaced with
    $\psi_r^2e(\bi)= Q_{i_r,i_{r+1}}(y_r,y_{r+1})e(\bi)$ and
    $(\psi_r\psi_{r+1}\psi_r-\psi_{r+1}\psi_r\psi_{r+1})e(\bi)$ is equal to
    $$\begin{cases}
        \frac{Q_{i_r,i_{r+1}}(y_r,y_{r+1})-Q_{i_r,i_{r+1}}(y_r,y_{r+1})}
        {y_{r+2}-y_r},&\text{if }i_{r+2}=i_r,\\
        0,&\text{otherwise}.
      \end{cases}
    $$
    The assumptions on $Q$ ensure that the last expression is a polynomial in
    the generators. In general, $y_re(\bi)$ is homogeneous of degree
    $(\alpha_{i_r},\alpha_{i_r})$, for $1\le r\le n$ and $\bi\in I^n$. Under
    some mild assumptions, the isomorphism type of~$\Rn[\beta]$ is independent of the
    choice of~$Q$ by \cite[Proposition~3.12]{Rouq:2KM}. We leave it to the
    reader to find a suitable matrix $Q$ for \autoref{D:QuiverRelations}.
  \end{Remark}

  For the rest of these notes for $w\in\Sym_n$ we arbitrarily fix a reduced
  expression $w=s_{r_1}\dots s_{r_k}$, with $1\le r_j<n$. Using this fixed
  reduced expression for $w$ define $\psi_w=\psi_{r_1}\dots\psi_{r_k}$.

  \begin{Example}\label{E:PsiReduced}
  As the $\psi$-generators of $\Rn$ do not satisfy the braid relations the
  element~$\psi_w$ will, in general, depend upon the choice of reduced
  expression for $w\in\Sym_n$.  For example, by \autoref{E:braid} if
  $e\ne2$, $n=3$ and $w=s_1s_2s_1=s_2s_1s_2$ then
  $\psi_1\psi_2\psi_1e(0,2,0)=\psi_2\psi_1\psi_2e(0,2,0)+e(0,2,0),$
  by \autoref{E:braid}. Therefore, these two reduced expressions
  determine different elements of~$\Rn$.
  \end{Example}

  Khovanov and Lauda~\cite[Theorem~2.5]{KhovLaud:diagI} and Rouquier
  \cite[Theorem~3.7]{Rouq:2KM} proved the following.

  \begin{Theorem}[\protect{Khovanov-Lauda~\cite{KhovLaud:diagI}
                 and Rouquier \cite{Rouq:2KM}}] \label{T:RnBasis}
   Suppose that $\beta\in Q^+$. Then $\Rn[\beta](\Zcal)$ is free as an
   $\Zcal$-algebra with homogeneous basis
   $\set{\psi_w y_1^{a_1}\dots y_n^{a_n}e(\bi)|w\in \Sym_n,\ a_1,\dots,a_n\in\N
         \text{ and }\bi\in I^\beta}.$
  \end{Theorem}

  Li~\cite[Theorem~4.3.10]{Li:PhD} has constructed a graded cellular basis
  of~$\Rn$.  In the special case when $e=\infty$, that Kleshchev, Loubert and
  Miemietz~\cite{KleshchevLoubertMiemietz} give a graded affine cellular basis
  of~$\Rn$, in the sense of K\"onig and Xi~\cite{KoenigXi:AffineCellular}.

  In these notes we are not directly concerned with the quiver Hecke
  algebras~$\Rn$. Rather, we are more interested in
  \textit{cyclotomic quotients} of these algebras.

  \begin{Definition}[Brundan-Kleshchev~\cite{BK:GradedKL}]\label{D:CycQuiverHecke}
    Suppose that $\Lambda\in P^+$. The \textbf{cyclotomic quiver Hecke
    algebra} of type $\Gamma_e$ and weight $\Lambda$ is the quotient algebra
    $\R=\Rn/\<y_1^{(\Lambda,\alpha_{i_1})}e(\bi)\mid\bi\in I^n\>.$
  \end{Definition}

We abuse notation and identify the KLR generators of $\Rn$ with
their images in~$\R$. That is, we consider the algebra $\R$ to be generated by
$\psi_1,\dots,\psi_{n-1},y_1,\dots,y_n$ and $e(\bi)$, for $\bi\in I^n$, subject
to the relations in \autoref{D:QuiverRelations} and \autoref{D:CycQuiverHecke}.
From this point onwards, $\Lambda\in P^+$.

When $\Lambda$ is a weight of level~$2$, the algebras~$\R$ first appeared in
the work of Brundan and Stroppel~\cite{BrundanStroppel:KhovanovIII} in their
series of papers on the Khovanov diagram algebras. In full generality, the
cyclotomic quotients of~$\Rn$ were introduced by
Khovanov-Lauda~\cite{KhovLaud:diagI} and Rouquier~\cite{Rouq:2KM}. Brundan and
Kleshchev were the first to systematically study the cyclotomic quiver Hecke
algebras~$\R$, for any $\Lambda\in P^+$.

Although we will not need this here we note that, rather than working
algebraically, it is often easier to work diagrammatically by identifying the
elements of~$\R$ with certain planar diagrams. In these diagrams, the end-points
of the strings are labeled by $\{1,2,\dots,n,1',2',\dots,n'\}$ and the strings
themselves are coloured by~$I^n$. For example, following~\cite{KhovLaud:diagI},
the KLR generators can be identified with the diagrams:
$$
  e(\bi)=\begin{braid}\tikzset{baseline=7mm}
    \draw (0,4)node[above]{$i_1$}--(0,0);
    \draw (1,4)node[above]{$i_2$}--(1,0);
    \draw[dots] (1.2,4)--(3.8,4);
    \draw[dots] (1.2,0)--(3.8,0);
    \draw (4,4)node[above]{$i_{n}$}--(4,0);
  \end{braid}
  \space  \psi_re(\bi)=
  \begin{braid}\tikzset{baseline=7mm}
    \draw (0,4)node[above]{$i_1$}--(0,0);
    \draw[dots] (0.2,4)--(1.8,4);
    \draw[dots] (0.2,0)--(1.8,0);
    \draw (2,4)node[above]{$i_{r-1}$}--(2,0);
    \draw (3,4)node[above]{$i_r$}--(4,0);
    \draw (4,4)node[above]{$i_{r+1}$}--(3,0);
    \draw (5,4)node[above]{}--(5,0);
    \draw[dots] (5.2,4)--(6.8,4);
    \draw[dots] (5.2,0)--(6.8,0);
    \draw (7,4)node[above]{$i_{n}$}--(7,0);
  \end{braid}
  \space  y_se(\bi)=
  \begin{braid}\tikzset{baseline=7mm}
    \draw (0,4)node[above]{$1$}--(0,0);
    \draw[dots] (0.2,4)--(1.8,4);
    \draw[dots] (0.2,0)--(1.8,0);
    \draw (2,4)node[above]{${s-1}$}--(2,0);
    \draw (3,4)node[above]{$s$}--(3,0);
    \greendot(3,2);
    \draw (4,4)node[above]{$s$}--(4,0);
    \draw[dots] (4.2,4)--(5.8,4);
    \draw[dots] (4.2,0)--(5.8,0);
    \draw (6,4)node[above]{${n}$}--(6,0);
  \end{braid}.
$$
Multiplication of diagrams is given by concatenation, read from top to bottom,
subject to the relations above that are also interpreted diagrammatically. As an
exercise, we leave it to the reader to identify the two relations in
\autoref{D:QuiverRelations} that correspond to the following `local' relations on
strings inside braid diagrams:
  \begin{align*}
    \begin{braid}[7]
    \draw(1,4) node[above]{$i$}--(3,0);
    \draw(3,4) node[above]{$j$}--(1,0);
    \greendot(2.5,1);
  \end{braid}
  =
  \begin{braid}[7]
    \draw(1,4) node[above]{$i$}--(3,0);
    \draw(3,4) node[above]{$j$}--(1,0);
    \greendot(1.5,3);
  \end{braid}
  +\delta_{ij}\begin{braid}[7]
    \YCrossing i 1 3 1 3
  \end{braid}
  \quad\text{and}\quad
    \begin{braid}[7]
    \draw(1,4)node[above]{$i$}--(3,0);
    \draw(2,4)node[above]{$i{\pm}1$}--(3,2)--(2,0);
    \draw(3,4)node[above]{$i$}--(1,0);
  \end{braid}
  =
  \begin{braid}[7]
    \draw(1,4)node[above]{$i$}--(3,0);
    \draw(2,4)node[above]{$i{\pm}1$}--(1,2)--(2,0);
    \draw(3,4)node[above]{$i$}--(1,0);
  \end{braid}
  \pm
  \begin{braid}[7]
    \draw(1,4)node[above]{$i$}--(1,0);
    \draw(2,4)node[above]{$i{\pm}1$}--(2,0);
    \draw(3,4)node[above]{$i$}--(3,0);
  \end{braid}.
  \end{align*}
(For the second relation, $e\ne2$.) For more rigorous definitions of such
diagrams, and non-trivial examples of their application, we refer the
reader to the papers
\cite{HoffnungLauda:KLRnilpotency,KMR:UniversalSpecht,LyleMathas:CycCP,Li:PhD}
for examples of these diagrams in action.

  \begin{Example}[Rank one algebras]
    Suppose that $n=1$ and $\Lambda\in P^+$. Then
    $\R[1]=\<y_1,e(i)\mid y_1e(i)=e(i)y_1\text{ and }
                        y_1^{\<\Lambda,\alpha_i\>}e(i)=0, \text{ for }i\in I\>,$
    with $\deg y_1=2$ and $\deg e(i)=0$, for $i\in I$. Therefore,
    there is an isomorphism of graded algebras
    $$\R[1]\cong\bigoplus_{\substack{i\in I\\(\Lambda,\alpha_i)>0}}
    \Zcal[y]/y^{(\Lambda,\alpha_i)}\Zcal[y],$$
    where $y=y_1$ is in degree~$2$. Armed with this
    description of~$\R$ it is now straightforward to show that
    $\H\cong\R$ when~$\Zcal$ is a field and~$n=1$.
  \end{Example}

  \subsection{Nilpotence and small representations}
  In this section and the next we use the KLR relations to prove some
  results about the cyclotomic quiver Hecke algebras~$\R$ for particular
  $\Lambda$ and~$n$.

  By \autoref{T:RnBasis} the algebra $\Rn$ is infinite dimensional, so it is not
  obvious from the relations that the cyclotomic Hecke algebra~$\R$ is finite
  dimensional --- or even that~$\R$ is non-zero. The following result shows
  that~$y_r$ is nilpotent, for $1\le r\le n$, which implies that $\R$ is finite
  dimensional.

  \begin{Lemma}[\protect{Brundan and Kleshchev~\cite[Lemma~2.1]{BK:GradedKL}}]
    \label{L:nilpotence}
    Suppose that $1\le r\le n$ and $\bi\in I^n$.
    Then $y_r^Ne(\bi)=0$ for $N\gg0$.
  \end{Lemma}

  \begin{proof}
    We argue by induction on~$r$. If $r=1$ then
    $y_1^{(\Lambda,\alpha_{i_1})}e(\bi)=0$ by \autoref{D:CycQuiverHecke},
    proving the base step of the induction. Now consider $y_{r+1}e(\bi)$. By
    induction, we may assume that there exists $N\gg0$ such that
    $y_r^Ne(\bj)=0$, for all $\bj\in I^n$.  There are three cases to consider.

    \Case{$i_{r+1}\noedge i_r$}
    By \autoref{E:quadratic} and \autoref{E:ypsi},
    $$y_{r+1}^Ne(\bi)=y_{r+1}^N\psi_r^2e(\bi)
           =\psi_r y_r^N\psi_re(\bi)=\psi_r y_r^Ne(s_r\cdot\bi)\psi_r=0,$$
    where the last equality follows by induction.

    \Case{$i_{r+1}=i_r\pm1$}
    Suppose first that $e\ne2$. This is a variation on the previous case, with a
    twist.  By \autoref{E:quadratic} and \autoref{E:ypsi}, again
    \begin{align*}
      y_{r+1}^{2N}e(\bi)&=y_{r+1}^{2N-1}y_re(\bi)+y_{r+1}^{2N-1}(y_{r+1}-y_r)e(\bi)\\
         &=y_ry_{r+1}^{2N-1}e(\bi)\pm y_{r+1}^{2N-1}\psi_r^2e(\bi)\\
         &=y_ry_{r+1}^{2N-1}e(\bi)\pm\psi_ry_r^{2N-1}e(s_r\cdot\bi)\psi_r\\
         &=y_ry_{r+1}^{2N-1}e(\bi)=\dots=y_r^Ny_{r+1}^Ne(\bi)=0.
    \end{align*}
    The case when $e=2$ is similar. First, observe that
    $y_{r+1}^2e(\bi)=(2y_ry_{r+1}-y_r^2-\psi_r^2)e(\bi)$
    by~\autoref{E:quadratic}.  Therefore, arguing as before,
    $$y_{r+1}^{3N}e(\bi)=y_r(2y_{r+1}-y_r)y_{r+1}^{3N-2}e(\bi)
       =\dots=y_r^N(2y_{r+1}-y_r)^Ny_{r+1}^Ne(\bi)=0.$$

    \Case{$i_{r+1}=i_r$}
    Let $\phi_r=\psi_r(y_r-y_{r+1})$.  Then $\phi_r\psi_re(\bi)=-2\psi_re(\bi)$
    by \autoref{E:ypsi}, so that $(1+\phi_r)^2e(\bi)=e(\bi)$. Moreover,
    $$
      (1+\phi_r)y_r(1+\phi_r)e(\bi)
       =(y_r+\phi_ry_r+y_r\phi_r+\phi_ry_r\phi_r)e(\bi)\\
       =y_{r+1}e(\bi),
    $$
    where the last equality uses \autoref{E:ypsi}. Now
    we are done because
    $$y_{r+1}^Ne(\bi)=\((1+\phi_r)y_r(1+\phi_r)\)^Ne(\bi)
               =(1+\phi_r)y_r^N(1+\phi_r)e(\bi)=0,$$
    since $\phi_r$ commutes with $e(\bi)$ and $y_r^Ne(\bi)=0$ by induction.
  \end{proof}

  We have marginally improved on Brundan and Kleshchev's original proof
  of \autoref{L:nilpotence} because the argument above gives an upper
  bound for the nilpotency index of~$y_r$.  In general, this bound is
  far from sharp.  For a better estimate of the nilpotency index
  of~$y_r$ see \cite[Corollary~4.6]{HuMathas:SeminormalQuiver} (and
  \cite{HoffnungLauda:KLRnilpotency} when $e=\infty$). See
  \cite[Lemma~4.4]{KangKashiwara:CatCycKLR} for another argument that
  applies to cyclotomic quiver Hecke algebras of arbitrary type.

  Combining \autoref{T:RnBasis} and \autoref{L:nilpotence} we have:

  \begin{Corollary}[\protect{Brundan and Kleshchev~\cite[Corollary~2.2]{BK:GradedKL}}]
    \label{C:FiniteDimensional}
    Suppose $\Zcal$ is an integral domain. Then $\R$ is finite dimensional.
  \end{Corollary}

%
%


  As our next exercise we classify the one dimensional representations of~$\R$
  when $\Zcal=F$ is a field. For $i\in I$ let $\bi_n^+=(i,i+1,\dots,i+n-1)$ and
  $\bi_n^-=(i,i-1,\dots,i-n+1)$. Then $\bi_n^\pm\in I^n$. If
  $(\Lambda,\alpha_i)=0$ then $e(\bi^\pm_n)=0$ by \autoref{D:CycQuiverHecke}.
  However, if $(\Lambda,\alpha_i)\ne0$ then using the relations it is easy to
  see that $\R$ has unique one dimensional representations
  $D^+_{i,n}=Fd^+_{i,n}$ and $D^-_{i,n}=Fd^-_{i,n}$ such that
  $$d_{i,n}^\pm e(\bi) = \delta_{\bi,\bi^\pm_n}d_{i,n}^\pm\quad\text{and}\quad
  d_{i,n}^+y_r=0=d_{i,n}^\pm\psi_s,$$
  for $\bi\in I^n$, $1\le r\le n$ and $1\le s<n$ and such that $\deg d_{i,n}^\pm=0$.
  In particular, this shows that $e(\bi^\pm_n)\ne0$ and hence that
  $\R\ne0$. If $e\ne2$ then $\set{D_{i,n}^\pm|i\in I\text{ and }(\Lambda,\alpha_i)\ne0}$
  are pairwise non-isomorphic irreducible representations of~$\R$. If
  $e=2$ then $\bi^+_n=\bi_n^-$ so that $D_{i,n}^+=D_{i,n}^-$.

  \begin{Proposition}\label{P:OneDimensionals}
    Suppose that $\Zcal=F$ is a field and that $D$ is a one dimensional graded
    $\R$-module. Then $D\cong D^\pm_{i,n}\<k\>$, for some $k\in\Z$ and
    $i\in I$ such that $(\Lambda,\alpha_i)\ne0$.
  \end{Proposition}

  \begin{proof}Let $d$ be a non-zero element of $D$ so that $D=Fd$. Then
    $d=\sum_{\bj\in I^n}de(\bj)$ so that $de(\bi)\ne0$ for some $\bi\in I^n$.
    Moreover, $de(\bj)=0$ if and only if $\bj=\bi$ since otherwise $de(\bi)$ and
    $de(\bj)$ are linearly independent elements of $D$, contradicting
    assumption that $D$ is one dimensional.  Now, $\deg dy_r=2+\deg d$, so
    $dy_r=0$, for $1\le r\le n$, since $D$ is one dimensional. Similarly,
    $d\psi_r=de(\bi)\psi_r=0$ if $i_r=i_{r+1}$ or $i_r=i_{r+1}\pm1$ since in these
    cases $\deg e(\bi)\psi_r\ne0$.

    It remains to show that $\bi=\bi^\pm_n$ and that $(\Lambda,\alpha_{i_1})\ne0$.
    First, since $0\ne d=d\,e(\bi)$ we have that $e(\bi)\ne0$ so that
    $(\Lambda,\alpha_{i_1})\ne0$ by \autoref{D:CycQuiverHecke}. To complete the
    proof we show that if $\bi\ne\bi^\pm_n$ then $d=0$, which is a
    contradiction.  First, suppose that $i_r=i_{r+1}$ for some $r$, with
    $1\le r<n$. Then $d=de(\bi)=d\(\psi_ry_{r+1}-y_r\psi_r\)e(\bi)=0$ by
    \autoref{E:ypsi}, which is not possible, so $i_r\ne i_{r+1}$. Next, suppose
    that $i_{r+1}\ne i_r\pm1$. Then
    $d=de(\bi)=d\psi_r^2e(\bi)=d\psi_re(s_r\cdot\bi)\psi_r=0$ because $D$ is one
    dimensional and $de(\bj)=0$ if $\bj\ne\bi$. This is another contradiction,
    so we must have $i_{r+1}=i_r\pm1$ for $1\le r<n$.  Therefore, if
    $\bi\ne\bi^\pm_n$ then $e\ne2$, $n>2$ and $i_r=i_{r+2}=i_{r+1}\pm1$ for some $r$.
    Applying the braid relation \autoref{E:braid},
    $$d=de(\bi)=\pm d\cdot(\psi_r\psi_{r+1}\psi_r-\psi_{r+1}\psi_r\psi_{r+1})e(\bi)
       =0,$$
    a contradiction. Hence, $D\cong D^\pm_{i,n}\<\deg d\>$, completing the proof.
  \end{proof}

  \subsection{Semisimple KLR algebras}
  Now that we understand the one dimensional representations of $\R$ we consider
  the semisimple representation theory of the cyclotomic quiver Hecke
  algebras.  These results do not appear in the literature, but there
  are few surprises here because everything we do can be easily deduced
  from results that are known. The main idea is to show by example how
  to use the quiver Hecke algebra relations.

  In this section we fix $e>n$ and $\Lambda\in P^+$ such that
  $(\Lambda,\alpha_{i,n})\le1$, for all~$i\in I$, and we study the
  algebras~$\R$. Notice that these conditions ensure that~$\H$ is
  semisimple by \autoref{C:IntegralSSimple}.

  Recall from \autoref{S:blocks} that $\bi^\t=(i^\t_1,\dots,i^\t_n)$ is the
  residue sequence of $\t\in\Std(\Parts)$, where $i^\t_r=c^\Z_r(\t)+e\Z$.  We
  caution the reader that if~$\t$ is a standard tableau then the contents
  $c^\Z_r(\t)\in\Z$ and the residues $i^\t_r\in I$ are in general different.

  If $i\in I$ then a node $A=(l,r,c)$ is an \textbf{$i$-node} if
  $i=\kappa_l+c-r+e\Z$. Therefore, extending the definitions of
  \autoref{S:Tableaux}, we can now talk of addable and removable $i$-nodes.

  \begin{Lemma}\label{L:UniqueResidues}
    Suppose that $e>n$ and $(\Lambda,\alpha_{i,n})\le1$, for all $i\in I$. Let
    $\s,\t\in\Std(\Parts)$. Then $\s=\t$ if and only if $\bi^\s=\bi^\t$.
  \end{Lemma}

  \begin{proof}Observe that if $i\in I$ and $\bmu\in\Parts[m]$,
    where $0\le m<n$, then $\bmu$ has at most one addable $i$-node since
    $(\Lambda,\alpha_{i,n})\le 1$. Hence, it follows easily by
    induction on~$n$ that $\s=\t$ if and only if\/ $\bi^s=\bi^\t$.
  \end{proof}

  \autoref{L:UniqueResidues} also follows from \autoref{T:SSimple} and
  \autoref{C:IntegralSSimple}.

  Let $I^n_\Lambda=\set{\bi^\t|\t\in\Std(\Parts)}$ be the set of residue
  sequences of all of the standard tableaux in $\Std(\Parts)$.  By the
  proof of \autoref{L:UniqueResidues}, if $\bi=\bi^\t\in I^n_\Lambda$
  and $i_{r+1}=i_r\pm1$ then $r$ and $r+1$ must be in either in the same
  row or in the same column of~$\t$.  Hence, we have the following
  useful fact.

  \begin{Corollary}\label{C:SameRowColumn}
    Suppose that $e>n$ and that $(\Lambda,\alpha_{i,n})\le1$, for all
    $i\in I$, and that $\bi\in I^n_\Lambda$ such that $i_{r+1}=i_r\pm1$.
    Then $s_r\cdot\bi\notin I^n_\Lambda$.
  \end{Corollary}

  When $\Lambda=\Lambda_0$ the next result is due to Brundan and
  Kleshchev~\cite[\S5.5]{BK:GradedKL}. More generally, Kleshchev and
  Ram~\cite[Theorem~3.4]{KleshRam:PureIrredRepsKLR} prove similar
  results for quiver Hecke algebras of simply laced type.

  \begin{Proposition}[Seminormal representations of $\R$]\label{P:KLRSemisimpleReps}
    Suppose that $\Zcal=F$ is a field, $e>n$ and that $\Lambda\in P^+$
    with $(\Lambda,\alpha_{i,n})\le1$, for all $i\in I$. Then for each
    $\blam\in\Parts$ there is a unique irreducible graded $\R$-module $S^\blam$
    with homogeneous basis $\set{\psi_\t|\t\in\Std(\blam)}$ such that $\deg
    \psi_\t=0$, for all $\t\in\Std(\blam)$, and where the $\R$-action is given
    by
    $$ \psi_\t e(\bi)=\delta_{\bi,\bi^\t} \psi_\t,\qquad
       \psi_\t y_r=0\qquad\text{and}\qquad
       \psi_\t \psi_r=\psi_{\t(r,r+1)},$$
    where we set $\psi_{\t(r,r+1)}=0$ if $\t(r,r+1)$ is not standard.
  \end{Proposition}

  \begin{proof}
    By \autoref{L:UniqueResidues}, if $\s,\t\in\Std(\blam)$ then $\s=\t$ if and
    only if $\bi^\s=\bi^\t$.  Moreover, $i^\t_{r+1}=i^\t_r\pm1$ if and only if
    $r$ and $r+1$ are in the same row or in the same column of~$\t$.  Similarly,
    $i^\t_r\ne i^\t_{r+1}$ for any $r$. Consequently, since $\psi_\t=\psi_\t
    e(\bi^\t)$ almost all of the relations in \autoref{D:QuiverRelations} are
    trivially satisfied. In fact, all that we need to check is that
    $\psi_1,\dots,\psi_{n-1}$ satisfy the braid relations of the symmetric
    group~$\Sym_n$ with $\psi_r^2$ acting as zero when $i^\t_{r+1}=i^\t_r\pm1$,
    which follows automatically by \autoref{C:SameRowColumn}.  By the same
    reasoning if $\t(r,r+1)$ is standard then $\deg e(\bi^\t)\psi_r=0$. Hence,
    we can set $\deg \psi_\t=0$, for all $\t\in\Std(\blam)$. This proves that
    $S^\blam$ is a graded $\R$-module.

    It remains to show that $S^\blam$ is irreducible. If $\s,\t\in\Std(\blam)$
    then $\s=\tlam d(\s) = \t d(\t)^{-1}d(\s)$, so $\psi_\s=\psi_\t
    \psi_{d(\t)^{-1}}\psi_{d(\s)}$. Suppose that $x=\sum_\t r_\t \psi_\t$ is a
    non-zero  element of~$S^\blam$. If $r_\t\ne0$ then
    $\psi_\t=\tfrac1{r_\t}xe(\bi^\t)$, so it follows that $\psi_\s\in x\R$, for any
    $\s\in\Std(\blam)$. Therefore, $S^\blam=x\R$ so that $S^\blam$ is
    irreducible as claimed.
  \end{proof}

  Hence, $e(\bi)\ne0$ in $\R$, for all $\bi\in I^\Lambda_n$. This
  was not clear until now.

  We want to show that \autoref{P:KLRSemisimpleReps} describes all of the graded
  irreducible representations of~$\R$, up to degree shift. To do this we need a
  better understanding of the set $I^n_\Lambda$. Okounkov and
  Vershik~\cite[Theorem~6.7]{OkounkovVershik} explicitly described the set of
  all \textit{content sequences} $(c^\Z_1(\t),\dots,c^\Z_n(\t))$ when $\ell=1$.
  This combinatorial result easily extends to higher levels and so suggests a
  description of~$I^n_\Lambda$.

  If $\bi\in I^n$ and $1\le m\le n$ set $\bi_{\downarrow m}=(i_1,\dots,i_m)$.
  Then $\bi_{\downarrow m}\in I^m$ and
  $I^m_\Lambda=\set{\bi_{\downarrow m}|\bi\in I^n_\Lambda}$.

  \begin{Lemma}[\protect{\textit{cf.}
    Ogievetsky-d'Andecy~\cite[Proposition~5]{OgievetskydAndecy}}]
    Suppose that $e>n$ and $(\Lambda,\alpha_{i,n})\le1$, for all $i\in I$. Let
    $\bi\in I^n$.  Then $\bi\in I^n_\Lambda$ if and only if it
    satisfies the following three conditions:
    \label{L:SSResidues}
    \begin{enumerate}
      \item $(\Lambda,\alpha_{i_1})\ne0$.
      \item If $1<r\le n$ and $(\Lambda,\alpha_{i_r})=0$ then
         $\{i_r-1,i_r+1\}\cap\{i_1,\dots,i_{r-1}\}\ne\emptyset$.
      \item If $1\le s<r\le n$ and $i_r=i_s$ then
         $\{i_r-1,i_r+1\}\subseteq\{i_{s+1},\dots,i_{r-1}\}$.
    \end{enumerate}
  \end{Lemma}

  \begin{proof}
    Suppose that $\t\in\Std(\Parts)$ and let $\bi=\bi^\t$. We prove by induction
    on~$r$ that $\bi_{\downarrow r}\in I^r_\Lambda$. By definition,
    $i_1=\kappa_t+e\Z$ for some~$t$ with $1\le t\le\ell$, so (a) holds. By
    induction we may assume that the subsequence $(i_1,\dots,i_{r-1})$ satisfies
    properties (a)--(c).  If $(\Lambda,\alpha_{i_r})=0$ then~$r$ is not in
    the first row or in the first column of any component of~$\t$, so $\t$ has an
    entry in the row directly above~$r$ or in the column immediately to the
    left of~$r$ --- or both! Hence, there exists an integer $s$ with $1\le s<r$
    such that $i^\t_s=i^\t_r\pm1$. Hence, (b) holds. Finally, suppose that
    $i_r=i_s$ as in~(c).  As  the residues of the nodes in different components
    of~$\t$ are disjoint it follows that~$s$ and $r$ are in same component
    of~$\t$ and on the same diagonal. In particular, $r$ is not in the first row
    or in the first column of its component in~$\t$. As~$\t$ is standard, the entries
    in~$\t$ that are immediately above or to the left of~$r$ are both larger
    than~$s$ and smaller than~$r$. Hence,~(c) holds.

    Conversely, suppose that $\bi\in I^n$ satisfies properties (a)--(c).
    We show by induction on $m$ that $\bi_{\downarrow m}\in
    I^m_\Lambda$, for $1\le m\le n$. If $m=1$ then $\bi_{\downarrow1}\in
    I^1_\Lambda$ by property~(a). Now suppose that $1<m<n$ and that
    $\bi_{\downarrow m}\in I^m_\Lambda$. By induction $\bi_{\downarrow
    m}=\bi^\s$, for some $\s\in\Std(\Parts[m])$. Let $\bnu=\Shape(\s)$.
    If $i\in I$ then $(\Lambda,\alpha_{i,n})\le1$, so the
    multipartition~$\bnu$ can have at most one addable $i$-node. On the
    other hand, reversing the argument of the last paragraph, using
    properties~(b) and (c) with $r=m+1$, shows that~$\bnu$ has at least
    one addable $i_{m+1}$-node.  Let $A$ be the unique addable
    $i_{m+1}$-node of $\bnu$. Then $\bi_{\downarrow(m+1)}=\bi^\t$ where
    $\t\in\Std(\Parts[m+1])$ is the unique standard tableau such that
    $\t_{\downarrow m}=\s$ and $\t(A)=m+1$. Hence, $\bi\in
    I^{m+1}_\Lambda$ as required.
  \end{proof}

  By \autoref{P:KLRSemisimpleReps}, if $\bi\in I^n_\Lambda$ then $e(\bi)\ne0$.
  We use \autoref{L:SSResidues} to show that $e(\bi)=0$ if $\bi\notin
  I^n_\Lambda$. First, a result that holds for all $\Lambda\in P^+$.

  \begin{Lemma}\label{L:gaps}
    Suppose that $\Lambda\in P^+$ and that $e(\bi)\ne0$, for $\bi\in I^n$.
    Then $(\Lambda,\alpha_{i_1})\ne0$.  Moreover,
    $\{i_r-1,i_r+1\}\cap\{i_1,\dots,i_{r-1}\}\ne\emptyset$ whenever
    $(\Lambda,\alpha_{i_r})=0$ for some $1<r\le n$.
  \end{Lemma}

  \begin{proof}By \autoref{D:CycQuiverHecke}, $e(\bi)=0$ whenever
    $(\Lambda,\alpha_{i_1})=0$. To prove the second claim
    suppose that  $(\Lambda,\alpha_{i_r})=0$ and $i_r\pm1\notin\{i_1,\dots,i_{r-1}\}$.
    We may assume that $i_r\ne i_s$ for $1\le s<r$. Applying
    \autoref{E:quadratic} $r$-times,
    \begin{align*}
     e(\bi)&=\psi_{r-1}^2e(\bi)=\psi_{r-1}e(i_1,\dots,i_r,i_{r-1},i_{r+1},\dots,i_n)
               \psi_{r-1}\\
         &=\dots=\psi_{r-1}\dots\psi_1
             e(i_r,i_1,\dots,i_{r-1},i_{r+1},\dots,i_n)\psi_1\dots\psi_{r-1}=0,
    \end{align*}
    where the last equality follows because $(\Lambda,\alpha_{i_r})=0$.
  \end{proof}

  \begin{Proposition}\label{P:SSimpleIdempotents}
    Suppose that $1\le m\le n$ and that $(\Lambda,\alpha_{i,m})\le1$, for all
    $i\in I$. Then $y_1=\dots=y_m=0$. Moreover, if $\bi\in I^n$ then $e(\bi)\ne0$ only
    if\/ $\bi_{\downarrow m}\in I^m_\Lambda$.
  \end{Proposition}

  \begin{proof}We argue by induction on~$r$ to show that $y_r=0$ and $e(\bi)=0$
    if $\bi_{\downarrow r}\notin I^r_\Lambda$, for $1\le r\le m$. If $r=1$ this is immediate
    because $y_1^{(\Lambda,\alpha_{i_1})}e(\bi)=0$ by \autoref{D:CycQuiverHecke}
    and $(\Lambda,\alpha_{i_1})\le1$ by assumption.
    Suppose then that $1<r\le m$.

    We first show that $e(\bi)=0$ if $\bi_{\downarrow r}\notin I^r_\Lambda$. By induction,
    \autoref{L:SSResidues} and \autoref{L:gaps}, it is enough to show that
    $e(\bi)=0$ whenever there exists an integer $1\le s<r$ such that $i_s=i_r$ and
    $\{i_r-1,i_r+1\}\not\subseteq\{i_{s+1},\dots,i_{r-1}\}$. We may assume that~$s$
    is maximal such that $i_s=i_r$ and $1\le s<r$. There are three cases to
    consider.

    \Case{$r=s+1$} By \autoref{E:ypsi},
    $e(\bi)=(y_{s+1}\psi_s-\psi_sy_s)e(\bi)=y_{s+1}\psi_se(\bi),$
    since $y_s=0$ by induction. Using this identity twice, reveals that
    $e(\bi)=y_{s+1}\psi_se(\bi)=y_{s+1}e(\bi)\psi_s
       =y_{s+1}^2\psi_se(\bi)\psi_s
       =y_{s+1}^2\psi_s^2e(\bi)=0,$
     where the last equality comes from \autoref{E:quadratic}. Therefore,
     $e(\bi)=0$ as we wanted to show.

    \Case{$s<r-1$ and $\{i_r-1,i_r+1\}\cap\{i_{s+1},\dots,i_{r-1}\}=\emptyset$}
    By the maximality of~$s$, $i_r\notin\{i_{s+1},\dots,i_{r-1}\}$. Therefore,
    arguing as in the proof of \autoref{L:gaps}, there exists a permutation
    $w\in\Sym_r$ such that
    $e(\bi)
      =\psi_we(i_1,\dots,i_s,i_r,i_{s+1},\dots,i_{r-1},i_{r+1},\dots,i_n)\psi_w.$
    Hence, $e(\bi)=0$ by Case~1.

    \Case{$s<r-1$ and $\{i_r-1,i_r+1\}\cap\{i_{s+1},\dots,i_{r-1}\}=\{j\}$,
    where $j=i_r\pm1$} Let $t$ be an index such that $i_t=j=i_r\pm1$ and
    $s<t<r$. Note that if there exists an integer~$t'$ such that $i_t=i_{t'}$
    and $s<t<t'<r$ then we may assume that $i_s\in\{i_{t+1},\dots,i_{t'-1}\}$ by
    \autoref{L:SSResidues}(c) and induction. Therefore, since~$s$ was chosen to
    be maximal, $t$ is the unique integer such that $i_t=j$ and $s<t<r$.  Hence,
    arguing as in Case~2, there exists a permutation $w\in\Sym_r$ such that
    $$e(\bi)=\psi_we(\dots,i_{s-1},i_{s+1},\dots,i_{t-1},i_s,i_t,i_r,i_{t+1},
             \dots, i_{r-1},i_{r+1},\dots)\psi_w.$$
    For convenience, we identify $e(i_1,\dots,i_s,i_t,i_r,\dots,i_n)$ with
    $e(i,j,i)$, where $i=i_s=i_r$ and $j=i\pm1$. Then, by \autoref{E:braid},
    \begin{align*}
      e(i,j,i)&=\pm\(\psi_1\psi_2\psi_1-\psi_2\psi_1\psi_2\)e(i,j,i)\\
        &=\pm\psi_1\psi_2e(j,i,i)\psi_1\mp\psi_2\psi_1e(i,i,j)\psi_2=0,
    \end{align*}
    where the last equality follows by Case~1.


    Combining Cases~1-3, if $e(\bi)\ne0$ then
    $\{i_r-1,i_r+1\}\subseteq\{i_{s+1},\dots,i_{r-1}\}$ whenever
    there exists an integer~$s$ such that $i_s=i_r$ and $1\le s<r$.
    Hence, as remarked above, induction, \autoref{L:gaps} and
    \autoref{L:SSResidues} show that $e(\bi)\ne0$ only if $\bi_{\downarrow r}\in I^r_\Lambda$.

    \medskip
    To complete the proof of the inductive step (and of the proposition), it
    remains to show that $y_r=0$. Using what we have just proved, it is enough
    to show that $y_re(\bi)=0$ whenever $\bi_{\downarrow r}\in I^r_\Lambda$. If
    $i_{r-1}=i_r\pm1$ then, by induction and \autoref{E:quadratic},
    $$y_re(\bi) = (y_r-y_{r-1})e(\bi) = \pm\psi_{r-1}^2e(\bi)
    = \pm\psi_{r-1}e(s_{r-1}\cdot\bi)\psi_{r-1}=0,$$
    where the last equality follows because $(s_{r-1}\cdot\bi)_{\downarrow r}\notin I^r_\Lambda$
    by \autoref{C:SameRowColumn}. If $i_{r-1}\ne i_r\pm1$
    then $i_{r-1}\noedge i_r$ by \autoref{L:SSResidues}
    since $\bi_{\downarrow r}\in I^r_\Lambda$. Therefore,
    $y_re(\bi) = y_r\psi_{r-1}^2e(\bi)=\psi_{r-1}y_{r-1}\psi_{r-1}e(\bi)=0$
    since $y_{r-1}=0$ by induction. This completes the proof.
  \end{proof}

  Before giving our main application of \autoref{P:SSimpleIdempotents}
  we interpret it means for the cyclotomic quiver Hecke algebras of the
  symmetric groups.

  \begin{Example}[Symmetric groups]
    Suppose that $\Lambda=\Lambda_0$, $n\ge0$ and set $f=\min\{e,n\}$. Then
    $(\Lambda,\alpha_{i,f-1})\le1$ for all $i\in I$. Therefore,
    \autoref{P:SSimpleIdempotents} shows that $y_r=0$ for
    $1\le r<f$ and that $e(\bi)\ne0$ only if
    $\bi_{\downarrow(f-1)}\in I^{f-1}_\Lambda$. In addition, we also have $\psi_1=0$ because
    if $\bi\in I^n$ then $\psi_1e(\bi)=e(s_1\cdot\bi)\psi_1=0$ because if
    $\bi_{\downarrow(f-1)}\in I^{f-1}_\Lambda$ then $(s_1\cdot\bi)_{f-1}\notin I^{f-1}_\Lambda$.

    Translating the proof of \autoref{P:SSimpleIdempotents} back to
    \autoref{L:UniqueResidues}, the reason why $\psi_1=0$ is that if
    $\bi=\bi^\t$ is the residue sequence of some standard tableau
    $\t\in\Std(\Parts)$ then $i_1=0$ and $i_2\ne0$, so that
    $s_1\cdot\bi\notin I^\Lambda_n$ is not a residue sequence and,
    consequently, $\psi_1 e(\bi)=e(s_1\cdot \bi)\psi_1=0$. By the same reasoning,
    $\psi_1\ne0$ if $\Lambda$ has level $\ell>1$.
  \end{Example}

  We now completely describe the structure of the KLR algebras $\R$ when
  $e>n$ and $\Lambda\in P^+$ such that $(\Lambda,\alpha_{i,n})\le 1$,
  for all $i\in I$. For $(\s,\t)\in\Std^2(\Parts)$ define
  $e_{\s\t}=\psi_{d(\s)^{-1}}e(\bi^{\blam})\psi_{d(\t)}$, where
  $\bi^\blam=\bi^{\tlam}$.

  \begin{Theorem}\label{T:SSKLRBasis}
    Suppose that $e>n$ and $\Lambda\in P^+$ with $(\Lambda,\alpha_{i,n})\le1$, for all
    $i\in I$. Then $\R$ is a graded cellular algebra with graded cellular basis
    $\set{e_{\s\t}|(\s,\t)\in\Std^2(\Parts)}$
    with $\deg e_{\s\t}=0$ for all $(\s,\t)\in\Std^2(\Parts)$.
  \end{Theorem}

  \begin{proof}
    By \autoref{P:SSimpleIdempotents}, $y_r=0$ for $1\le r\le n$ and $e(\bi)=0$
    if $\bi\notin I^\Lambda_n$. In particular, this implies that
    $\psi_1,\dots,\psi_{n-1}$ satisfy the braid relations for the symmetric
    group~$\Sym_n$ because, by \autoref{L:SSResidues}, if $\bi\in I^n_\Lambda$ then
    $(i,i\pm 1,i)$ is not a subsequence of $\bi$, for any $i\in I$. Therefore,
    $\R$ is spanned by the elements $\psi_v e(\bi)\psi_w$, where $v,w\in\Sym_n$
    and $\bi\in I^n_\Lambda$. Moreover, if $\bj\in I^n$ then $e(\bj)\psi_v
    e(\bi)\psi_w=0$ unless $\bj=v\cdot\bi\in I^n_\Lambda$. Therefore, $\R$ is
    spanned by the elements $\set{e_{\s\t}|(\s,\t)\in\Std^2(\Parts)}$ as
    required by the statement of the theorem. Hence, $\R$ has rank at most
    $\sum_{\blam\in\Parts}|\Std(\blam)|^2=\ell^nn!$, where this combinatorial
    identity comes from \autoref{T:SeminormalBasis}.

    Let $K$ be the algebraic closure of the field of fractions
    of~$\Zcal$. Then $\R(K)\cong\R(\Zcal)\otimes_\Zcal K$.  By the last
    paragraph, the dimension of $\R$ is at most $\ell^nn!$.  Let
    $\rad\R(K)$ be the Jacobson radical of~$\R(K)$.  For each
    multipartition $\blam\in\Parts$, \autoref{P:KLRSemisimpleReps}
    constructs an irreducible graded Specht module~$S^\blam$. By
    \autoref{L:UniqueResidues}, if $\blam,\bmu\in\Parts$ and $d\in\Z$
    then $S^\blam\cong S^\bmu\<d\>$ if and only if~$\blam=\bmu$ and
    $d=0$. By the Wedderburn theorem,
    \begin{align*}
      \ell^nn!&\ge\dim \R(K)/\rad\R(K)\ge\sum_{\blam\in\Parts} (\dim S^\blam)^2\\
              &=\sum_{\blam\in\Parts} |\Std(\blam)|^2 = \ell^nn!.
    \end{align*}
    Hence, we have equality throughout, so
    $\set{e_{\s\t}|(\s,\t)\in\Std^2(\Parts)}$ is a basis of~$\R(K)$. As  the
    elements $\{e_{\s\t}\}$ span $\R(\Zcal)$, and their images in~$\R(K)$ are
    linearly independent, so $\{e_{\s\t}\}$ is a basis of~$\R(\Zcal)$.

    It remains to prove that $\{e_{\s\t}\}$ is a graded cellular basis of~$\R$.
    The orthogonality of the KLR idempotents implies that
    $e_{\s\t}e_{\u\v}=\delta_{\t\u}e_{\s\v}$. Therefore, $\{e_{\s\t}\}$ is a
    basis of matrix units for $\R$. Consequently, $\R$ is a direct sum of matrix
    rings, for any integral domain~$\Zcal$, and $\{e_{\s\t}\}$ is a cellular basis
    of~$\R$.

    Finally, we need to show that $e_{\s\t}$ is homogeneous of
    degree zero. This will follow if we show that $\deg\psi_r e(\bi)=0$, for
    $1\le r<n$ and $\bi\in I^n_\Lambda$. In fact, this is already clear because if
    $\bi\in I^n_\Lambda$ then $i_r\ne i_{r+1}$, by
    \autoref{L:SSResidues}, and if $i_{r+1}=i_r\pm1$ then $\psi_re(\bi)=0$ by
    \autoref{C:SameRowColumn} and \autoref{P:SSimpleIdempotents}.
  \end{proof}

  By definition, $e_{\s\t}e_{\u\v}=\delta_{\t\v}e_{\s\v}$. Let $\Mat_d(\Zcal)$
  be the ring of $d\times d$ matrices over~$\Zcal$. Hence, the proof of
  \autoref{T:SSKLRBasis} also yields the following.

  \begin{Corollary}\label{C:KLRMatrixRing}
    Suppose that $\Zcal$ is an integral domain $e>n$ and that
    $\Lambda\in P^+$ with $(\Lambda,\alpha_{i,n})\le1$, for all $i\in I$. Then
  $$\R(\Zcal)\cong\bigoplus_{\blam\in\Parts}\Mat_{s_\blam}(\Zcal),$$
  where $s_\blam=\#\Std(\blam)$ for $\blam\in\Parts$.
  \end{Corollary}

  Another consequence of \autoref{T:SSKLRBasis} is that the KLR relations
  simplify in the semisimple case --- giving a non-standard
  presentation for a direct sum of matrix rings.

  \begin{Corollary}\label{C:SSKLRRelations}
    Suppose that $\Zcal$ is an integral domain, $e>n$ and that $\Lambda\in P^+$ with
    $(\Lambda,\alpha_{i,n})\le1$, for all $i\in I$. Then $\R$ is the unital
    associative $\Z$-graded algebra generated by $\psi_1,\dots,\psi_{n-1}$ and
    $e(\bi)$, for $\bi\in I^n$, subject to the relations
    {\setlength{\abovedisplayskip}{2pt}
     \setlength{\belowdisplayskip}{1pt}
      \begin{align*}
        e(\bi)^{(\Lambda,\alpha_{i_1})}&=0  &
        {\textstyle\sum_{\bi \in I^n}} e(\bi)&= 1, &
        e(\bi) e(\bj) &= \delta_{\bi\bj} e(\bi),\\
        \psi_re(\bi)&=e(s_r\cdot\bi)\psi_r  &
        e(\bi)&=0\text{ if }i_r=i_{r+1}, &
        \psi_r^2e(\bi)&=e(\bi)
      \end{align*}
      \begin{align*}
        \psi_r\psi_s&=\psi_s\psi_r, & \text{if }|r-s|>1,
      \end{align*}
      \begin{align*}
        \psi_r \psi_{r+1} \psi_re(\bi) &= \begin{cases}
           (\psi_{r+1} \psi_{r} \psi_{r+1}-1)e(\bi),
                 &\text{if }i_{r+2}=i_r\rightarrow i_{r+1},\\
           (\psi_{r+1} \psi_{r} \psi_{r+1}+1)e(\bi),
                 &\text{if }i_{r+2}=i_r\leftarrow i_{r+1},\\
           \psi_{r+1} \psi_{r} \psi_{r+1}e(\bi),
                 &\text{otherwise,}
      \end{cases}
      \end{align*}
    }%
    for all $\bi,\bj\in I^n$ and admissible $r$ and~$s$. Moreover, $\R$
    is concentrated in degree zero.
  \end{Corollary}

  The reader is encouraged to check the details here.  Note that these
  relations, together with the argument of
  \autoref{P:SSimpleIdempotents}, imply that $e(\bi)\ne0$ if only if
  $\bi\in I_\Lambda^n$. In particular, the combinatorics of tableau
  content sequences is partially encoded in the failure of the braid
  relations for the $\psi_r$.

  As a final application, we prove Brundan and
  Kleshchev's graded isomorphism theorem in this special case.

  \begin{Corollary}\label{C:SSBKIso}
    Suppose that $\Zcal=K$ is a field,  $e>n$, and that $\Lambda\in P^+$
    with $(\Lambda,\alpha_{i,n})\le1$, for all $i\in I$. Then $\R\cong\H$.
  \end{Corollary}

  \begin{proof}
    By \autoref{C:SSKLRRelations} and
    \autoref{T:SeminormalBasis}, there is a well-defined homomorphism
    $\Theta\map{\R}{\H}$ determined by
    $e(\bi^\s)\mapsto F_\s$ and
    $$\psi_re(\bi^\s)\mapsto\begin{cases}
        \frac1{\alpha_r(\s)}\Big(T_r
          -\frac{c^\Zcal_{r+1}(\s)-c^\Zcal_r(\s)}{1+(v-v^{-1})c^\Zcal_{r+1}(\s)}\Big)F_\s,
             &\text{if }\alpha_r(\s)\ne0,\\
          0,&\text{otherwise,}
        \end{cases}
    $$
    for $\s\in\Std(\Parts)$ and $1\le r<n$. Using
    \autoref{T:SSKLRBasis}, or \autoref{P:KLRSemisimpleReps}, it follows
    that $\Theta$ is an isomorphism.
  \end{proof}

  We emphasize that it is essential to work over a field in \autoref{C:SSBKIso}
  because \autoref{C:KLRMatrixRing} says that $\R$ is always a direct sum of
  matrix rings whereas if $n>1$ this is only true of~$\H$ when it is defined
  over a field.

  These results suggest that~$\R$ should be considered as the ``idempotent
  completion'' of the algebra~$\H$ obtained by adjoining idempotents $e(\bi)$,
  for $\bi\in I^n$. We will see how to make sense of the idempotents
  $e(\bi)\in\H$ for any $\bi\in I^n$ in \autoref{T:BKiso} and
  \autoref{L:ResidueIdempotents} below.


  \subsection{The nil-Hecke algebra}\label{S:NilHecke} Still working just with
  the relations we now consider the shadow of the nil-Hecke algebra in the
  cyclotomic KLR setting. For the affine KLR algebras the nil-Hecke algebras
  case has been well-studied~\cite{KhovLaud:diagI,Rouq:2KM}. For the cyclotomic
  quotients (in type~$A$) the story is similar.

  For this section fix $i\in I$ and set $\beta=n\alpha_i$ and
  $\Lambda=n\Lambda_i$.  Following \autoref{E:RnBlocks}, set $\R[\beta]=e(\bi)\R
  e(\bi)$, where $\bi=\bi^\beta=(i^n)$. Then $\R[\beta]$ is a direct summand
  of~$\R$ and, moreover, it is a non-unital subalgebra with identity element
  $e(\bi)$. As $e(\bi)$ is the unique non-zero KLR idempotent in
  $\R[\beta]$, $\psi_r=\psi_r e(\bi)$ and $y_s=y_se(\bi)$.  Therefore,
  $\R[\beta]$ is the unital associative graded algebra generated by
  $\psi_r$ and $y_s$, for $1\le r<n$ and $1\le s\le n$, with relations
  \begin{gather*}
      y_1^n=0,  \qquad \psi_r^2=0,  \qquad y_r y_s = y_s y_r,\\
      \psi_r y_{r+1} =y_r\psi_r+1, \qquad y_{r+1}\psi_r=\psi_r y_r+1,\\
      \psi_r \psi_s  = \psi_s \psi_r \quad\text{if }|r-s|>1,\qquad
      \psi_r y_s  = y_s \psi_r \quad\text{if }s \neq r,r+1,\\
  \psi_{r}\psi_{r+1} \psi_{r}  = \psi_{r+1} \psi_{r} \psi_{r+1}.
  \end{gather*}
  The grading on $\R[\beta]$ is determined by $\deg\psi_r=-2$ and $\deg y_s=2$.
  Some readers will recognize this presentation as defining as a cyclotomic
  quotient of the \textbf{nil-Hecke algebra} of
  type~$A$~\cite{KostantKumar:NilHecke}. Note that the argument from Case~3 of
  \autoref{L:nilpotence} shows that $y_r^\ell=0$ for $1\le r\le\ell$.

  Let $\blam=(1|1|\dots|1)\in\Parts[\beta]$. Then the map $\t\mapsto d(\t)$
  defines a bijection between the set of standard $\blam$-tableaux and the
  symmetric group~$\Sym_n$. For convenience, we identify the standard
  $\blam$-tableaux with the set of (non-standard) tableaux of partition
  shape~$(n)$ by concatenating their components. In other words, if $d=d(\t)$
  then $\t=\otab(d_1,d_2,\Dots,d_n)$, where $d=d_1\dots d_n$ is the
  permutation written in one-line notation.

  If $\v,\s\in\Std(\blam)$ then write $\s\GDom\v$ if $\s\gdom\v$ and
  $\ell(d(\v))=\ell(d(\s))+1$. To make this more explicit write $t\prec_\v m$ if
  $t$ is in an earlier component of~$\v$ than~$m$ --- that is, $t$ is to the left
  of~$m$ in~$\v$. The reader can check that $\s\GDom\v$ if and only if
  there exist integers $1\le m<t\le n$ such that $\s=\v(m,t)$, $m\prec_\v t$ and
  if $m<l<t$ then either $l\prec_\v m$ or $t\prec_\v l$.

  \begin{Example}
    Suppose that $n=6$. Let $\v=\Otab[4/0](4,6,5,3,1,2)$ and take $t=3$. Then
    $$ \Big\{\hspace*{1mm}
                     \Otab[4/0,1/0](3,6,5,4,1,2), \hspace*{1.5mm}
                     \Otab[4/0,3/0](4,6,3,5,1,2), \hspace*{1.5mm}
                     \Otab[4/0,6/0](4,6,5,2,1,3), \hspace*{1.5mm}
                     \Otab[4/0,5/0](4,6,5,1,3,2)
       \hspace*{1mm}\Big\}$$
       is the set of $\blam$-tableaux
       $\set{\s|\s=\v(3,r)\GDom\v\text{ for }1\le r\le n}$.
  \end{Example}

  We can now state the main result of the section.

  \begin{Proposition}\label{P:SchubertPolynomials}
    Suppose that $\beta=n\alpha_i$ and $\Lambda=n\Lambda_i$, for $i\in I$.  Then
    there is a unique graded $\R[\beta]$-module $S^\blam$ with homogeneous basis
    $\set{\psi_\s|\s\in\Std(\blam)}$ such that
    $\deg \psi_\s=\binom n2-2\ell(d(\s))$ and
    \begin{align*}
      \psi_\s\psi_r&=\begin{cases}
          \psi_{\s(r,r+1)},& \text{if }\s\gdom\s(r,r+1)\in\Std(\blam),\\
          0,&\text{otherwise},
        \end{cases}\\
        \psi_\v y_t&=
            \sum_{\substack{1\le k<t\\\u=\v(k,t)\GDom\u}} \psi_\u
            \quad-\quad \sum_{\substack{t<k\le n\\\u=\v(k,t)\GDom\v}} \psi_\u,
    \end{align*}
    for $\s,\v\in\Std(\blam)$, $1\le r<n$ and $1\le t\le n$.
    Moreover, if $\Zcal$ is a field then $S^\blam$ is irreducible.
  \end{Proposition}

  \begin{proof} The uniqueness is clear. To show that $S^\blam$ is an
    $\R[\beta]$-module we check that the action respects the
    relations of~$\R[\beta]$. By definition, if $\v\in\Std(\blam)$ then
    $\psi_\v=\psi_{\tlam}\psi_{d(\v)}$ and $\psi_\v\psi_r^2=0$ since
    $\psi_\v\psi_r=0$ if $\v(r,r+1)\gdom\v$. In particular, this implies that
    the action of $\psi_1\dots,\psi_{n-1}$ on~$S^\blam$ respects the braid
    relations of~$\Sym_n$ and that~$\psi_\v$ has the specified degree.
    Further, note that if $\u\GDom\v$ then $\ell(d(\v))=\ell(d(\u))+1$ so
    that $\deg \psi_\u=\deg \psi_\v+2$.

    By the last paragraph, the action of $\R[\beta]$ is compatible with the
    grading on~$S^\blam$, but we still need to check the
    relations involving $y_1,\dots,y_n$. First consider $\psi_\v y_ry_t=\psi_\v
    y_ty_r$, for $1\le r,t\le n$ and $\v\in\Std(\blam)$. If $r=t$ there is
    nothing to prove so suppose $r\ne t$.  By definition,
    $$\psi_\v y_ty_r = \sum_{\u\GDom\v}\sum_{\s\GDom\u}
         \varepsilon_t(\v,\u)\varepsilon_r(\u,\s) \psi_\s,$$
    for appropriate choices of the signs $\varepsilon_t(\v,\u)$
    and $\varepsilon_r(\u,\s)$. Suppose that $\psi_\s$ appears with non-zero
    coefficient in this sum. Then we can write $\u=\v(m,t)$ and
    $\s=\v(l,r)$, for some $l,m$ such that $\s\GDom\u\GDom\v$. Suppose first that $l\ne m$.
    Then the permutations $(m,t)$ and $(l,r)$ commute and, as their lengths add, we
    have $\s\GDom\v(l,r)\GDom\v$. Therefore, $\psi_\s$ appears with the
    same coefficient in $\psi_\v y_t y_r$ and~$\psi_\v y_r y_t$. If $l=m$ then
    $\s\gdom\u\gedom\v$ only if $m$ is in between $r$
    and~$t$ in~$\v$. That is, either $r\prec_\v m\prec_\v t$ or
    $t\prec_\v m\prec_\v r$. However, this implies that either $\s\notGDom\u$ or
    $\u\notGDom\v$, so that~$\psi_\s$ does not appear in either $\psi_\v y_ry_t$
    or in~$\psi_\v y_ty_r$. Hence, the actions of $y_r$ and $y_t$ on~$S^\blam$
    commute.

    Similar, but easier, calculations with tableaux show that the action
    defined on~$S^\blam$ respects the three relations
    $\psi_r y_{r+1} =y_r\psi_r+1$, $y_{r+1}\psi_r=\psi_r y_r+1$ and
    $\psi_r y_t=y_t\psi_r$ when $t\ne r,r+1$.  To complete the
    verification of the relations in~$\R[\beta]$ it remains to show that
    $\psi_\v y_1^n=0$, for all $\v\in\Std(\blam)$. This is clear,
    however, because $\psi_\v y_1$ is equal to a linear combination of
    terms~$\psi_\s$ where~$1$ appears in an earlier component of~$\s$
    than it does in~$\v$.

    Finally, it remains to prove that $S^\blam$ is irreducible over a field.
    First we need some more notation. Let $\t_\blam=\otab(n,\Dots,\Dots,2,1)$
    and set $w_\blam=d(\t_\blam)$.  Then $w_\blam$ is the unique element of
    maximal length in~$\Sym_n$. Recall from \autoref{S:Tableaux}, that $d'(\t)$
    is the unique permutation such that $\t=\t_\blam d'(\t)$ and, moreover,
    $d(\t)d'(\t)^{-1}=w_\blam$ with the lengths adding.  Therefore, if
    $\ell(d(\s))\ge\ell(d(\t))$ then
    $\psi_\s\psi_{d'(\t)}^\star=\delta_{\s\t}\psi_{\t_\blam}$.

    We are now ready to show that $S^\blam$ is irreducible.  Suppose that
    $x=\sum_\s r_\s \psi_\s$ is a non-zero element of~$S^\blam$. Let~$\t$ be any
    tableau such that $r_\t\ne0$ and $\ell(d(\t))$ is minimal.
    Then, by the last paragraph, $x\psi_{d'(\t)}^\star=r_\t \psi_{\t_\blam}$, so
    $\psi_{\t_\blam}\in x\R[\beta]$. We have already observed that $y_1$ acts by
    moving $1$ to an earlier component. Therefore,
    $\psi_{\t_\blam}y_1^{n-1}=(-1)^{n-1}\psi_{\t_{\blam,1}}$, where
    $\t_{\blam,1}=\otab(1,n,\Dots,3,2)$. Similarly,
    $\psi_{\t_\blam}y_1^{n-1}y_2^{n-2}=(-1)^{2n-3}\psi_{\t_{\blam,2}}$, where
    $\t_{\blam,2}=\otab(1,2,n,\Dots,3)$. Continuing in this way shows that
    $\psi_{\t_\blam}y_1^{n-1}y_2^{n-2}\dots y_{n-1}
              =(-1)^{\frac12n(n-1)}\psi_{\tlam}$.
    Hence, $x\R[\beta]=S^\blam$, so that $S^\blam$ is irreducible as claimed.
  \end{proof}

  The proof of \autoref{P:SchubertPolynomials} shows that
  $y_1^{n-1}y_2^{n-2}\dots y_{n-1}$ is a non-zero element of~$\R[\beta]$.
  Using the relations, and a bit of ingenuity, it is possible to show that
  $\set{\psi_wy_1^{a_1}\dots y_n^{a_n}|w\in\Sym_n\text{ and }0\le a_r\le n-r,
        \text{ for }1\le r\le n}$
  is a basis of~$\R[\beta]$. Alternatively, it follows from
  \cite[Theorem~4.20]{BK:GradedDecomp} that $\dim\R[\beta]=(n!)^2$.  Hence,
  we obtain the following.

  \begin{Corollary}\label{C:NilHeckeBasis}
    Suppose that $\beta=n\alpha_i$ and $\Lambda=n\Lambda_i$, for $i\in I$.
    Let $\blam=(1|1|\dots|1)$ and for $\s,\t\in\Std(\blam)$ define
    $\psi_{\s\t}=\psi_{d(\s)}^\star e(\ilam)y^\blam\psi_{d(\t)}$, where
    $\ilam=\bi^{\tlam}$ and $y^\blam=y_1^{n-1}y_2^{n-2}\dots y_{n-1}$.
    Then $\set{\psi_{\s\t}|\s,\t\in\Std(\blam)}$ is a graded cellular basis of
    $\R[\beta]$.
  \end{Corollary}

  The basis of the Specht module $S^\blam$ in \autoref{P:SchubertPolynomials} is
  well-known because it is really a disguised version of the basis of
  \textit{Schubert polynomials} of the coinvariant algebra of the symmetric
  group~$\Sym_n$~\cite{LascouxSch:SchubertPolys,Manivel:SchurbertPolys}.  The
  \textbf{coinvariant algebra}~$\Poly$ is the quotient of the polynomial ring
  $\Z[\bx]=\Z[x_1,\dots,x_n]$ by the deal generated by the symmetric
  polynomials in $x_1,\dots,x_n$ of \textit{positive degree}.  Then
  $\Poly$ is free of rank~$n!$. As we have quotiented out by a
  homogeneous ideal, $\Poly$ inherits a grading from~$\Z[\bx]$, where we
  set $\deg x_r=2$ for $1\le r\le n$. Identify $x_r$ with its image in
  $\Poly$, for $1\le r\le n$. There is a well-defined action
  of~$\R[\beta]$ on~$\Poly$ where~$y_r$ acts as multiplication by~$x_r$,
  and~$\psi_r$ acts from the right as a \textit{divided difference
  operator}:
  $$f(\bx)\psi_r= \partial_rf(\bx)
        = \frac{ f(\bx)-f(s_r\cdot\bx)}{x_r-x_{r+1}},$$
  where
  $\bx=(x_1,\dots,x_n)$ and $s_r\cdot\bx=(x_1,\dots,x_{r+1},x_r,\dots,x_n)$
  for $1\le r<n$. Here we are secretly thinking of $\R[\beta]$ as being a
  quotient of the nil-Hecke algebra, where this action is well-known.

  For $d\in\Sym_n$ define $\sigma_d=(x_1^{n-1}x_2^{n-2}\dots x_{n-1})\psi_{w_0d}$.
  Then $\set{\sigma_d|d\in\Sym_n}$ is the basis of \textbf{Schubert polynomials}
  of~$\Poly$.  The Specht module is isomorphic to $\Poly$ as an $\R[\beta]$-module,
  where an isomorphism is given by $\psi_\t\mapsto\sigma_{d'(\t)}$.  To see
  this it is enough to know that the Schubert polynomials satisfy the identity
  $$\partial_r \sigma_d = \begin{cases}
                \sigma_{s_rd}, & \text{if }\ell(s_rd)=\ell(d)-1,\\
                0,&\text{otherwise}.
              \end{cases}
  $$
  By the last paragraph of the proof of \autoref{P:SchubertPolynomials}, if
  $\t\in\Sym_n$ then
  $$\psi_\t=\psi_{\tlam}\psi_{d(\t)}
           =\psi_{\t_\blam}y_1^{n-1}y_2^{n-1}\dots y_{n-1}\psi_{d(\t)}.$$
  Therefore, our claim follows by identifying $\psi_{\t_\blam}$ with the polynomial
  $1\in\Poly$.

  Finally, we remark that the formula for the action of $y_1,\dots,y_n$ in
  \autoref{P:SchubertPolynomials} is a well-known corollary of Monk's rule;
  see, for example, \cite[Exercise~2.7.3]{Manivel:SchurbertPolys}.

\section{Isomorphisms, Specht modules and categorification}
  \label{S:GradedIso}

  In the last section we proved that the algebras $\R$ and $\H$ are isomorphic
  when $e>n$ and $(\Lambda,\alpha_{i,n})\le1$, for all $i\in I$. In this
  section we state Brundan and Kleshchev's Graded Isomorphism Theorem,
  which says that $\R\cong\H$, and we start to investigate the
  consequences of this result for both algebras.

  \subsection{The Graded Isomorphism Theorem}
  One of the most fundamental results for the cyclotomic Hecke algebras~$\H$ is
  Brundan and Kleshchev's spectacular isomorphism theorem.

  \begin{Theorem}[Graded Isomorphism Theorem~\cite{BK:GradedKL,Rouq:2KM}]
  \label{T:BKiso}
  Suppose that $\Zcal=F$ is a field, $v\in F$ has quantum characteristic~$e$ and
  that $\Lambda\in P^+$. Then there is an isomorphism of algebras $\R\cong\H$.
  \end{Theorem}

  Suppose that $F$ is a field of characteristic $p>0$ and that $e=pf$, where
  $f>1$. Then $F$ cannot contain an element $v$ of quantum characteristic~$e$,
  so \autoref{T:BKiso} says nothing about the quiver Hecke algebra $\R(F)$.

  As a first consequence of \autoref{T:BKiso}, by identifying $\H$ and $\R$ we
  can consider $\H$ as a graded algebra.

  \begin{Corollary}
    Suppose that $\Lambda\in P^+$ and $\Zcal=F$ is a field. Then there is a unique
    grading on~$\H$ such that $\deg e(\bi)=0$, $\deg y_r=2$ and
    $\deg\psi_se(\bi)=-c_{i_s,i_{s+1}}$, for $1\le r\le n$, $1\le s<n$
    and $\bi\in I^n$.
  \end{Corollary}

  Brundan and Kleshchev prove \autoref{T:BKiso} by constructing  family of
  isomorphisms $\R\longrightarrow\H$, together with their inverses, and then
  painstakingly checking that these isomorphisms respect the relations of both
  algebras. Their argument starts with the well-known fact that $\H$ decomposes
  into a direct sum of simultaneous generalized eigenspaces for the Jucys-Murphy
  elements $L_1,\dots,L_n$. These eigenspaces are indexed by~$I^n$, so for each
  $\bi\in I^n$ there is an element $e(\bi)\in\H$, possibly zero, such that
  $e(\bi)e(\bj)=\delta_{\bi\bj}e(\bi)$.  We describe these idempotents
  explicitly in \autoref{L:ResidueIdempotents} below.

  Translating through \autoref{D:HeckeAlgebras}, Brundan and Kleshchev's
  isomorphism is given by $e(\bi)\mapsto e(\bi)$ and
  $$ y_r\mapsto \sum_{\bi\in I^n}v^{-i_r}\(L_r -[i_r]_v\)e(\bi),\quad\text{and}\quad
     \psi_s\mapsto \sum_{\bi\in I^n}\(T_s+P_s(\bi)\)\frac1{Q_s(\bi)}e(\bi),
  $$
  for $1\le r\le n$, $1\le s<n$, $\bi\in I^n$ and where $P_r(\bi)$
  and $Q_r(\bi)$ are certain rational functions in $y_r$ and $y_{r+1}$
  that are well-defined because $(L_t-[i_t]_v)e(\bi)$ is nilpotent
  in~$\H$, for $1\le t\le n$; see \cite[\S3.3 and \S4.3]{BK:GradedKL}.
  We are abusing notation by identifying the KLR generators with their
  images in~$\H$. The inverse isomorphism is given by $e(\bi)\mapsto
  e(\bi)$,
  $$ L_r\mapsto\sum_{\bi\in I^n}\(v^{i_r}y_r+[i_r]_v\)e(\bi)\quad\text{and}\quad
     T_s\mapsto \sum_{\bi\in I^n}\(\psi_sQ_s(\bi)-P_s(\bi)\)e(\bi),
  $$
  for $1\le r\le n$, $1\le s<n$ and $\bi\in I^n$.

  Rouquier \cite[Corollary~3.20]{Rouq:2KM} has given a more direct proof of
  \autoref{T:BKiso} by first showing that the (non-cyclotomic) quiver Hecke
  algebra~$\Rn$ is isomorphic to the (extended) affine Hecke algebra of
  type~$A$.  Following~\cite{HuMathas:SeminormalQuiver}, we sketch another
  approach to \autoref{T:BKiso} in \autoref{S:OKLR} below.

  The following easy but important application of \autoref{T:BKiso} was a
  surprise (at least to the author!).

  \begin{Corollary}\label{C:independence}
    Suppose that $\Zcal=F$ is a field and that $v,v'\in F$ are
    two elements of quantum characteristic~$e$.
    Then $\Hcal^\Lambda_n(F,v)\cong\Hcal^\Lambda_n(F,v')$.
  \end{Corollary}

  \begin{proof}
    By \autoref{T:BKiso}, $\Hcal^\Lambda_n(F,v)\cong\R(F)\cong\Hcal^\Lambda_n(F,v')$.
  \end{proof}

  Consequently, up to isomorphism, the algebra $\H$ depends only
  on~$e$,~$\Lambda$ and the field~$F$. Therefore, because $\H$ is cellular, the
  decomposition matrices of~$\H$ depend only on~$e$,~$\Lambda$ and~$p$,
  where~$p$ is the characteristic of~$F$. In the special case of the symmetric
  group, when $\Lambda=\Lambda_0$, this weaker statement for the decomposition
  matrices was conjectured in \cite[Conjecture~6.38]{M:Ulect}.

  When $F=\C$ it is easy to prove \autoref{C:independence} because there is a
  Galois automorphism of~$\Q(v)$, as an extension of~$\Q$, which
  interchanges~$v$ and~$v'$. It is not difficult to see that this automorphism
  induces an isomorphism $\Hcal^\Lambda_n(F,v)\cong\Hcal^\Lambda_n(F,v')$. This
  argument fails for fields of positive characteristic because such fields have
  fewer automorphisms.

  \subsection{Graded Specht modules}\label{S:GradedSpecht} As we noted in
  \autoref{S:GradedCellular}, if we impose a grading on an algebra
  $\underline{A}$ then it is not true that every (ungraded)
  $\underline{A}$-module has a graded lift, so there is no reason to expect that
  graded lifts of Specht modules $\UnS^\lambda$ always exist. Of course,
  graded Specht modules do exist and this section describes one way to define
  them.

  Recall from \autoref{S:Murphy} that the ungraded Specht module $\UnS^\blam$,
  for $\blam\in\Parts$, has basis $\set{m_\t|\t\in\Std(\blam)}$. By
  construction, $\UnS^\blam=m_{\tlam}\H$.  Brundan, Kleshchev and
  Wang~\cite{BKW:GradedSpecht} proved that $\UnS^\blam$ has a graded lift
  essentially by declaring that $m_{\tlam}$ should be homogeneous and then
  showing that this induces a grading on the Specht
  module~$\UnS^\blam=m_{\tlam}\R$.

  Partly inspired by \cite{BKW:GradedSpecht}, Jun Hu and the
  author~\cite{HuMathas:GradedCellular} showed that $\H$ is a graded cellular
  algebra. The graded cell modules constructed from this cellular basis coincide
  exactly with those of~\cite{BKW:GradedSpecht}. Perhaps most significantly, the
  construction of the graded Specht modules using cellular algebra techniques
  endows the graded Specht modules with a homogeneous bilinear form of degree
  zero.

  Following Brundan, Kleshchev and Wang~\cite[\S3.5]{BKW:GradedSpecht}
  we now define the degree of a standard tableau. Suppose that $\bmu\in\Parts$.
  For $i\in I$ let $\Add_i(\bmu)$ be the set of addable $i$-nodes of~$\bmu$ and
  let $\Rem_i(\bmu)$ be its set of removable $i$-nodes.

  \begin{Definition}\label{E:dA}
  If~$A$ is an addable or
  removable $i$-node of~$\bmu$ define:
    \begin{align*}
      d^A(\bmu)&=\#\set{B\in\Add_i(\bmu)|A>B}-\#\set{B\in\Rem_i(\bmu)|A>B},\\
      d_A(\bmu)&=\#\set{B\in\Add_i(\bmu)|A<B}-\#\set{B\in\Rem_i(\bmu)|A<B},\\
      d_i(\bmu)&=\#\Add_i(\bmu)-\#\Rem_i(\bmu).
    \end{align*}
  \end{Definition}

  If $\t$ is a standard $\bmu$-tableau then its
  \textbf{codegree} and \textbf{degree} are defined inductively by setting
  $\codeg_e\t=0=\deg_e\t$ if $n=0$ and if~$n>0$ then
  $$ \codeg_e\t=\codeg_e\t_{\downarrow(n-1)}+d^A(\bmu)
  \quad\text{and}\quad
     \deg_e\t=\deg_e\t_{\downarrow(n-1)}+d_A(\bmu),
  $$
  where $A=\t^{-1}(n)$. If $e$ is fixed write
  $\codeg\t=\codeg_e\t$ and $\deg\t=\deg_e\t$.

  Implicitly, all of these definitions depend on the choice of
  multicharge~$\charge$. The definition of the degree and codegree of a
  standard tableau is due to Brundan, Kleshchev and
  Wang~\cite{BKW:GradedSpecht}, however, the underlying combinatorics
  dates back to Misra and Miwa~\cite{MisraMiwa} and their work on the
  crystal graph and Fock space representations of~$\Uq$.

  Recall that we have fixed an arbitrary  reduced expression for each
  permutation $w\in\Sym_n$.  In \autoref{S:Tableaux} for each standard
  tableau $\t\in\Std(\blam)$ we have defined permutations
  $d'(\t),d(\t)\in\Sym_n$ by $\t_\blam d'(\t)=\t=\tlam d(\t)$.

  \begin{Definition}[\protect{%
    \cite[Definitions~4.9 and 5.1]{HuMathas:GradedCellular}}]\label{D:psis}
  Suppose that $\bmu\in\Parts$. Define non-negative integers
  $d^\bmu_1,\dots,d^\bmu_n$ and
  $d_\bmu^1,\dots,d_\bmu^n$ recursively by requiring that
  $d_\bmu^1+\dots+d_\bmu^k=\codeg(\tmu_{\downarrow k})$ and
  $d^\bmu_1+\dots+d^\bmu_k=\deg(\tmu_{\downarrow k})$,
  for $1\le k\le n$. Now set $\bi_\bmu=\bi^{\t_\bmu}$, $\bi^\bmu=\bi^{\tmu}$,
  $y_\bmu=y_1^{d_\bmu^1}\dots y_n^{d_\bmu^n}$ and
  $y^\bmu=y_1^{d^\bmu_1}\dots y_n^{d^\bmu_n}$.
  For  $(\s,\t)\in\Std^2(\bmu)$ define
  $$\psi'_{\s\t}=\psi_{d'(\s)}^\star e(\bi_\bmu) y_\bmu\psi_{d'(\t)}
  \quad\text{and}\quad
    \psi_{\s\t}=\psi_{d(\s)}^\star e(\imu) y^\bmu\psi_{d(\t)},$$
  where $\star$ is the unique (homogeneous) anti-isomorphism of $\R$ that fixes
  the KLR generators.
  \end{Definition}

  \begin{Example}
    Suppose that $e=3$, $\Lambda=\Lambda_0+\Lambda_2$ and $\bmu=(7,6,3,2|4,3,1)$,
    with multicharge $\charge=(0,2)$. Then
    $$\tmu=\left(\rule{0cm}{9.8mm}\right.\space
      \begin{tikzpicture}[baseline=-7.5mm,scale=0.44,draw/.append style={thick,black}]
        \foreach \r/\c in {2/0,2/-1,2/-2,5/0,5/-1,10/0,10/-1} {
          \filldraw[pattern=my north west lines,pattern color=blue!55,
                    line space=4pt, pattern line width=1.5pt,opacity=0.8]
              (\r,\c)+(-.5,-.5)rectangle++(.5,.5);
        }
        \foreach \r/\c in {0/-1,1/-2,2/0,3/-1,5/0} {
          \filldraw[pattern=my north east lines,pattern color=red!75,
                    line space=4pt, pattern line width=1.5pt,opacity=0.8]
              (\r,\c)+(-.5,-.5)rectangle++(.5,.5);
        }
        \newcount\row
        \newcount\col
        \foreach\tab/\Col in %
              {{{1,2,3,4,5,6,7},{8,9,10,11,12,13},{14,15,16},{17,18}}/0,%
               {{19,20,21,22},{23,24,25},{26}}/8} {
           \row=0
           \foreach\tabRow in \tab {
             \global\col=\Col
             \foreach\k in \tabRow {
                \draw(\the\col,\row)+(-.5,-.5)rectangle++(.5,.5);
                \draw(\the\col,\row)node{\k};
                \global\advance\col by 1
             }
             \global\advance\row by -1
           }
        }
        \draw[thin](7,0.5)--+(0,-4);
      \end{tikzpicture}%
      \space\left.\rule{0cm}{9.8mm}\right)
    $$
    The reader may check that
    $e(\bi^\bmu)=e(01201202012011200120121200)$.  We have shaded the
    nodes in~$\tmu$ when they have column index divisible by~$e$ and
    when they have residue~$2=\res_{\tmu}(19)$. This should convince the
    reader that
    $y^\bmu=y_3^2y_6^2y_8y_{10}y_{11}y_{13}y_{15}y_{16}y_{21}y_{25}$.
    With analogous shadings,
    $$\t_\bmu=\left(\rule{0cm}{9.8mm}\right.\space
      \begin{tikzpicture}[baseline=-7.5mm,scale=0.44,draw/.append style={thick,black}]
        \foreach \r/\c in {0/-2,1/-2,2/-2,8/-2} {
          \filldraw[pattern=my north west lines,pattern color=blue!55,
                    line space=4pt, pattern line width=1.5pt,opacity=0.8]
              (\r,\c)+(-.5,-.5)rectangle++(.5,.5);
        }
        \foreach \r/\c in {8/-2,9/0,10/-1} {
          \filldraw[pattern=my north east lines,pattern color=red!75,
                    line space=4pt, pattern line width=1.5pt,opacity=0.8]
              (\r,\c)+(-.5,-.5)rectangle++(.5,.5);
        }
        \newcount\row
        \newcount\col
        \foreach\tab/\Col in %
              {{{9,13,17,20,22,24,26},{10,14,18,21,23,25},{11,15,19},{12,16}}/0,%
               {{1,4,6,8},{2,5,7},{3}}/8} {
           \row=0
           \foreach\tabRow in \tab {
             \global\col=\Col
             \foreach\k in \tabRow {
                \draw(\the\col,\row)+(-.5,-.5)rectangle++(.5,.5);
                \draw(\the\col,\row)node{\k};
                \global\advance\col by 1
             }
             \global\advance\row by -1
           }
        }
        \draw[thin](7,0.5)--+(0,-4);
      \end{tikzpicture}%
      \space\left.\rule{0cm}{9.5mm}\right).
      $$
      Hence, reading right to left, $y_\bmu=y_3^2y_4y_7y_{11}y_{15}y_{19}$.
      Note that $\res_{\t_\bmu}(9)=0$.
    \end{Example}

  \begin{Theorem}[\protect{Hu-Mathas\cite[Theorem~5.8]{HuMathas:GradedCellular}}]
    \label{T:PsiBasis}
    Suppose that $\Zcal=F$ is a field. Then
    $\set{\psi_{\s\t}|(\s,\t)\in\Std^2(\Parts)}$
    is a graded cellular basis of~$\R$ with $\psi_{\s\t}^\star=\psi_{\t\s}$ and
    $\deg\psi_{\s\t}=\deg\s+\deg\t$, for $(\s,\t)\in\Std^2(\Parts)$.
  \end{Theorem}

  \begin{Example}\label{Ex:NilHeckeylam}
      Let $\beta=n\alpha_i$ and $\Lambda=n\Lambda_i$, for some $i\in I$, so
      that~$\R[\beta]$ is the nil-Hecke algebra $\R[\beta]$ of \autoref{S:NilHecke}.
      Let $\blam=(1|1|\dots|1)$. Then the definitions give
      $y^\blam=y_1^{n-1}\dots y_{n-2}^2y_{n-1}$.
      Hence, the basis $\{\psi_{\s\t}\}$ of $\R[\beta]$ coincides with that of
      \autoref{C:NilHeckeBasis}.
  \end{Example}

  \begin{Example}
    As in \autoref{E:PsiReduced}, in general, the basis element~$\psi_{\s\t}$
    depends on the choices of reduced expressions that we have fixed for
    the permutations $d(\s)$ and $d(\t)$.  For example, let
    $\Lambda=2\Lambda_0+\Lambda_1$, $\charge=(0,1,0)$ and $\bmu=(1|1|1)$ and
    consider the standard $\bmu$-tableaux
    $\tmu=\onetab(1|2|3)$ and~$\t_\bmu=\onetab(3|2|1)$.
    Then $d(\tmu)=1$ and $d(\t_\bmu)=(1,3)=s_1s_2s_1=s_2s_1s_2$ has two different
    reduced expressions. Let
    $\psi_{\t_\bmu\t_\bmu}=\psi_1\psi_2\psi_1e(\bi^\bmu)y^\bmu\psi_1\psi_2\psi_1$
    and
    $\hat\psi_{\t_\bmu\t_\bmu}=\psi_2\psi_1\psi_2e(\bi^\bmu)y^\bmu\psi_2\psi_1\psi_2$.
    Then the calculation in \autoref{E:PsiReduced} implies that
    $$\hat\psi_{\t_\bmu\t_\bmu} =\psi_{\t_\bmu\t_\bmu}
             +\psi_{\t_\bmu\tmu}+\psi_{\tmu\t_\bmu}+\psi_{\tmu\tmu}.$$
    This is probably the simplest example where different reduced expressions
    lead to different $\psi$-basis elements, but examples occur for almost
    all~$\R$. This said, in view of \autoref{P:KLRSemisimpleReps}, $\psi_{\s\t}$
    is independent of the choice of reduced expressions for $d(\s)$ and $d(\t)$
    whenever $e>n$ and $(\Lambda,\alpha_{i,n})\le1$, for all $i\in I$.
    The $\psi$-basis can be independent of the choice of reduced
    expressions even when~$\R$ is not semisimple. For example, this is
    always the case when $e>n$ and $\ell=2$ by
    \cite[Appendix]{HuMathas:QuiverSchurI}, yet these algebras are typically
    not semisimple.
  \end{Example}

  Using the theory of graded cellular algebras from
  \autoref{S:GradedCellular}, \autoref{T:PsiBasis} allows us to
  construct a family $\set{S_F^\blam|\blam\in\Parts}$ of graded Specht
  modules for~$\H$.  By \cite[Corollary~5.10]{HuMathas:GradedCellular}
  the graded Specht modules attached to the $\psi$-basis coincide with
  those constructed by Brundan, Kleshchev and
  Wang~\cite{BKW:GradedSpecht}. When $e>n$ and
  $(\Lambda,\alpha_{i,n})\le1$, for $i\in I$, it is not hard to show
  that these Specht modules coincide with those we constructed in
  \autoref{P:KLRSemisimpleReps} above.  Similarly, for the nil-Hecke
  algebra considered in \autoref{S:NilHecke}, the graded Specht module
  $S_F^\blam$, with $\blam=(1|1|\dots|1)$, is isomorphic to the graded
  module constructed in \autoref{P:SchubertPolynomials}. Moreover, on
  forgetting the grading $S_F^\blam$ coincides exactly with the ungraded
  Specht module~$\UnS_F^\blam$ constructed in~\autoref{S:Murphy}, for
  $\blam\in\Parts$.

  If $\blam\in\Parts$ the graded Specht module $S_F^\blam$ has basis
  $\set{\psi_\t|\t\in\Std(\blam)}$, with $\deg\psi_\t=\deg\t$. The reader should
  be careful not to confuse $\psi_\t\in S_F^\blam$ with $\psi_{d(\t)}\in\R$!
  By \autoref{T:PsiBasis} we recover
  \cite[Theorem~4.20]{BK:GradedDecomp}:
  $$
     \gdim\H=\sum_{(\s,\t)\in\Std(\blam)}q^{\deg\s+\deg\t}
            =\sum_{\blam\in\Parts} \(\gdim S_F^\blam\)^2.$$

  In essence, \autoref{T:PsiBasis} is proved in much the same way that Brundan,
  Kleshchev and Wang~\cite{BKW:GradedSpecht} constructed a grading on the Specht
  modules: we proved that the transition matrix between the $\psi$-basis and the
  Murphy basis of \autoref{T:MurphyBasis} is triangular. In order to do this we
  needed the correct definition of the elements~$y^\bmu$, which we discovered by
  first looking at the one dimensional two-sided ideals of~$\H$ (which are
  necessarily homogeneous). We then used Brundan and Kleshchev's Graded
  Isomorphism \autoref{T:BKiso}, together with the seminormal forms
  (\autoref{T:SeminormalBasis}), to show that $e(\bi^\bmu)y^\bmu\ne0$. This
  established that  the basis of \autoref{T:PsiBasis} is a graded cellular
  basis. Finally, the combinatorial results of~\cite{BKW:GradedSpecht} are used
  to determine the degree of $\psi$-basis elements.

  Following the recipe in \autoref{S:GradedCellular}, for $\bmu\in\Parts$ define
  $D_F^\bmu=S_F^\bmu/\rad S_F^\bmu$, where $\rad S_F^\bmu$ is the radical of the
  homogeneous bilinear form on~$S_F^\bmu$. This yields the classification of the
  graded irreducible $\H$-modules. The main point of the next result is that the
  labelling of the graded irreducible $\H$-modules agrees with
  \autoref{C:UngradedSimples}.

  \begin{Corollary}[\!\protect{\cite[Theorem~5.13]{BK:GradedDecomp},
    \cite[Corollary~5.11]{HuMathas:GradedCellular}}]\label{C:GradedSimples}
    Suppose that $\Lambda\in P^+$ and that $\Zcal=F$ is a field. Then
    $\set{D_F^\bmu\<d\>|\bmu\in\Klesh\text{ and }d\in\Z}$
    is a complete set of pairwise non-isomorphic graded $\H$-modules. Moreover,
    $(D_F^\bmu)^\circledast\cong D_F^\bmu$ and $D_F^\bmu$ is absolutely irreducible,
    for all $\bmu\in\Klesh$.
  \end{Corollary}

  The \textbf{graded decomposition numbers} are the Laurent
  polynomials
  \begin{equation}\label{E:GradedDec}
   d^F_{\blam\bmu}(q)=[S_F^\blam:D_F^\bmu]_q
            =\sum_{d\in\Z}[S^\blam_F:D^\bmu_F\<d\>)]\,q^d,
  \end{equation}
  for $\blam\in\Parts$ and $\bmu\in\Klesh$.
  Write $S^\blam=S^\blam_F$, $D^\bmu=D^\bmu_F$
  and $d_{\blam\bmu}(q)=d^F_{\blam\bmu}(q)$ when~$F$ is understood. By definition,
  $d_{\blam\bmu}(q)\in\N[q,q^{-1}]$ is a Laurent polynomial with
  non-negative coefficients. Let
  $\dec=(d_{\blam\bmu}(q))_{\blam\in\Parts,\bmu\in\Klesh}$ be the
  \textbf{graded decomposition matrix} of~$\H$.

  The KLR algebra $\Rn$ is always $\Z$-free, however, it is not clear
  whether the same is true for the cyclotomic KLR algebra~$\R$. To prove
  this you cannot use the Graded Isomorphism \autoref{T:BKiso} because
  this result holds only over a field. Using extremely sophisticated
  diagram calculus, Li~\cite{Li:PhD} proved the following.

  \begin{Theorem}[Li~\cite{Li:PhD}]\label{T:Zfree}
    Suppose that $\Lambda\in P^+$. Then the quiver Hecke algebra $\R(\Z)$ is
    free as a $\Z$-module of rank $\ell^nn!$.  Moreover, $\R(\Z)$ is a graded
    cellular algebra with graded cellular basis
    $\set{\psi_{\s\t}|(\s,\t)\in\Std^2(\Parts)}$.
  \end{Theorem}

  Therefore,~$\R$ is free over any commutative ring and any field is a splitting
  field for~$\R$. Moreover, the graded Specht modules, together with their
  homogeneous bilinear forms, are defined over~$\Z$. The integrality of the
  graded Specht modules can also be proved using \autoref{T:Garnir} below.

  The next result lists some important properties of the $\psi$-basis.

  {\samepage
  \begin{Proposition}\label{P:PsiProperties}
    Suppose that $(\s,\t)\in\Std^2(\Parts)$ and that $\Zcal$ is an integral
    domain. Then:
    \begin{enumerate}%
      \item\cite[Lemma~5.2]{HuMathas:GradedCellular}
      If $\bi,\bj\in I^n$ then
      $\psi_{\s\t}=\delta_{\bi,\bi^\s}\delta_{\bj,\bi^\t}
                   e(\bi)\psi_{\s\t}e(\bj)$.
      \item\cite[Lemma~3.17]{HuMathas:QuiverSchurI}
      Suppose that $\psi_{\s\t}$ and $\hat\psi_{\s\t}$ are defined using
      different reduced expressions for the permutations $d(\s),d(\t)\in\Sym_n$.
      Then there exist $a_{\u\v}\in\Zcal$ such that
      $$\hat\psi_{\s\t}=\psi_{\s\t}+\sum_{(\u,\v)\Gdom(\s,\t)}a_{\u\v}\psi_{\u\v},$$
      where $a_{\u\v}\ne0$ only if $\bi^\u=\bi^\s$, $\bi^\v=\bi^\t$ and
      $\deg\u+\deg\v=\deg\s+\deg\t$.
      \item\relax\cite[Corollary~3.11]{HuMathas:GradedInduction}
      If $1\le r\le n$ then there exist $b_{\u\v}\in\Zcal$ such that
      $$\psi_{\s\t}y_r=\sum_{(\u,\v)\Gdom(\s,\t)}b_{\u\v}\psi_{\u\v},$$
      where $b_{\u\v}\ne0$ only if $\bi^\u=\bi^\s$, $\bi^\v=\bi^\t$ and
      $\deg\u+\deg\v=\deg\s+\deg\t+2$.
    \end{enumerate}
  \end{Proposition}
  }

  Part~(a) follows quickly using the relations in
  \autoref{D:QuiverRelations} and the definition of the $\psi$-basis. In
  contrast, parts~(b) and~(c) are proved by using \autoref{T:BKiso} to
  reduce to the seminormal basis. With part~(c), it is fairly easy to
  show that $b_{\u\v}\ne0$ only if $\u\gedom\s$. The difficult part is
  showing that $b_{\u\v}\ne0$ only if $\v\gedom\t$. Again, this is done
  using seminormal bases.

  Finally, we note that \autoref{T:Zfree} implies that $e(\bi)\ne0$ in $\R$ if
  and only if $\bi\in I^n_\Lambda=\set{\bi^\t|\t\in\Std(\Parts)}$, generalizing
  \autoref{P:SSimpleIdempotents}. In fact, if $F$ is a field and
  $\H(F)\cong\R(F)$ then it is shown in
  \cite[Lemma~4.1]{HuMathas:GradedCellular} that the non-zero KLR idempotents
  are a complete set of primitive (central) idempotents in the Gelfand-Zetlin
  algebra $\L(F)$ and that $\L(F)=\<y_1,\dots,y_n,e(\bi)\mid\bi\in I^n\>$. It
  follows that $\L(F)$ is a positively graded commutative algebra with one
  dimensional irreducible modules indexed by~$I^n_\Lambda$, up to shift. It would be
  interesting to find a (homogeneous) basis of~$\L(F)$. The author would also
  like to know whether~$\R$ is projective as a graded~$\L$-module.

  \subsection{Blocks and dual Specht modules}
  This section shows that the blocks of~$\H$ are graded symmetric algebras and
  it sketches the proof of an analogous statement that relates the graded
  Specht modules and their graded duals.

  \autoref{T:Blocks} describes the block decomposition of~$\H$ so, by
  \autoref{T:BKiso}, it gives the block decomposition of~$\R$.
  As in \autoref{E:RnBlocks}, set
  $$\R[\beta]=\R e_\beta,\qquad\text{where }
              e_\beta=\sum_{\bi\in I^\beta}e(\bi).$$
  It follows from \autoref{D:QuiverRelations} that $e_\beta$ is central in
  $\R$, so $\R[\beta]=e_\beta\R e_\beta$ is a two-sided ideal of~$\R$.
  Let $Q^+_n=Q^+_n(\Lambda)=\set{\beta\in Q^+|e_\beta\ne0}$ in $\R$. Similarly,
  let $\Parts[\beta]=\set{\blam\in\Parts|\bi^\blam\in I^\beta}
                    =\set{\blam\in\Parts|\beta^\blam=\beta}$.

  Combining \autoref{T:Zfree}, \autoref{T:BKiso} and \autoref{C:IntegralBlocks}
  we obtain the following.

  \begin{Theorem}
     Suppose that $\Lambda\in P^+$. Then
     $\R=\bigoplus_{\beta\in Q^+_n}\R[\beta]$
     is the decomposition of~$\R$ into indecomposable two-sided ideals.
     Moreover, $\R[\beta]$ is a graded cellular algebra with  cellular basis
     $\set{\psi_{\s\t}|(\s,\t)\in\Std^2(\Parts[\beta])}$ and weight
     poset~$\Parts[\beta]$.
  \end{Theorem}

  By virtue of \autoref{T:Zfree}, the block decomposition of~$\R$ holds
  over~$\Z$, even though we cannot think about the blocks as linkage classes of
  simple modules in this case. Compare with \autoref{T:SSKLRBasis} in the
  semisimple case.

  Suppose that $A$ is a graded $\Zcal$-algebra. Then $A$ is a
  \textbf{graded symmetric algebra} if there exists a homogeneous
  non-degenerate trace form $\tau\map A\Zcal$, where $\Zcal$ is in
  degree zero. That is, $\tau(ab)=\tau(ba)$ and if $0\ne a\in A$ then
  there exists $b\in A$ such that $\tau(ab)\ne0$. The map $\tau$ is
  \textbf{homogeneous of degree}~$d$ if $\tau(a)\ne0$ only if $\deg a=-d$.

  Fix $\beta\in Q^+$. The \textbf{defect} of $\beta$ is the non-negative integer
  \begin{equation*}\label{E:defect}
    \defect\beta =(\Lambda,\beta) - \frac12(\beta,\beta)
           =\frac12\Big((\Lambda,\Lambda)-(\Lambda-\beta,\Lambda-\beta)\Big).
  \end{equation*}
  Notice that $\defect\beta=\defect_\Lambda\beta$ depends on~$\Lambda$.
  If $\blam\in\Parts$ set $\defect\blam=\defect\beta^\blam$ (see
  \autoref{C:IntegralBlocks}). If $\lambda\in\Pcal_{1,n}$ is a partition then
  $\defect\lambda$ is equal to its $e$-weight; see, for
  example,~\cite[Proposition~2.1]{Fayers:AKweight} or the proof of
  \cite[Lemma~7.6]{LLT}.

  \autoref{E:dA}, and the definition of defect, readily implies the
  following combinatorial relationships between degrees, codegrees and
  defects.

  \begin{Lemma}\label{L:defect}
    Suppose that $\blam\in\Parts$.
    \begin{enumerate}
      \item\cite[Lemma~3.11]{BKW:GradedSpecht}
        If $A\in\Add_i(\blam)$ then
        $d_A(\blam)+1+d^A(\blam)=d_i(\blam)$ and
        $\defect(\blam{+}A)=\defect\blam+d_i(\blam)-1$.
      \item\cite[Lemma~3.12]{BKW:GradedSpecht}
        If $\s\in\Std(\blam)$ then $\deg\s+\codeg\s=\defect\blam$.
    \end{enumerate}
  \end{Lemma}

  In \autoref{D:psis} we defined two sets of elements $\{\psi_{\s\t}\}$ and
  $\{\psi'_{\s\t}\}$ in~$\R$.  Just as there are two versions of the Murphy
  basis, $\{m_{\s\t}\}$ and $\{m'_{\s\t}\}$, that are built from the trivial and
  sign representations of~$\H$~\cite{M:gendeg}, respectively, there are two
  versions of the $\psi$-basis. By \cite[Theorem~6.17]{HuMathas:GradedCellular},
  $\set{\psi'_{\s\t}|(\s,\t)\in\Std^2(\Parts)}$ is a second graded cellular basis
  of~$\H$ with weight poset $(\Parts,\ledom)$ and with
  $\deg\psi'_{\s\t}=\codeg\s+\codeg\t$. We warn the reader that we are following
  the conventions of \cite{HuMathas:QuiverSchurI}, rather than the notation of
  \cite{HuMathas:GradedCellular}. See \cite[Lemma~3.15 and
  Remark~3.12]{HuMathas:QuiverSchurI} for the translation.

  The bases $\{\psi_{\s\t}\}$ and $\{\psi'_{\u\v}\}$ of $\R$ are dual in the
  sense that if $(\s,\t),(\u,\v)\in\Std^2(\Parts[\beta])$ then, by
  \cite[Theorem~6.17]{HuMathas:GradedInduction},
  \begin{equation}\label{E:psidual}
    \psi_{\s\t}\psi'_{\t\s}\ne0\quad\text{and}\quad
    \psi_{\s\t}\psi'_{\u\v}\ne0\quad\text\quad\text{only if}\quad
          \bi^\t=\bi^\u\text{ and }\u\gedom\t.
  \end{equation}
  Let $\tau$ be the usual non-degenerate trace form
  on~$\H$~\cite{BK:HigherSchurWeyl,MM:trace}.  We can write
  $\tau=\sum_d\tau_d$, where $\tau_d$ is homogeneous of degree~$d\in\Z$.
  Let $\tau_\beta=\tau_{-2\defect\beta}$ be the homogeneous component
  of~$\tau$ of degree $-2\defect\beta$. By
  \cite[Theorem~6.17]{HuMathas:GradedInduction}, if
  $(\s,\t)\in\Std^2(\Parts)$ then $\tau_\beta(\psi_{\s\t}\psi'_{\t\s})\ne0$,
  so \autoref{E:psidual} implies the following.

  \begin{Theorem}[\protect{Hu-Mathas\cite[Corollary~6.18]{HuMathas:GradedCellular}}]
    \label{T:GradedSymmetric}
    Let $\beta\in Q^+_n$. Then~$\R[\beta]$ a graded symmetric algebra with
    homogeneous trace form $\tau_\beta$ of degree $-2\defect\beta$.
  \end{Theorem}

  It would be better to have an  intrinsic definition of~$\tau_\beta$ for
  $\R(\Z)$. Webster~\cite[Remark~2.27]{Webster:HigherRep} has given a
  diagrammatic description of a trace form on an arbitrary cyclotomic KLR
  algebra. It is unclear to the author how these two forms on~$\R$ are related.

  The $\psi'$-basis is a graded cellular basis of $\H$ so it defines a
  collection of graded cell modules. For $\blam\in\Parts[\beta]$, the
  \textbf{dual graded Specht module}~$S_\blam$ is the corresponding
  graded cell module determined by the $\psi'$-basis. The dual Specht
  module $S_\blam$ has basis $\set{\psi'_\t|\t\in\Std(\blam)}$, with
  $\deg\psi'_\t=\codeg\t$, and
  $$\gdim S_\blam=\sum_{\t\in\Std(\blam)}q^{\codeg\t}.$$
  We can identify $S_\blam\<\codeg\t_\blam\>$ with
  $(\psi'_{\t_\blam\t_\blam}+\Hlamp)\H$, where $\Hlamp$ is the two-sided ideal of
  $\H$ spanned by $\psi'_{\s\t}$ where $(\s,\t)\in\Std^2(\bmu)$ for some
  multipartition $\bmu$ such that $\blam\gdom\bmu$. Similarly, we can identify
  $S^\blam\<\deg\tlam\>$ with $(\psi_{\tlam\tlam}+\Hlam)\H$. By
  \autoref{E:psidual} there is a non-degenerate pairing
  $$\{\ ,\ \}\map{S^\blam\<\deg\tlam\>\times S_\blam\<\codeg\t_\blam\>}\Z$$
  given by $\{a+\Hlam, b+\Hlamp\}=\tau_\beta(ab^\star).$
  Hence, \autoref{L:defect} implies:

  \begin{Corollary}[\protect{%
    Hu-Mathas~\cite[Proposition 6.19]{HuMathas:GradedCellular}}]
    \leavevmode\newline\label{C:DualSpechts}%
    Suppose that $\blam\in\Parts$. Then
    $S^\blam\cong S_\blam^\circledast\<\defect\blam\>$
    and $S_\blam=(S^\blam)^\circledast\<\defect\blam\>$.
  \end{Corollary}

  This result holds for the Specht modules defined over~$\Z$ by
  \autoref{T:Zfree} or by \cite[Theorem~7.25]{KMR:UniversalSpecht}.

  There is an interesting byproduct of the proof of \autoref{C:DualSpechts}.
  In the ungraded setting the Specht module~$\UnS^\blam$ is isomorphic to the
  submodule of $\H$ generated by an element~$m_\blam T_{w_\blam}m_\blam'$; see
  \cite[Definition~2.1 and Theorem~2.9]{DuRui:branching}. By
  \cite[Corollary~6.21]{HuMathas:GradedCellular},
  $m_\blam T_{w_\blam}m_\blam'$ is homogeneous. In fact,
  $\psi_{\tlam\tlam}\psi_{w_\blam}\psi'_{\t_\blam\t_\blam}
                 =m_\blam T_{w_\blam}m_\blam'$
  and $\psi_{\tlam\tlam}\psi_{w_\blam}\psi'_{\t_\blam\t_\blam}\R\cong
  S^\blam\<\defect\blam+\codeg\t_\blam\>$.

  \subsection{Induction and restriction}
  The cyclotomic Hecke algebra $\H$ is naturally a subalgebra of $\H[n+1]$,
  and~$\H[n+1]$ is free as an $\H$-module by \autoref{E:AKBasis}. This
  gives rise to the usual induction and restriction functors. These
  functors can be decomposed into the ``classical'' $i$-induction and $i$-restriction
  functors, for $i\in I$, by projecting onto the blocks of these two
  algebras. As we will see, these functors are implicitly built into the
  graded setting.

  Recall that $I=\Z/e\Z$ and $\Lambda\in P^+$. For each $i\in I$ define
  $$e_{n,i}=\sum_{\bj\in I^n}e(\bj\vee i)\in\R[n+1].$$ The relations
  for~$\R[n+1]$ in \autoref{D:QuiverRelations} imply that $e_{n,i}$ is an
  idempotent and that $\sum_{i\in I}e_{n,i}=\sum_{\bi\in I^{n+1}}e(\bi)$ is the
  identity element of $\R[n+1]$.

  Let $\rep(\R)$, and $\rep(\R[\beta])$ for $\beta\in Q^+$, be the
  category of finite dimension (graded) $\R$-modules, respectively,
  $\R[\beta]$-modules.  Similarly, let $\proj(\R)$ and
  $\proj(\R[\beta])$ be the categories of finitely generated projective
  modules for these algebras.

  \begin{Lemma}\label{L:embedding}
    Suppose that $i\in I$ and that $\Zcal$ is an integral domain. Then there is a
    (non-unital) embedding of graded algebras $\R\hookrightarrow\R[n+1]$ given by
    $$e(\bj)\mapsto e(\bj\vee i),\quad
      y_r\mapsto e_{n,i}y_r\quad\text{and}\quad
      \psi_s\mapsto e_{n,i}\psi_s,$$
      for $\bj\in I^n$, $1\le r\le n$ and $1\le s<n$. This  map induces an
      exact functor
      $$\iInd\map{\rep(\R)}{\rep(\R[n+1])}; M\mapsto M\otimes_{\R}e_{n,i}\R[n+1].$$
      Moreover, $\Ind=\bigoplus_{i\in I}\iInd$ is the graded induction
      functor from $\rep(\R)$ to $\rep(\R[n+1])$.
  \end{Lemma}

  \begin{proof}
    The images of the homogeneous generators of $\R$ under this embedding
    commute with $e_{n,i}$, which implies that this map defines a non-unital
    degree preserving homomorphism from $\R$ to~$\R[n+1]$.  This map is an
    embedding by \autoref{T:Zfree}. The remaining claims follow because, by
    definition, $e_{n,i}$ is an idempotent and $\sum_{i\in I}e_{n,i}$ is the
    identity element of~$\R[n+1]$.
  \end{proof}

  The \textbf{$i$-induction functor} $\iInd$ functor is obviously a left adjoint
  to the \textbf{$i$-restriction} functor $\iRes$, which sends an
  $\R[n+1]$-module~$M$ to
  $$\iRes M= Me_{n,i}\cong \ZHom_{\R}(e_{n,i}\R,M).$$
  A much harder fact is that these functors are two-sided adjoints.

  \begin{Theorem}[%
    \protect{Kashiwara~\cite[Theorem~3.5]{Kashiwara:KLRBiadjointness}}]
    \label{T:Biadjoint}
    Suppose $i\in I$. Then $(E_i,F_i)$ is a biadjoint pair.
  \end{Theorem}

  Kashiwara proves this theorem by  constructing explicit homogeneous
  adjunctions. He does this for any cyclotomic quiver Hecke algebras
  defined by a symmetrizable Cartan matrix. As we do not need this
  result we feel justified in stating it now, even though its proof
  builds upon Kang and Kashiwara's proof that the cyclotomic quiver
  Hecke algebras of arbitrary type categorify the integrable highest
  weight modules of the corresponding quantum
  group~\cite{KangKashiwara:CatCycKLR}; compare with \autoref{P:Fock}
  and \autoref{C:GradedDecomp} below. This biadjointness property is
  also a consequence of Rouquier's Kac-Moody categorification
  axioms~\cite[Theorem~5.16]{Rouq:2KM}.  \autoref{T:Biadjoint} was
  conjectured by Khovanov-Lauda~\cite{KhovLaud:diagI}.

  Recall from \autoref{E:dual} that $\circledast$ defines a graded
  duality on $\rep(\R)$. Similarly,  define
  $\#$ to be the graded functor given by
  \begin{equation}\label{E:HashDual}
    M^\#=\ZHom_{\R}(M,\R),\qquad\text{ for $M\in\rep(\R)$},
  \end{equation}
  where the action of $\R$ on $M^\#$ is given by
  $$(f\cdot h)(m)=h^\star f(m),\qquad
       \text{for  $f\in M^\#$, $h\in\R$ and $m\in M$}.
  $$
  We consider
  $\circledast$ and $\#$ as endofunctors of
  $\rep(\R)=\bigoplus_\beta\rep(\R[\beta])$ and
  $\proj(\R)=\bigoplus_\beta\proj(\R[\beta])$.  As noted in
  \cite[Remark~4.7]{BK:GradedDecomp}, \autoref{T:GradedSymmetric}
  implies that these two functors agree up to shift.

  \begin{Lemma}\label{L:SameDuals}
    As endofunctors of $\rep(\R[\beta])$,
    there is an isomorphism of functors
    $\#\cong\<2\defect\beta\>\circ\circledast$.
  \end{Lemma}

  \begin{proof}
    By \autoref{T:GradedSymmetric},
    $\R[\beta]\cong(\R[\beta])^\circledast\<2\defect\beta\>$.
    If $M\in\rep(\R[\beta])$ then
    \begin{align*}
       M^\# & =\ZHom_{\R[\beta]}(M,\R[\beta])
           =\ZHom_{\R[\beta]}\(M,(\R[\beta])^\circledast\<2\defect\beta\>\)\\
      &\cong\ZHom_{\R[\beta]}\(M,\ZHom_{\Zcal}(\R[\beta],\Zcal)\)\<2\defect\beta\>\\
      &\cong\ZHom_{\Zcal}\(M\otimes_{\R[\beta]}\R[\beta],\Zcal\)\<2\defect\beta\>\\
      &\cong M^\circledast\<2\defect\beta\>,
    \end{align*}
    where the third isomorphism is the standard adjointness of tensor and hom.
    As all of these isomorphisms are functorial, the lemma follows.
  \end{proof}

  By well-known arguments, $(M^\circledast)^\circledast\cong M$ for all
  $M\in\rep(\H)$. Hence, $(M^\#)^\#\cong M$ by \autoref{L:SameDuals}.
  Therefore, $\circledast$ and $\#$ define self-dual equivalences on
  $\rep(\R)$ and $\proj(\R)$.

  \begin{Proposition}\label{P:CommutingDuals}
    Suppose that $\beta\in Q^+$ and $i\in I$. Then there are functorial isomorphisms
    \begin{align*}
    \circledast\circ\iRes&\cong\iRes\circ\circledast
          \map{\rep(\R[n+1])}\rep(\R),\\
      \#\circ\iInd&\cong\iInd\circ\#\map{\proj(\R)}\proj(\R[n+1]).
    \end{align*}
  \end{Proposition}

  \begin{proof}
    The isomorphism $\circledast\circ\iRes\cong\iRes\circ\circledast$
    is immediate from the definitions. For the second isomorphism,
    recall that if $P\in\proj(\R[\beta])$ then
    $\ZHom_{\R}(P,M)\cong\ZHom_{\R}(M,\R)\otimes_{\R}M$, for any
    $\R$-module~$M$. Now,
    $$ (e_{n,i}\R[n+1])^\#
        =\ZHom_{\R[n+1]}(e_{n,1}\R[n+1],\R[n+1])\cong e_{n,i}\R[n+1],$$
    the last isomorphism following because $e_{n,i}^\star=e_{n,i}$. Therefore,
    \begin{align*}
      \iInd(P^\#)
        &=\ZHom_{\R}(P,\R)\otimes_{\R}e_{n,i}\R[n+1]
         \cong\ZHom_{\R}(P,e_{n,i}\R[n+1])\\
         &\cong\ZHom_{\R}\(P,\ZHom_{\R[n+1]}(e_{n,1}\R[n+1],\R[n+1])\)\\
        &\cong\ZHom_{\R[n+1]}(P\otimes_{\R}e_{n,i}\R[n+1],\R[n+1])\\
        &\cong(\iInd P)^\#,
    \end{align*}
    where the second last isomorphism is the usual tensor-hom adjointness.
  \end{proof}

  It follows from \autoref{P:CommutingDuals} and \autoref{L:SameDuals}
  that the functors $\circledast$ and $\iInd$, and $\#$ and $\iRes$,
  commute only up to shift. This difference in degree shift is what
  makes \autoref{L:UqCommute} work below.

  We next describe the effect of the $i$-induction and $i$-restriction
  functors on the graded Specht modules, for $i\in I$. This result
  generalizes the well-known (ungraded) branching rules for the
  symmetric group~\cite[Example~17.16]{James} and the cyclotomic Hecke
  algebras~\cite{AM:simples,RyomHansen:GradedTranslation,M:SpechtInduction}.

  Recall the integers $d^A(\blam)$ and $d_A(\blam)$ from \autoref{E:dA}.

  {\samepage
  \begin{Theorem}\label{T:Induction}
    Suppose that $\Zcal$ is an integral domain and $\blam\in\Parts$.
    \begin{enumerate}
    \item\cite[Main theorem]{HuMathas:GradedInduction}
      Let $A_1<A_2\dots<A_z$ be the addable $i$-nodes of~$\blam$.
    Then $\iInd S^\blam$ has a graded Specht
    filtration
    $$0=I_0\subset I_1\subset\dots\subset I_z=\iInd S^\blam,$$
    such that $I_j/I_{j-1}\cong S^{\blam+A_j}\<d_{A_j}(\blam)\>$,
    \item \cite[Theorem~4.11]{BKW:GradedSpecht}
      Let $B_1>B_2>\dots>B_y$ be the removable $i$-nodes of~$\blam$.
      Then $\iRes S^\blam$ has a graded Specht filtration
    $$0=R_0\subset R_1\subset\dots\subset R_y=\iRes S^\blam,$$
    such that $R_j/R_{j-1}\cong S^{\blam-B_j}\<d_{B_j}(\blam)\>$,
    for $1\le j\le y$.
   \item\cite[Corollary~4.7]{HuMathas:GradedInduction}
     Let $A_1>A_2>\dots>A_z$ be the addable $i$-nodes of~$\blam$.
   Then $\iInd S_\blam$ has a graded Specht
   filtration
   $$0=I_0\subset I_1\subset\dots\subset I_z=\iInd S_\blam,$$
   such that $I_j/I_{j-1}\cong S_{\blam+A_j}\<d^{A_j}(\blam)\>$,
   for $1\le j\le z$.
   \item Let $B_1<B_2<\dots<B_y$ be the removable $i$-nodes of~$\blam$.
      Then $\iRes S^\blam$ has a graded Specht filtration
    $$0=R_0\subset R_1\subset\dots\subset R_y=\iRes S_\blam,$$
    such that $R_j/R_{j-1}\cong S_{\blam-B_j}\<d^{B_j}(\blam)\>$,
    for $1\le j\le y$.
  \end{enumerate}
  \end{Theorem}
  }

  Observe that parts (a) and~(c), and parts~(b) and~(d), are equivalent
  by \autoref{C:DualSpechts} (and \autoref{L:defect}).

  Part~(b) is proved using the fact that the action of~$\H$ on the
  $\psi$-basis is compatible with restriction. Part~(a), which was
  conjectured by Brundan, Kleshchev and
  Wang~\cite[Remark~4.12]{BKW:GradedSpecht}, is proved by extending
  elegant ideas of Ryom-Hansen~\cite{RyomHansen:GradedTranslation} to
  the graded setting using \cite{HuMathas:GradedCellular}.

  \subsection{Grading Ariki's Categorification
  Theorem}\label{S:Categorification}

  The aim of this section is to prove the Ariki-Brundan-Kleshchev
  Categorification Theorem~\cite{Ariki:can} that connects the canonical
  bases of $\Uq$-modules with the simple and projective indecomposable
  $\R$-modules in characteristic zero.  We give a new proof of Brundan
  and Kleshchev's theorem~\cite{BK:GradedDecomp} that the cyclotomic KLR
  algebras of type~$A$ categorify the integrable highest weight modules
  of~$\Uq$.  Our argument runs parallel to Brundan and Kleshchev's with
  the key difference being that we use the representation theory of~$\H$,
  and in particular the graded branching rules, to construct a
  \textit{bar involution} on the Fock space. In this way we are able to
  show that the canonical basis is categorified by the basis of simple
  $\H$-modules if and only if the graded decomposition numbers are
  polynomials. As a consequence, Ariki's categorification
  theorem~\cite{Ariki:can} lifts to the graded setting.

  Throughout this section we assume that the Hecke algebra~$\H$ is
  defined over a field~$\F$.  In the end we will assume that~$\F$ is a
  field of characteristic zero, however, almost all of the results in
  this section hold over an arbitrary field. We delay introducing
  the quantum group $\Uq$ until we actually need it because we want to
  emphasize the role that the quantum group is playing in the
  representation theory of~$\H$.

  For the time being fix an integer $n\ge0$. Very soon we will vary~$n$.
  Let $\A=\Z[q,q^{-1}]$ be the ring of Laurent polynomials in~$q$
  over~$\Z$.

  Let $\rep(\H)$ be the category of finitely generated graded
  $\H$-modules and let $\proj(\H)$ be the category of finitely generated
  projective graded $\H$-modules.  Let $\repH$ and $\projH$ be the
  Grothendieck groups of these categories. If $M$ is a finitely
  generated $\H$-module let $[M]$ be its image in $\repH$.  Abusing
  notation slightly, if $M$ is projective we also let~$[M]$ be its image
  in~$\projH$.  Consider $\repH$ and $\projH$ as $\A$-modules by
  letting~$q$ act as the grading shift functor: $[M\<d\>]=q^d[M]$, for
  $d\in\Z$.

  \begin{Definition}
    Suppose that $\bmu\in\Klesh$. Let $Y^\bmu$ be the projective cover of
    $D^\bmu$ in $\rep(\H)$.
  \end{Definition}

  Importantly, the module $Y^\bmu$ is graded. Since $Y^\bmu$ is
  indecomposable, the grading on~$Y^\bmu$ is uniquely determined by the
  surjection $Y^\bmu\surjection D^\bmu$, for~$\bmu\in\Klesh$. We use the
  notation~$Y^\bmu$ because these modules are special cases of the
  graded lifts of the \textit{Young modules} constructed
  in~\cite{M:tilting}; see \cite[\S5.1]{HuMathas:QuiverSchurI} and
  \cite[\S2.6]{Maksimau:QuiverSchur}. (The symbol $P^\bmu$ is usually
  reserved for the projective indecomposable modules of the cyclotomic Schur
  algebras~\cite{BK:HigherSchurWeyl,DJM:cyc,HuMathas:QuiverSchurI,%
  StroppelWebster:QuiverSchur}.)

  By definition, the Grothendieck groups $\repH$ and
  $\projH$ are free $\A$-modules that come equipped with distinguished bases:
  $$\repH=\bigoplus_{\bmu\in\Klesh}\A[D^\bmu]\quad\text{and}\quad
    \projH=\bigoplus_{\bmu\in\Klesh}\A[Y^\bmu],
  $$
  respectively. Recall from \autoref{E:GradedDec} that
  $\dec=\(\d_{\blam\bmu}(q)\)$ is the graded decomposition matrix
  of~$\H$. If $\blam\in\Parts$ and $\bmu\in\Klesh$ then
  in~$\repH$,
  $$ [S^\blam] = \sum_{\substack{\btau\in\Klesh\\\blam\gedom\btau}}
                      d_{\blam\btau}(q)[D^\btau]
  \quad\text{and}\quad
  [Y^\bmu] = \sum_{\substack{\bsig\in\Parts\\\bsig\gedom\bmu}}
         d_{\bsig\bmu}(q)[S^\bsig],$$
  where the second formula comes from \autoref{C:CartanSymmetric}.
  By \autoref{T:GradedSimples}(c), the submatrix
  $\deck=\(d_{\blam\bmu}(q)\)_{\blam,\bmu\in\Klesh}$ of the graded
  decomposition matrix~$\dec$ is invertible over~$\A$ with inverse
  $$\dinv  = (\deck)^{-1}= \(e_{\blam\bmu}(-q)\)_{\blam,\bmu\in\Klesh}.$$
  (The reason why we consider $e_{\blam\bmu}(-q)$ as a Laurent
  polynomial in $-q$ is explained after \autoref{C:GradedDecomp} below.)
  Hence, if $\blam\in\Klesh$ then
  $$[D^\blam] = \sum_{\bmu\in\Klesh}e_{\blam\bmu}(-q)[S^\bmu].$$
  Consequently, $\set{[S^\bmu]|\bmu\in\Klesh}$ is a second $\A$-basis of $\repH$.

  The set of projective indecomposable $\H$-modules $\{[Y^\bmu]\}$ is the only
  natural basis of the split Grothendieck group $\projH$. Somewhat
  artificially, but motivated by the formulas above, for $\bmu\in\Klesh[]$
  define
  $$X_\bmu=\sum_{\blam\in\Klesh}e_{\blam\bmu}(-q)[Y^\blam]\in\projH.$$
  Then $\set{X_\bmu|\bmu\in\Klesh[]}$ is an $\A$-basis of $\projH$.
  We will use the bases $\set{[S^\bmu]|\bmu\in\Klesh[]}$ and
  $\set{X_\bmu|\Klesh[]}$ of $\repH$ and $\projH$, respectively, to
  construct new distinguished bases of the Grothendieck groups.

  The \textbf{bar involution} on $\A=\Z[q,q^{-1}]$ is the unique
  $\Z$-linear map such that $\overline{q^d}=q^{-d}$, for $d\in\Z$. In
  particular, $\gdim M^\circledast=\overline{\gdim M}$, for any
  $\R$-module $M$. A \textbf{semilinear} map of $\A$-modules is a
  $\Z$-linear map $\theta\map MN$ such that
  $\theta(f(q)m)=\overline{f(q)}\theta(m)$, for all $f(q)\in\A$ and
  $m\in M$.

  A \textbf{sesquilinear} map $f\map{M\times N}\A$, where $M$ and $N$
  are $\A$-modules, is a function that is semilinear in the first
  variable and linear in the second.  Let
  $\<\ ,\ \>\map{\projH\times\repH}\A$ be the sesquilinear pairing
  \begin{equation}\label{E:CartanPairing}
     \<[P],[M]\>=\gdim\ZHom_{\H}(P,M),
  \end{equation}
  for $P\in\proj(\H)$ and $M\in\rep(\H)$.
  This pairing is naturally \textbf{sesquilinear} because
  $\ZHom_{\H}(P\<k\>,M)\cong\ZHom_{\H}(P,M\<-k\>)$, for any $k\in\Z$.

  The functors $\circledast$ and $\#$, of \autoref{E:dual} and
  \autoref{E:HashDual}, induce semilinear automorphisms of
  the Grothendieck groups $\repH$ and $\projH$:
  $$ [P]^\#=[P^\#], \quad\text{and}\quad
     [M]^\circledast=[M^\circledast]
  $$
  for $M\in\rep(\H)$ and $P\in\proj(\H)$.
  The next result is fundamental.

  \begin{Lemma}\label{L:Shapovalov}
    Suppose that $[P]\in\Proj[\A]$ and $[M]\in\Rep[\A]$. Then
    $$\<[P],[M]^\circledast\>=\overline{\<[P]^\#,[M]\>}.$$
  \end{Lemma}

  \begin{proof}Applying the definitions, and tensor-hom adjointness,
    \begin{align*}
      \<[P],[M]^\circledast\>
      &=\gdim\ZHom_{\R}(P,M^\circledast)
       =\gdim\ZHom_{\R}\(P,\ZHom_{\R}(M,\F)\)\\
      &=\gdim\ZHom_{\R}(P\otimes_{\R}M,\F)
       =\gdim(P\otimes_{\R}M)^\circledast\\
      &=\overline{\gdim P\otimes_{\R}M}
       =\overline{\gdim\ZHom_{\R}(P^\#,\R)\otimes_{\R}M}\\
      &=\overline{\gdim\ZHom_{\R}(P^\#,M)}
       =\overline{\<[P]^\#,[M]\>}.
    \end{align*}
    For the second last line, note that
    $\ZHom_{\R}(Q,M)\cong\ZHom_{\R}(Q,\R)\otimes_{\R}M$ whenever~$Q$
    is projective.
  \end{proof}

  \begin{Lemma}\label{L:DualBases}
    Suppose that $\blam,\bmu\in\Klesh$. Then
    $$\<[Y^\blam],[D^\bmu]\>=\delta_{\blam\bmu}
                            =\<X_\blam,[S^\bmu]^\circledast\>.
    $$
  \end{Lemma}

  \begin{proof}
    The first equality is immediate from the definition of the sesquilinear form
    $\<\ ,\ \>$ because $Y^\blam$ is the projective
    cover of~$D^\blam$, for $\blam\in\Klesh$. For the second equality,
    using the fact that $\circledast$ is semilinear,
    \begin{align*}
      \<X_\blam,[S^\bmu]^\circledast\>
          &=\sum_{\bsig\in\Klesh} \overline{e_{\bsig\blam}(-q)}
              \<[Y^\bsig],[S^\bmu]^\circledast\>\\
          &=\sum_{\bsig,\btau\in\Klesh} \overline{e_{\bsig\blam}(-q)}\,
              \overline{d_{\bmu\btau}(q)}\<[Y^\bsig],[D^\btau]\>\\
          &=\sum_{\bsig\in\Klesh}
              \overline{d_{\bmu\bsig}(q)}\,\overline{e_{\bsig\blam}(-q)}
           =\delta_{\blam\bmu},
    \end{align*}
    where the last equality follows because $\dinv=(\deck)^{-1}$.
  \end{proof}

  \begin{Lemma}\label{L:BarInvariance}
    Suppose that $\bmu\in\Klesh$. Then
      $[Y^\bmu]^\#=[Y^\bmu]$,
      $[D^\bmu]^\circledast=[D^\bmu]$,
    $$
      (X_\bmu)^\# = X_\bmu+
      \sum_{\substack{\bsig\in\Klesh\\\bsig\gdom\bmu}}a_{\bsig\bmu}(q)X_\bsig
      \quad\text{and}\quad
      [S^\bmu]^\circledast = [S^\bmu] +
      \sum_{\substack{\btau\in\Klesh\\\bmu\gdom\btau}}a^{\bmu\btau}(q)[S^\btau],
    $$
    for some Laurent polynomials $a_{\bsig\bmu}(q),a^{\btau\bmu}(q)\in\A$.
  \end{Lemma}

  \begin{proof}
    That $[D^\bmu]^\circledast=[D^\bmu]$ is immediate by
    \autoref{C:GradedSimples}, whereas $[Y^\bmu]^\#=Y^\bmu$
    because $Y^\bmu$ is a direct summand of~$\H$ --- alternatively, use
    \autoref{L:DualBases} and \autoref{L:Shapovalov}. If $\bmu\in\Klesh$
    then, by \autoref{T:GradedSimples},
    \begin{align*}
    [S^\bmu]^\circledast
    &=\Big(\sum_{\substack{\bnu\in\Klesh\\\bmu\gedom\bnu}}
                d_{\bmu\bnu}(q)[D^\bnu]\Big)^\circledast
         =\sum_{\substack{\bnu\in\Klesh\\\bmu\gedom\bnu}}
              \overline{d_{\bmu\bnu}(q)}\,[D^\bnu]\\
        &=[S^\bmu]+\sum_{\substack{\btau\in\Klesh\\\bmu\gdom\btau}}
           \Big(\sum_{\substack{\bnu\in\Klesh\\\bmu\gedom\bnu\gedom\btau}}
           \overline{d_{\bmu\bnu}(q)}\,e_{\bnu\btau}(-q)\Big)[S^\btau]
    \end{align*}
    as claimed. Note that $d_{\bmu\bmu}(q)=1=e_{\bmu\bmu}(-q)$.

    Finally, we can compute $(X_\bmu)^\#$ by writing $X_\bmu=\sum_\bmu
    e_{\bmu\blam}(-q)[Y^\blam]$ and then using essentially the same
    argument to show that $(X_\bmu)^\#$ can be written in the required
    form.  Alternatively, use \autoref{L:DualBases} and
    \autoref{L:Shapovalov}.
  \end{proof}

  The triangularity of the action of $\circledast$ and $\#$ on $\Rep[\A]$ and
  $\Proj[\A]$, respectively, has the following easy but important consequence.

  \begin{Proposition}\label{P:CanonicalBasis}
    Suppose that $\F$ is a field.
    Then there exist unique bases $\set{B_\bmu|\bmu\in\Klesh}$ and
    $\set{B^\bmu|\bmu\in\Klesh}$ of~$\projH$ and $\repH$, respectively,
    such that $(B_\bmu)^\#=B_\bmu$, $(B^\bmu)^\circledast=B^\bmu$
    $$
    B_\bmu = X_\bmu +
    \sum_{\substack{\bsig\in\Klesh\\\bsig\gdom\bmu}}b_{\bsig\bmu}(q)X_\bsig
    \quad\text{and}\quad
    B^\bmu = [S^\bmu] +
    \sum_{\substack{\btau\in\Klesh\\\bmu\gdom\btau}}b^{\bmu\btau}(q)[S^\btau]
    $$
    for \underline{polynomials}
    $b^{\bmu\bsig}(q), b_{\bsig\bmu}(q)\in\delta_{\bsig\bmu}+q\Z[q]$.
  \end{Proposition}

  \begin{proof}The existence and uniqueness of these two bases follows
    immediately from \autoref{L:BarInvariance} by \textit{Lusztig's
    Lemma}~\cite[Lemma 24.2.1]{Lusztig:QuantBook}.  We give a variation
    of Lusztig's argument for the basis $\{B^\bmu\}$.

    Fix a multipartition $\bmu\in\Klesh$, for some $n\ge0$, and suppose
    that~$B^\bmu$ and~$\dot B^\bmu$ are two elements of $\repH$ with the
    required properties. By assumption the element $B^\bmu-\dot B^\bmu$ is
    $\circledast$-invariant and we can write
    $$B^\bmu-\dot B^\bmu=\sum_{\bmu\gdom\btau}\dot b^{\bmu\btau}(q)[S^\btau],$$
    for some polynomials $\dot b^{\bmu\btau}(q)\in q\Z[q]$. As the left-hand
    side is $\circledast$-invariant and
    $\overline{\dot b^{\bmu\btau}(q)}\in q^{-1}\Z[q^{-1}]$, so
    \autoref{L:BarInvariance} forces $B^\bmu=\dot B^\bmu$.

    To prove existence, we argue by induction on dominance. If~$\bmu$ is minimal
    in~$\Klesh$ then we can take $B^\bmu=[S^\bmu]=[D^\bmu]$ by
    \autoref{L:BarInvariance}. If~$\bmu\in\Klesh$ is not minimal with respect to
    dominance then set $\dot B^\bmu=[D^\bmu]$. Then
    $$(\dot B^\bmu)^\circledast=\dot B^\bmu\quad\text{and}\quad
         \dot B^\bmu=[S^\bmu]+\sum_{\substack{\btau\in\Klesh\\\bmu\gdom\btau}}
                     \dot b^{\bmu\btau}(q)[S^\btau],$$
    for some Laurent polynomials $\dot b^{\bmu\btau}(q)\in\Z[q,q^{-1}]$. If
    $\dot b^{\bmu\btau}(q)\in q\Z[q]$, for all~$\bmu\gdom\btau$, then
    $B^\bmu=\dot B^\bmu$ has all of the required properties. Otherwise,
    pick any multipartition~$\bmu\gdom\bnu$ that is maximal with respect to
    dominance such that $\dot b^{\bmu\bnu}(q)\notin q\Z[q]$. By induction,
    there exists an element $B^\bnu$ with all of the required properties.
    Replace $\dot B^\bmu$ with the element $\dot B^\bmu-p^{\bmu\bnu}(q)B^\bnu$, where
    $p^{\bmu\bnu}(q)$ is the unique Laurent polynomial such that
    $\overline{p^{\bmu\bnu}(q)}=p^{\bmu\bnu}(q)$ and
    $\dot b^{\bmu\bnu}(q)-p^{\bmu\bnu}(q)\in q\Z[q]$.  Then
    $(\dot B^\bmu)^\circledast=\dot B^\bmu$ and the
    coefficient of~$[S^\bnu]$ in~$\dot B^\bmu$ belongs
    to~$q\Z[q]$. Continuing in this way, after finitely many steps we
    will construct an element $B^\bmu$ with the required properties.
  \end{proof}

  \begin{Corollary}\label{C:ShapovalovDual}
    Suppose that $\blam,\bmu\in\Klesh[]$. Then
    $$\<B_\bmu,B^\blam\>
          =\sum_{\substack{\bsig\in\Klesh[]\\\blam\gedom\bsig\gedom\bmu}}
                b^{\blam\bsig}(q)b_{\bsig\bmu}(q)
          =\delta_{\blam\bmu}.$$
  \end{Corollary}

  \begin{proof}
    If $\bsig,\btau\in\Klesh$ then
    $\<X_\bsig,[S^\btau]^\circledast\>=\delta_{\bsig\btau}$ by
    \autoref{L:DualBases}. Therefore, since the form
    $\<\ ,\ \>$ is sesquilinear and $B^{\blam\circledast}=B^\blam$,
    \begin{align*}
       \<B_\bmu,B^\blam\>&= \<B_\bmu,B^{\blam\circledast}\>
          =\sum_{\bsig\gedom\bmu}\sum_{\blam\gedom\btau}
                 \overline{b_{\bsig\bmu}(q)}\,\overline{b^{\blam\btau}(q)}
                          \<X_\bsig,[S^\btau]^\circledast\>\\
          &=\sum_{\blam\gedom\bsig\gedom\bmu}
          \overline{b^{\blam\bsig}(q)}\,\overline{b_{\bsig\bmu}(q)}.
    \end{align*}
    In particular, $(B_\bmu,B^\blam)\in \delta_{\blam\bmu}+q^{-1}\Z[q^{-1}]$.
    On the other hand,
    $$\<B_\bmu,B^\blam\>=\<B_\bmu^\#,B^\blam\>
                      =\overline{\<B_\bmu,B^{\blam\circledast}\>}
                      =\overline{\<B_\bmu,B^\blam\>}$$
    by \autoref{L:Shapovalov},
    Therefore, $\<B_\bmu,B^\blam\>=\delta_{\blam\bmu}$ as this is the only bar
    invariant polynomial in $\delta_{\blam\bmu}+q^{-1}\Z[q^{-1}]$.
  \end{proof}

  Applying \autoref{L:SameDuals} to \autoref{P:CanonicalBasis} we obtain.

  \begin{Corollary}\label{C:SecondDuals}
    Suppose that $\bmu\in\Klesh[]$. Then
    $(q^{-\defect\bmu}B_\bmu)^\circledast=q^{-\defect\bmu}B_\bmu$ and
    $(q^{\defect\bmu}B^\bmu)^\#=q^{\defect\bmu}B^\bmu$
  \end{Corollary}

  In order to link the bases $\{B_\bmu\}$ and $\{B^\bmu\}$ with the
  representation theory of~$\H$ we need to introduce the quantum
  group~$\Uq$.

  The \textbf{quantum group} $\Uq$ associated with the
  quiver~$\Gamma_e$ is the $\Q(q)$-algebra generated by
  $\set{ E_i, F_i, K^\pm_i|i\in I}$, subject to the
  relations:
  \bgroup
     \setlength{\abovedisplayskip}{1pt}
     \setlength{\belowdisplayskip}{1pt}
    \begin{align*}
      K_iK_j &=K_jK_i,
         &K_iK_i^{-1}&=1,
         &[E_i,F_j]&=\delta_{ij}\frac{K_i-K_i^{-1}}{q-q^{-1}},
    \end{align*}
    \begin{align*}
        K_iE_jK_i^{-1}&=q^{c_{ij}}E_j,
        &K_iF_jK_i^{-1}&=q^{-c_{ij}}F_j,
    \end{align*}
    \begin{align*}
     \sum_{0\le c\le 1-c_{ij}}(-1)^c\qbinom{1-c_{ij}}{c}E_i^{1-c_{ij}-c}E_jE_i^c&=0,\\
     \sum_{0\le c\le 1-c_{ij}}(-1)^c\qbinom{1-c_{ij}}{c}F_i^{1-c_{ij}-c}F_jF_i^c&=0,
    \end{align*}
  \egroup
  where $\qbinom dc=\diag d!/\diag c!\diag{d-c}!$ and
  $\diag m!=\prod_{k=1}^m(q^k-q^{-k})/(q-q^{-1})$, for integers $c<d,m\in\N$.
  Then $\Uq$ is a Hopf algebra with coproduct determined by
  $\Delta(K_i)=K_i\otimes K_i$,
  $\Delta(E_i)=E_i\otimes K_i+1\otimes E_i$ and
  $\Delta(K_i)=F_i\otimes1+ K_i^{-1}\otimes F_i$, for $i\in I$.
  A self-contained account of much of what we need can be found in
  Ariki's book~\cite{Ariki:book}. See also
  \cite[\S3.1]{Lusztig:QuantBook} and \cite{BK:GradedDecomp}.

  The \textbf{combinatorial Fock space} $\Fock[\A]$ is the free
  $\A$-module with basis the set of symbols
  $\set{\fock|\blam\in\Parts[]}$, where
  $\Parts[]=\bigcup_{n\ge0}\Parts$. For future use, let
  $\Klesh[]=\bigcup_{n\ge0}\Klesh$ .  Set
  $\Fock=\Fock[\A]\otimes_\A\Q(q)$.  Then, $\Fock$ is an infinite
  dimensional $\Q(q)$-vector space. We consider
  $\set{\fock|\blam\in\Parts[]}$ as a basis of~$\Fock$ by identifying
  $\fock$ and $\fock\otimes1_{\Q(q)}$.

  Recall the integers $d^A(\blam)$, $d_B(\blam)$ and
  $d_i(\blam)$ from \autoref{E:dA}.

  \begin{Theorem}[Hayashi \cite{Hay,MisraMiwa}]
    Suppose that $\Lambda\in P^+$. Then $\Fock$ is an integrable
    $\Uq$-module with $\Uq$-action determined by
    $$E_i\fock = \sum_{B\in\Rem_i(\blam)} q^{d_B(\blam)}\fock[\blam{-}B]
    \quad\text{and}\quad
    F_i\fock=\sum_{A\in\Add_i(\blam)}q^{-d^A(\blam)}\fock[\blam{+}A],$$
    and $K_i\fock=q^{d_i(\blam)}\fock$,
    for all $i\in I$ and $\blam\in\Parts$.
    \label{T:FockSpace}
  \end{Theorem}

  \begin{Remark}\label{R:FockSpaceAction}
    A slightly different action of $\Uq$ on the Fock space is used in
    many places in the literature, such as \cite{LLT,Ariki:book,M:Ulect}.
    As is already evident, and will be made precise in
    \autoref{P:FormComparision} below, the $\Uq$-action on the Fock
    space is closely related to induction and restriction for the graded
    Specht modules. The $\Uq$-action on the Fock space used in
    \cite{LLT,Ariki:book,M:Ulect} corresponds to the action of the
    induction and restriction functors on the \textit{dual} graded Specht
    modules. Equivalently, this difference in the $\Uq$-action arises
    because, ultimately, we will work with an action of~$\Uq$
    on the Grothendieck groups of the finitely generated
    $\R$-modules, whereas these other sources consider the corresponding adjoint
    action on the projective Grothendieck groups.
  \end{Remark}

  Hayashi~\cite{Hay} considered only the special case when
  $\Lambda=\Lambda_0$, however, this implies the general case using the
  coproduct of~$\Uq$ because
  $$\Fock\cong\mathscr{F}^{\Lambda_{\overline{\kappa_1}}}_{\Q(q)}\otimes\dots
  \otimes\mathscr{F}^{\Lambda_{\overline{\kappa_\ell}}}_{\Q(q)}$$ as
  $\Uq$-modules. The crystal and canonical bases of~$\Fock$, which were
  first studied in \cite{JMMO,MisraMiwa,Uglov}, play an important role
  in what follows.  A self-contained proof of \autoref{T:FockSpace},
  stated with similar language, can be found in Ariki's
  book~\cite[Theorem~10.10]{Ariki:branching}.

  An element $x\in\Fock$ has \textbf{weight}
  $\wt(x)=\Gamma$ if
  $K_ix=q^{(\Gamma,\alpha_i)}x$, for~$i\in I$. In particular, if
  $\Zero=(0|0|\dots|0)\in\Parts[]$ is the
  empty multipartition of level~$\ell$ then
  $K_i\fock[\Zero]=q^{(\Lambda,\alpha_i)}\fock[\Zero]$, for $i\in I$, so that
  $\fock[\Zero]$ has weight $\Lambda$. More generally, if $\beta\in Q^+$
  then writing $\blam=\bmu+A$ it follows by induction that
  \begin{equation}\label{E:FockWeight}
    \text{if $\blam\in\Parts[\beta]$ then $d_i(\blam)=(\Lambda-\beta,\alpha_i)$,
          for all $i\in I$.}
  \end{equation}
  Therefore, $\wt(\fock)=\Lambda-\beta$ by \autoref{T:FockSpace}.
  Set $d_i(\beta)=(\Lambda-\beta,\alpha_i)$.

  For each dominant weight $\Lambda\in P^+$ let $L(\Lambda)=\Uq
  v_\Lambda$ be the irreducible integrable highest weight module of
  highest weight~$\Lambda$, where~$v_\Lambda$ is a highest weight vector
  of weight~$\Lambda$. By \autoref{T:FockSpace}, $\fock[\Zero]$ is a
  highest vector of weight~$\Lambda$ in~$\Fock$. In fact, it follows
  from \autoref{T:FockSpace} that
  $$L(\Lambda)\cong\Uq\fock[\Zero].$$
  For example, see \cite[Theorem~10.10]{Ariki:book}. Henceforth, we set
  $v_\Lambda=\fock[\Zero]$.

  To compare the Grothendieck groups $\repH$ and $\projH$ with the Fock
  space we need to consider all~$n\ge0$ simultaneously. Define
  $$\Rep[\A]=\bigoplus_{n\ge0}\repH\quad\text{and}\quad
    \Proj[\A]=\bigoplus_{n\ge0}\projH.
  $$
  Set $\Rep=\Rep[\A]\otimes_\A\Q(q)$ and $\Proj=\Proj[\A]\otimes_\A\Q(q)$.

   \begin{Proposition}\label{P:Fock}
    Suppose that $\Lambda\in P^+$. Then the $i$-induction and $i$-restriction
    functors of $\Rep$ induce a $\Uq$-module structure on $\Proj$ and
    $\Rep$ such that, as $\Uq$-modules,
    $$\Proj\cong L(\Lambda)\cong \Rep.$$
  \end{Proposition}

  \begin{proof}
    Recall that $\dec$ is the graded decomposition matrix of~$\H$ and
    $\dec^T$ is its transpose. Abusing notation slightly by simultaneously
    using these matrices for all $n\ge0$, define linear maps
    $$\begin{tikzpicture}
     \matrix[matrix of math nodes,row sep=1cm,column sep=16mm]{
     |(P)| \Proj &|(F)|\Fock\\
                 &|(R)|\Rep\\
     };
     \draw[right hook->](P) -- node[above]{$\dec^T$} (F);
     \draw[->>](F) -- node[right]{$\dec$} (R);
     \draw[->](P) -- node[below]{$\car$} (R);
    \end{tikzpicture}$$
    where $\dec^T([Y^\bmu])=\sum_\blam d_{\blam\bmu}(q)\fock$,
    $\dec(\fock)=\sum_\bmu d_{\blam\bmu}(q)[D^\bmu]$ and where
    $\car=\dec\circ\dec^T$ is the Cartan map. As vector
    space homomorphisms, $\dec^T$ is injective and $\dec$ is surjective.
    As defined these maps are only vector space homomorphisms, however,
    we claim that both maps can be made into~$\Uq$-module homomorphisms.

    The $i$-induction and $i$-restriction functors are exact, for $i\in
    I$, because they are exact when we forget the
    grading~\cite[Corollary8.9]{Groj:control}.  Therefore,
    they send projective modules to projectives and they induce
    endomorphisms of the Grothendieck groups~$\Rep$ and~$\Proj$. By
    \autoref{T:Induction},
    \begin{align*}
      [\iRes S^\blam]&= \sum_{B\in\Rem_i(\blam)}q^{d_B(\blam)}[S^{\blam-B}],\\
      [\iInd S^\blam\<1-d_i(\blam\>]
        &=\sum_{A\in\Add_i(\blam)}q^{d_A(\blam)+1-d_i(\blam)}[S^{\blam+A}]\\
        &=\sum_{A\in\Add_i(\blam)}q^{-d^A(\blam)}[S^{\blam+A}],
    \end{align*}
    where the last equality uses \autoref{L:defect}(a).
    Identifying $E_i$ with $\iRes$, and~$F_i$ with $q\iInd K_i^{-1}$, the linear
    maps~$\dec$ and~$\dec^T$ become well-defined $\Uq$-module
    homomorphisms by \autoref{T:FockSpace}. As $\Uq$-modules,
    $\Rep$ and $\Proj$ are both cyclic because they are both generated
    by~$[Y^\Zero]=[S^\Zero]=[D^\Zero]$. By definition,
    $\dec^\T([Y^\Zero])=\zero$ and $\dec(\zero)=[S^\Zero]$, so the
    proposition follows because $L(\Lambda)\cong\Uq\zero$ is
    irreducible.
  \end{proof}

  Let $\UA$ be Lusztig's $\A$-form of $\Uq$, which is the
  $\A$-subalgebra of~$\Uq$ generated by the quantised divided powers
  $E_i^{(k)}=E_i^k/\diag k!$ and~$F_i^{(k)}=F_i^k/\diag k!$, for
  $i\in I$ and~$k\ge0$. \autoref{T:FockSpace} implies that $\UA$ acts on the
  $\A$-submodule~$\Fock[\A]$ of~$\Fock$; compare with
  \cite[Lemma~6.2]{LLT} and with \cite[Lemma~6.16]{M:Ulect}.  Therefore, by
  \autoref{P:Fock}, $\Rep[\A]$ and~$\Proj[\A]$ are $\UA$-modules.
  Moreover, if we set $L_\A(\Lambda)=\UA v_\Lambda$ then there are
  $\UA$-module homomorphisms $\Proj[\A]\hookrightarrow
  L_\A(\Lambda)\twoheadrightarrow \Rep[\A]$.  In particular,
  $L_\A(\Lambda)\cong\Proj[\A]$ as $\UA$-modules by \autoref{P:Fock}.

  \begin{Lemma}\label{L:UqCommute}
    Suppose that $i\in I$. The involution $\circledast$
    commutes with the actions of $E_i$ and $F_i$ on $\Rep[\A]$ and
    on~$\Proj[\A]$.
  \end{Lemma}

  \begin{proof}
    By \autoref{P:CommutingDuals}, there are isomorphisms
    $\iRes\circ\circledast\cong\circledast\circ \iRes$ and
    $\iInd\circ\#\cong\#\circ\iInd$. In particular, the actions
    of~$E_i=\iRes$ and~$\circledast$ commute. Fix $\beta\in Q^+$. Recall from after
    \autoref{E:FockWeight} that $d_i(\beta)=(\Lambda-\beta,\alpha_i)$.
    Identifying~$F_i$ with the functor
    $q\circ\iInd\circ K_i^{-1}=q^{1-d_i(\beta)}\iInd$ on $\rep(\H[\beta])$,
    there are isomorphisms
    \begin{align*}
      F_i\circ\circledast
        &\cong q\iInd K^{-1}_i\circ q^{-2\defect\beta}\#
           && \text{by \autoref{L:SameDuals}},\\
        &\cong q^{1-d_i(\beta)-2\defect\beta}\iInd\circ\#\\
        &\cong q^{1-d_i(\beta)-2\defect\beta}\#\circ\iInd
           && \text{by \autoref{P:CommutingDuals}},\\
        &\cong q^{-2\defect(\beta+\alpha_i)}\#\circ q^{d_i(\beta)-1}\circ\iInd
           &&\text{by \autoref{L:defect}(a),}\\
        &\cong\circledast\circ q^{-1}\iInd K_i\cong\circledast\circ F_i,
           &&\text{by \autoref{L:SameDuals}}.
    \end{align*}
    Hence, $E_i$ and $F_i$ commute with $\circledast$ on $\Rep[\A]$ and $\Proj[\A]$
    as claimed.
  \end{proof}

  In contrast,~$E_i$ and~$F_i$ on $\Rep[\A]$ and $\Proj[\A]$ do not
  commute with~$\#$.

  We want to relate the Cartan pairing $\<\ ,\ \>$ on
  $\Proj[\A]\times\Rep[\A]$ with the representation theory of
  $L_\A(\Lambda)$. Define a non-degenerate symmetric bilinear
  form $(\ ,\ )$ on the Fock space $\Fock[\A]$ by
  \begin{equation}\label{E:FockPairing}
    (\fock,\fock[\bmu])=\delta_{\blam\bmu}q^{\defect\blam},
      \qquad \text{for }\blam,\bmu\in\Parts[].
  \end{equation}
  By \autoref{T:FockSpace}, $(K_ix,y)=q^{\wt(y)}(x,y)=(x,K_iy)$, for
  weight vectors $x,y\in\Fock[\A]$ and $i\in I$.  By restriction,
  we also consider $(\ ,\ )$ as a bilinear form on~$L_\A(\Lambda)$.

  \begin{Lemma}\label{L:ContrBiadjoint}
    The bilinear form $(\ ,\ )$ on $L_\A(\Lambda)$ is characterised by
    the three properties: $(\zero,\zero)=1$,
    $(E_ix,y)=(x,F_iy)$ and $(F_ix,y)=(x,E_iy)$, for all $i\in I$ and
    $x,y\in L_\A(\Lambda)$.
  \end{Lemma}

  \begin{proof}
     By definition, $(\zero,\zero)=1$. If $i\in I$ then in order to check
     that~$E_i$ and~$F_i$ are biadjoint with respect to $(\ ,\ )$ it is enough
     to consider the cases when $x=\fock$ and $y=\fock[\bmu]$, for
     $\blam,\bmu\in\Parts[]$. By \autoref{T:FockSpace},
     $(F_i\fock,\fock[\bmu])=0=(\fock,E_i\fock[\bmu])$ unless $\bmu=\blam+A$ for some
     $A\in\Add_i(\blam)$. On the other hand, if~$A\in\Add_i(\blam)$ and
     $\bmu=\blam+A$ then, using \autoref{L:defect}(a) for the second equality,
     \begin{align*}
       (F_i\fock,\fock[\bmu])
                &=q^{\defect\bmu-d^A(\blam)}
                 =q^{\defect\blam+d_i(\blam)-1-d^A(\blam)}\\
                &=q^{\defect\blam+d_A(\bmu)}
                 =(\fock,E_i\fock[\bmu]).
     \end{align*}
    A similar calcuation shows that
    $(E_i\fock,\fock[\bmu])=(\fock,F_i\fock[\bmu])$, for all
    $\blam,\bmu\in\Parts$.  As~$\zero$ is the
    highest weight vector in the irreducible module $L_\A(\Lambda)$, these
    three properties uniquely determine the bilinear form $(\ ,\ )$
    on~$L_\A(\Lambda)$ by induction on weight.
  \end{proof}

  By restriction, the next result categorifies the pairing $(\ ,\ )$ on
  $L_\A(\Lambda)$.

  \begin{Proposition}\label{P:FormComparision}
      Let $x\in\Proj[\A]$ and $y\in\Fock[\A]$ with
      $\wt(\blam)=\beta$. Then
      $$\<x,\dec(y)\> = q^{-\defect\beta}(\dec^T(x^\#),y).$$
  \end{Proposition}

  \begin{proof}
    As $(\ ,\ )$ is bilinear and $\<\ ,\ \>$ is sesquilinear
    it is enough to verify this identity when $x=X_\bmu$ and $y=\fock$,
    for $\bmu\in\Klesh[]$ and $\blam\in\Parts[\beta]$. Then
    \begin{align*}
      \<x^\#,\dec(y)\>&=\<X_\bmu,[S^\blam]\>
         =\sum_{\bsig\in\Klesh[\beta]}d_{\blam\bsig}(q)\<X_\bmu,[D^\bsig]\>\\
        &=\sum_{\btau\in\Klesh[\beta]}\sum_{\bsig\in\Klesh[\beta]}
            d_{\blam\bsig}(q)\overline{e_{\bsig\btau}(-q)}\<X_\bmu,
               [S^\btau]^\circledast\>\\
        &=\sum_{\bsig\in\Klesh[\beta]}
            d_{\blam\bsig}(q)\overline{e_{\bsig\bmu}(-q)},
    \end{align*}
    where the last equality uses \autoref{L:DualBases}. For the right
    hand side,
    \begin{align*}
      (\dec^T(x^\#),y) &= \(\dec^T(X_\bmu^\#),\fock[\blam]\)
          =\sum_{\bsig\in\Klesh[\beta]}\overline{e_{\bsig\bmu}(-q)}
               \(\dec^T([Y^\bsig]),\fock\)\\
         &=\sum_{\substack{\bnu\in\Parts\\\bnu\gedom\bmu}}\Big(
       \sum_{\substack{\bsig\in\Klesh\\\bnu\gedom\bsig\gedom\bmu}}
       d_{\btau\bsig}(q)\overline{e_{\bsig\bmu}(-q)}\Big)
             \(\fock[\btau],\fock\)\\
        &=q^{\defect\beta}\<x,\dec(y)\>,
    \end{align*}
    by \autoref{E:FockPairing} and calculation above. The proof is
    complete.
  \end{proof}

  Now we can prove the results that we are really interested in.

  \begin{Corollary}\label{P:biadjointness}
      Let  $P\in\proj(\H)$, $y\in\rep(\H[n+1])$,and $i\in I$. Then
      $$ \<\iInd x,y\>=\<x,\iRes y\> \text{ and }
         \<\iRes x,y\>=\<x,\iInd y\>
      $$
  \end{Corollary}

  \begin{proof}
    By \autoref{T:Biadjoint}, $(\iRes,\iInd)$ is a biadjoint pair so the
    corollary follows directly from the definition of the Cartan pairing
    in~\autoref{E:CartanPairing}.  As it is non-trivial to show that
    $\iRes$ is left adjoint to $\iInd$ we prove this at the level of
    Grothendieck groups. Write $\dot x=\dec^T(x)$ and $y=\dec(\dot y)$
    where $\dot x,\dot y\in L_\A(\Lambda)$ and $\wt(\dot
    y)=\Lambda-\beta$. Then $\<\iRes x,y\>=0$ unless
    $\wt(x)=\Lambda-(\beta+\alpha_i)$. To improve readability, identify
    $x$ and $\dot x=\dec^T(x)$ below. Then,
    \begin{align*}
      \<\iRes x,y\> &=q^{-\defect\beta}\( (E_ix)^\#,\dot y\)
                   &&\text{by \autoref{P:FormComparision}},\\
             &=q^{\defect\beta}\(E_ix^\circledast,\dot y\),
             &&\text{by \ref{L:SameDuals} and \ref{L:UqCommute}},\\
             &=q^{\defect\beta}\(x^\circledast,F_i\dot y\),
                  &&\text{by \autoref{L:ContrBiadjoint}},\\
             &=q^{\defect\beta-2\defect(\beta+\alpha_i)}\(x^\#,F_i\dot y\),
                  &&\text{by \autoref{L:SameDuals}},\\
             &=q^{\defect\beta-\defect(\beta+\alpha_i)}\<x^\#,F_iy\>,
                  &&\text{by \autoref{P:FormComparision}},\\
             &=\<x,\iInd y\>,
    \end{align*}
    where the last equality uses \autoref{L:defect} and the identification
    of~$F_i$ and~$q\iInd\circ K^{-1}_i$ on $\Rep[\A]$, via
    \autoref{P:Fock}.
  \end{proof}

  Let $\tau$ be the unique semilinear anti-isomorphism of~$\Uq$ such
  that $\tau(K_i)=K_i^{-1}$, $\tau(E_i)=qF_iK_i^{-1}$ and
  $\tau(F_i)=q^{-1}K_iE_i$, for all $i\in I$. Then the biadjointness of
  induction and restriction with respect to the Cartan pairing translates
  into the following more Lie theoretic statement.

  \begin{Corollary}\label{C:TauAdjoint}
    Suppose that $x\in\Proj[\A]$ and  $y\in\Rep[\A]$. Then
    $$\<ux,y\>=\<x,\tau(u)y\>,\qquad\text{ for all }u\in\UA.$$
  \end{Corollary}

  The bar involution of~$\A$ extends to a semilinear involution
  of~$\UA$ determined by $\overline{K_i}=K_i^{-1}$, $\overline{E_i}=E_i$
  and $\overline{F_i}=F_i$, for all $i\in I$. Similarly, define a bar
  involution on~$L_\A(\Lambda)$ by
  $$\overline{v}_\Lambda=\zero\quad\text{and}\quad
    \overline{ux}=\overline{u}\,\overline{x},
    \qquad\text{for }u\in\UA\text{ and }x\in L_\A(\Lambda).
  $$
  As noted in \cite[\S3.1]{BK:GradedDecomp}, it follows from the
  relations that $\tau\circ\overline{\ \ }=\overline{\ \ }\circ\tau^{-1}$.

  As in \cite[\S3.3]{BK:GradedDecomp}, the \textbf{Shapovalov form} on $L(\Lambda)$
  is the sesquilinear map
  $$\<x,y\>=q^{\defect\beta}(\overline{x},y),$$
  for $x,y\in L(\Lambda)$ with $\wt(y)=\Lambda-\beta$, for~$\beta\in Q^+$.
  As our notation suggests, the Shapovalov form is categorified by
  the Cartan pairing.

  \begin{Corollary}\label{C:ShapovalovCartan}
    Suppose that $x\in\Proj[\A]$ and $y\in L_\A(\Lambda)$. Then
    $$\<\dec^T(x),y\>=\<x,\dec(y)\>.$$
  \end{Corollary}

  \begin{proof}By \autoref{L:ContrBiadjoint}, the pairing $\(\ ,\ )$ on
    $L_\A(\Lambda)$ is unique symmetric bilinear map on $L_\A(\Lambda)$
    that is biadjoint with respect to $E_i$ and $F_i$ and such that
    $(v_\Lambda,v_\Lambda)=1$. This implies that the Shapovalov form is
    the unique sesquilinear form on $L_\A(\Lambda)$ such that
    $\<v_\Lambda,v_\Lambda\>=1$ and $\<ux,y\>=\<x,\tau(u)y\>$, for
    $x,y\in L_\A(\Lambda)$ and $u\in\UA$. Hence, the result follows from
    \autoref{C:TauAdjoint}.
  \end{proof}

  The module $L_\A(\Lambda)=\UA v_\Lambda$ is the \textit{standard
  $\A$-form} of the irreducible $\Uq$-module $L(\Lambda)$. The
  \textit{costandard $\A$-form} of $L(\Lambda)$ is the dual lattice
  \begin{align*}
    L_\A(\Lambda)^*&= \set{y\in L(\Lambda)|(x,y)\in\A\text{ for all }
                                                 x\in L_\A(\Lambda)}\\
                   &= \set{y\in L(\Lambda)|\<x,y\>\in\A\text{ for all }
                                                                x\in L_\A(\Lambda)}
  \end{align*}
  We can now identify both $\Proj[\A]$ and $\Rep[\A]$ as $\UA$-modules.

  \begin{Corollary}
    Suppose that $\Lambda\in Q^+$. Then, as $\UA$-modules,
    $$\Proj[\A]\cong L_\A(\Lambda)
    \qquad\text{ and }\qquad\Rep[\A]\cong L_\A(\Lambda)^*.$$
  \end{Corollary}

  \begin{proof}
    The first isomorphism we noted already after \autoref{P:Fock}. The
    second isomorphism follows from \autoref{C:ShapovalovCartan}
    and \autoref{L:DualBases}.
  \end{proof}

  By \autoref{L:UqCommute}, the action of $F_i$ on $\Rep[\A]$ and
  $\Proj[\A]$, for $i\in I$, commutes with $\circledast$. In the
  language of \cite[\S3.1]{BK:GradedDecomp}, $\circledast$ is a
  \textit{compatible bar-involution}. As is easily proved by induction
  on weight, every integrable $\UA$-module has a unique
  bar-compatible involution, so
  \begin{align}\label{E:SameBar}
      \dec(\overline{y})=\dec(y)^\circledast\qquad\text{for all }y\in\Fock.
  \end{align}
  It follows that $\set{B^\bmu|\bmu\in\Klesh[]}$ is
  Kashiwara's \textbf{upper global basis} at $q=0$~\cite{Kashiwara:CrystalBases},
  or Lusztig's \textbf{dual canonical basis}~\cite[\S14.4]{L:can},
  of~$L(\Lambda)$.  By \autoref{C:SecondDuals},
  $q^{-\defect\bmu}B_\bmu$ is bar invariant and, thinking of~$(\ ,\ )$
  as a pairing from $\Proj[A]\times\Rep[\A]$ to $\A$, we have
  $$(q^{-\defect\bmu}B_\bmu,B^\blam) = \<B_\bmu^\#,B^\blam\>
          =\<B_\bmu, B^\blam\>=\delta_{\blam\bmu},$$
  by \autoref{P:FormComparision} and \autoref{C:ShapovalovCartan}.
  Hence, $\set{q^{-\defect\bmu}B_\bmu|\bmu\in\Klesh[]}$ is the
  \textbf{canonical basis}, or the \textbf{upper global basis}, of~$L(\Lambda)$.


  We could have proved the equivalence of parts (a)--(d) of the next
  result before we introduced the quantum group~$\Uq$. For~(e), however,
  we need Kashiwara's theory of crystal
  bases~\cite{Kashiwara:CrystalBases,Kashiwara:GlobalCrystalBases} and
  the work of Misra and Miwa~\cite{MisraMiwa} relate crystal
  bases of the Fock space and crystal bases of~$L(\Lambda)$.

  \begin{samepage}
  \begin{Proposition}\label{P:Positivity}
    Suppose that $\F$ is an arbitrary field and that $n\ge0$. Then
    the following are equivalent:
    \begin{enumerate}\topskip=0pt
      \item For all $\bmu\in\Klesh$, $B^\bmu=[D^\bmu]$.
      \item For all $\blam,\bmu\in\Klesh$,
                $e_{\blam\bmu}(-q)\in\delta_{\blam\bmu}+q\N[q]$.
      \item For all $\bmu\in\Klesh$, $B_\bmu=[Y^\bmu]$.
      \item For all $\blam,\bmu\in\Klesh$,
                $d_{\blam\bmu}(q)\in\delta_{\blam\bmu}+q\N[q]$.
      \item For all $\blam\in\Parts$ and $\bmu\in\Klesh$,
                $d_{\blam\bmu}(q)\in\delta_{\blam\bmu}+q\N[q]$.
    \end{enumerate}
  \end{Proposition}
  \end{samepage}

  \begin{proof}In the Grothendieck groups, $[D^\bmu] = [S^\bmu] +
    \sum_{\bmu\gdom\btau}e_{\bmu\tau}(-q)[S^\btau]$ and $[Y^\bmu] =
    X_\bmu + \sum_{\bmu\gdom\bsig}d_{\bsig\bmu}(q)X_\bsig$, where in the
    sums $\bsig,\btau\in\Klesh$.  Moreover, by
    \autoref{L:BarInvariance}, $[Y^\bmu]^\#=[Y^\bmu]$ and
    $[D^\bmu]^\circledast=[D^\bmu]$, for all $\bmu\in\Klesh$. By definition,
    $d_{\blam\bmu}(q)\in\N[q,q^{-1}]$ and
    $e_{\blam\bmu}(-q)\in\Z[q,q^{-1}]$. Hence, parts~(a) and~(b), and
    parts~(c) and~(d), are equivalent by \autoref{P:CanonicalBasis}. Moreover,
    $\dinv=(\deck)^{-1}$, $d_{\bmu\bmu}(1)=1=e_{\bmu\bmu}(-q)$ and the
    Laurent polynomials $d_{\blam\bmu}(q)$ and
    $e_{\blam\bmu}(-q)$ are non-zero only if $\blam\gedom\bmu$ by
    \autoref{T:CellularSimples}, so parts~(b) and~(d) are also
    equivalent. Certainly, (e) implies (d) so to complete the proof it
    is enough to show that (a) implies (e).

    Suppose that~(a) holds so that $B_\bmu=D^\bmu$, for all
    $\bmu\in\Klesh[]$. To prove that~(e) holds we need the machinery of
    \textit{crystal
    bases}~\cite{Kashiwara:CrystalBases,Kashiwara:GlobalCrystalBases} in
    the special case of the Fock space~$\Fock$. We will refer the reader
    to the literature for the definitions and results that we need.

    Following~\cite[\S2]{Kashiwara:GlobalCrystalBases} define rings
    $\AA=\Q[q,q^{-1}]=\Q\otimes_\Z\A$, $\AA_0=\AA_{(q)}$ and
    $\AA_\infty=\AA_{(q^{-1})}$, so that $\AA_0$ and $\AA_\infty$ are
    the rational functions in~$\Q(q)$ that are regular at~$0$
    and~$\infty$, respectively. Set
    $$ L_0(\Lambda)=\bigoplus_{\bmu\in\Klesh[]}\AA_0B^\bmu
                   =\bigoplus_{\bmu\in\Klesh[]}\AA_0[S^\bmu]$$
    and $B_0(\Lambda)=\set{[S^\bmu]+qL_0(\Lambda)|\bmu\in\Klesh[]}$.
    As $\{B^\bmu\}$ is the upper crystal basis, the
    pair $\(L_0(\Lambda),B_0(\Lambda)\)$ is an upper crystal
    base at $q=0$ for $L(\Lambda)$ as defined by
    Kashiwara~\cite[\S2]{Kashiwara:CrystalBases}. Similarly, in the Fock space
    define
    $$\Fock[0]=\bigoplus_{\blam\in\Parts[]}\AA_0\fock
    \quad\text{and}\quad
      C^\Lambda_0=\set{\fock+q\Fock[0]|\blam\in\Parts[]}.
    $$
    Misra and Miwa~\cite{MisraMiwa} showed that
    $(\Fock[0],C^\Lambda_0)$ is an upper crystal basis for~$\Fock$. (As
    we discuss below they explicitly described the crystal graph of
    $\Fock$.) By \cite[Theorem~7]{Kashiwara:CrystalBases} the Fock space
    $\Fock$ has a unique basis $\set{C^\blam|\blam\in\Parts[]}$, Kashiwara's
    upper global basis, such that
    \begin{equation}\label{E:CanonicalBasis}
    \overline{C^\blam}=C^\blam\text{ and }
      C^\blam\equiv\fock\quad\pmod{q\Fock[0]},
      \qquad\text{ for }\blam\in\Parts[].
    \end{equation}
    Let $\Fock[\bullet]=\(\Fock[\AA],\Fock[0],\overline{\Fock[0]}\)$ and
    $L_\bullet(\Lambda)=\(L_\AA(\Lambda),L_0(\Lambda),L_0(\Lambda)^\circledast\)$,
    where $\Fock[\AA]=\AA\otimes_{\AA_0}\Fock[0]$ and
    $L_\AA(\Lambda)=\AA\otimes_{\AA_0}L_0(\Lambda)$. Then
    $\Fock[\bullet]$ and $L_\bullet(\Lambda)$ are
    balanced triples in the sense of
    \cite[\S2]{Kashiwara:GlobalCrystalBases}. The $\Lambda$-weight space
    of~$\Fock$ is~$\Q(q)\zero$ so, up to a scalar, the decomposition map
    $\dec$ is the unique $\Uq$-module homomorphism $\dec\map{\Fock}\Rep$.
    By \cite[Proposition~5.2.1]{Kashiwara:CrystalBases},
    the image of~$\Fock[\bullet]$ under~$\dec$ is a balanced
    triple contained in~$L(\Lambda)$. In fact, we have
    $\dec(\Fock[\bullet])=L_\bullet(\Lambda)$ by
    \cite[Proposition~5.2.2]{Kashiwara:CrystalBases} because~$\dec$
    sends $\zero=\fock[\Zero]$ to~$[S^{\Zero}]$. Consequently,
    if~$\blam\in\Parts[]$ then $\dec(\fock)\in L_0(\Lambda)$. That is,
    $$[S^\blam] = \dec(\fock) \in\bigoplus_{\bmu\in\Klesh[]}\AA_0[S^\bmu]
                   =\bigoplus_{\bmu\in\Klesh[]}\AA_0[D^\bmu].$$
    As $[S^\blam]=\sum_\bmu d_{\blam\bmu}(q)[D^\bmu]$ it follows that
    $d_{\blam\bmu}(q)\in\N[q,q^{-1}]\cap\AA_0=\N[q]$. Moreover, because
    of \autoref{E:SameBar}, $\dec$ sends canonical basis elements
    in~$\Fock[0]$ to canonical basis elements in $L_0(\Lambda)$, or to
    zero. It follows that
    $$\dec(C^\blam)=\begin{cases}
          B^\blam=[D^\blam],&\text{if }\blam\in\Klesh[]\\
          0,&\text{otherwise}.
    \end{cases}$$
    By \autoref{E:CanonicalBasis}, $\fock-C^\blam\in q\Fock[0]$ for
    all $\blam\in\Parts[]$.
    Consequently, if $\blam\notin\Klesh[]$ then
    $$[S^\blam]=\dec(\fock)=\dec(\fock-C^\blam)\in \dec(q\Fock[0])
           =qL_0(\Lambda).$$
    Hence, $d_{\blam\bmu}(q)\in\delta_{\blam\bmu}+q\N[q]$, for all
    $\blam\in\Parts[]$ and all $\bmu\in\Klesh[]$. Thus, (e) holds and
    the proposition is proved.
  \end{proof}

  \begin{Remark}
    The difference between the upper and lower crystal bases, or the dual canonical
    and canonical bases, can be interpreted as changing between the
    bases of Specht modules and dual Specht modules. The global bases and
    their crystal lattices are:
    \begin{align*}
      \text{upper}&& q&=0 & B^\bmu&\equiv[S^\bmu]
          &\pmod{\textstyle\sum_{\blam\in\Klesh[]}\AA_0[S^\blam]}\\
      \text{lower}&& q&=\infty &
      \car(q^{-\defect\bmu}B_\bmu)&\equiv[S_{\Mull(\bmu)'}]
          &\pmod{\textstyle\sum_{\blam\in\Klesh[]}\AA_\infty[S_\blam]}
    \end{align*}
    where $\Mull$ is an
    involution on~$\Klesh$ that generalises the well-known Mullineux map
    for the symmetric groups. See \autoref{T:Mullineux} below.
  \end{Remark}

  \begin{Remark}
    As mentioned in \autoref{R:FockSpaceAction}, a
    different action on the Fock space is commonly used in the
    literature. With respect to the Cartan pairing, as
    in~\autoref{C:TauAdjoint}, this action is the adjoint of the action
    in \autoref{T:FockSpace}. As a consequence, the papers that use a
    different $\Uq$-action also use a different coproduct for~$\Uq$, as
    they have to if they want Kashiwara's tensor product rule to connect
    the crystal bases at different levels for a fixed~$\Lambda\in P^+$.
    In the dual set up, $\#$ categorifies the bar involution
    on~$L(\Lambda)$, $\set{B_\bmu|\bmu\in\Klesh[]}$ is the canonical
    basis, or lower global crystal basis at~$q=0$ for~$L(\Lambda)$ and
    $\set{q^{\defect\bmu}B^\bmu|\bmu\in\Klesh[]}$ is its dual canonical
    basis.
  \end{Remark}

  It is natural to ask when the equivalent conditions of
  \autoref{P:Positivity} are satisfied. In general, this is a difficult
  open problem. The next result shows that these properties hold
  whenever~$\F$ is a field of characteristic zero.

  We can now state Ariki's celebrated Categorification Theorem. By specializing
  $q=1$ the quantum group $\UA\otimes\Q$ becomes the Kac-Moody algebra $\Usl$.
  Let $L_1(\Lambda)$ be the irreducible integrable highest weight $\Usl$-module
  of high weight~$\Lambda$. The canonical bases of $L_1(\Lambda)$ are obtained
  by specializing $q=1$ in the canonical bases of~$L_\A(\Lambda)$. Forgetting
  the grading in the results above, $\UnRep\cong L_1(\Lambda)\cong\UnProj$,
  where $\UnRep=\bigoplus_n\rep(\UnH)\otimes_\Z\Q$ and
  $\UnProj=\bigoplus_n\proj(\UnH)\otimes_\Z\Q$.

  \begin{Theorem}[\protect{Ariki's Categorification Theorem~\cite{Ariki:can,BK:DegenAriki}}]
     Suppose that $\F$ is a field of characteristic zero. Then the
     canonical basis of $L_1(\Lambda)$ coincides with the basis of
     $($ungraded$)$ projective indecomposable $\H$-modules
     $\set{[\underline{Y}^\bmu]|\bmu\in\Klesh[]}$ of\/ $\UnProj$.
     \label{T:Ariki}
  \end{Theorem}

  This theorem was proved by Ariki~\cite[Theorem~4.4]{Ariki:can} when $v^2\ne1$
  and by Brundan and Kleshchev when $v^2=1$ \cite[Theorem~3.10]{BK:DegenAriki}.
  For a detailed proof of this important result when $v^2\ne1$ see
  \cite[Theorem~12.5]{Ariki:book}. For an overview and historical account of
  Ariki's theorem see \cite{Geck:Ariki}.

  Combining \autoref{T:Ariki} with \autoref{P:Positivity} we obtain the main
  result of this section.

  \begin{samepage}
  \begin{Corollary}[\protect{%
    Brundan and Kleshchev~\cite[Theorem~5.14]{BK:GradedDecomp}}]
    Suppose that $\F$ is a field of characteristic zero. Then the canonical basis
    of $L_\A(\Lambda)$ coincides with the basis
    $\set{q^{-\defect\bmu}[Y^\bmu]|\bmu\in\Klesh[]}$
    of\/ $\Proj$. In particular, $d_{\blam\bmu}(q)\in\delta_{\blam\bmu}+q\N[q]$,
    for all $\blam\in\Parts[]$ and $\bmu\in\Klesh[]$.
    \label{C:GradedDecomp}
  \end{Corollary}
  \end{samepage}

  When $\Lambda$ is a weight of level~$2$ and $e=\infty$ this was first proved by
  Brundan and Stroppel~\cite[Theorem~9.2]{BrundanStroppel:KhovanovIII}. For
  extensions of this result to cyclotomic quiver Hecke algebras of arbitrary type
  see \cite{KangKashiwara:CatCycKLR,LaudaVazirani:CatCrystals,Rouquier:QuiverHecke2Lie,Webster:HigherRep}.

  \autoref{C:GradedDecomp} implies that the graded decomposition numbers
  $d_{\blam\bmu}(q)=[S^\blam:D^\bmu]_q=b_{\blam\bmu}(q)$ are
  \textit{parabolic Kazhdan-Lusztig polynomials}. Explicit formulas are
  given in \cite[Appendix~A]{HuMathas:QuiverSchurI} and
  \cite[Lemma~2.46]{Maksimau:QuiverSchur}.

  For the canonical basis $\{B_\bmu\}$ it is immediate that the Laurent
  polynomials $b_{\blam\bmu}(q)\in\Z[q]$ are polynomials, for
  $\blam\in\Parts$ and $\bmu\in\Klesh$, however, it is a deep fact that
  their coefficients are \textit{non-negative} integers. In contrast, it
  is immediate that $d_{\blam\bmu}(q)\in\N[q,q^{-1}]$  but it is a deep
  fact that the graded decomposition numbers are polynomials rather than
  \textit{Laurent polynomials}.  Thus, the difficult result changes from
  positivity of coefficients in the ungraded setting, to positivity of
  exponents in the graded setting. In fact, it is also true when
  $\F=\C$ that the inverse graded decomposition numbers
  $e_{\blam\bmu}(q)=b_{\blam\bmu}(-q)$ are polynomials in~$q$ with
  non-negative integer coefficients. This is perhaps best explained by
  passing to the Koszul dual of the corresponding graded cyclotomic
  Schur
  algebras~\cite{HuMathas:QuiverSchurI,StroppelWebster:QuiverSchur,Ariki:GradedSchur}
  using \cite{HuMathas:QuiverSchurI,Maksimau:QuiverSchur}; see
  \cite[Lemma~2.15]{HuMathas:QuiverSchurI} where this is stated
  explicitly.

  Brundan and Kleshchev's proof of \autoref{C:GradedDecomp} is quite
  different to the one given here. They have to work quite hard to
  define triangular bar involutions on~$L_\A(\Lambda)$ whereas we have done
  this by exploiting the representation theory of~$\H$. One benefit of
  Brundan and Kleshchev's approach is that they have an explicit
  description of the bar involution on~$\Fock[\A]$. In contrast, we
  have no hope of working with our bar involution unless we already know
  the graded decomposition matrices.  On the other hand, the approach
  here works for an arbitrary multicharge~$\charge$.

  To complete the proof of \autoref{C:GradedDecomp}, Brundan and
  Kleshchev lift Grojnowski's approach~\cite{Groj:control} to the
  representation theory of $\H$ to the graded setting. As a result they
  obtain graded analogues of Kleshchev's modular branching
  rules~\cite{Klesh:III,Klesh:IV,Brundan:KM}.  Under categorification,
  these branching rules correspond to the action of the crystal
  operators on the crystal graph of~$L(\Lambda)$; see
  \cite[Theorem~4.12]{BK:GradedDecomp}. By invoking Ariki's theorem they
  deduce an analogue of \autoref{C:GradedDecomp}, although with a
  possibly different labelling of the irreducible modules. Finally, they
  prove that the labelling of the irreducible $\H$-modules coming from
  the branching rules agrees with the labelling of
  \autoref{C:UngradedSimples}; compare
  with~\cite{Ariki:branching,ArikiJaconLecouvey:FactorizationCanonicalBaes}.

  We have not yet given an explicit description of the labelling of the
  irreducible $\H$-modules because we defined
  $\Klesh=\set{\bmu\in\Parts|\UnD^\bmu\ne0}$. Extending \autoref{E:dA},
  if $\bmu\in\Parts$ and given nodes $A<C$ define
  $$d_A^C(\bmu) = \#\set{B\in\Add_i(\bmu)|A{<}B{<}C}
                      -\#\set{B\in\Rem_i(\bmu)|A{<}B{<}C}.$$
  Following Misra and Miwa~\cite{MisraMiwa} (but using Kleshchev's
  terminology~\cite{Klesh:II}), a removable $i$-node~$A$ is \textbf{normal} if
  $d_A(\bmu)\le0$ and $d_A^C(\bmu)<0$ whenever $C\in\Rem_i(\bmu)$ and $A<C$.  A normal
  $i$-node $A$ is \textbf{good} if $A\le B$ whenever~$B$ is a normal $i$-node.
  Write $\blam\goodarrow\bmu$ if~$\mu=\lambda{+}A$ for some good node $A$. Misra
  and Miwa~\cite[Theorem~3.2]{MisraMiwa} show that the crystal graph of
  $L_\A(\Lambda)$, considered as a submodule~$\Fock[\A]$, is the graph with
  vertex set $$\Cry=\set{\bmu\in\Parts[]|\bmu=\zero\text{ or }
  \blam\goodarrow\bmu\text{ for some $\blam\in \Cry$}}, $$ and with
  labelled edges $\blam\iarrow\bmu$ whenever $\bmu$ is obtained from $\blam$ by
  adding a good $i$-node, for some $i\in I$. See
  \cite[Theorem~11.11]{Ariki:book} for a self-contained proof of this
  result, couched in similar language.

  \begin{Corollary}[\protect{Ariki~\cite{Ariki:class,BK:GradedDecomp,%
             ArikiJaconLecouvey:ModularBranching}}]\label{C:SimpleClass}
    Suppose that $\F$ is an arbitrary field and that $\bmu\in\Parts$. Then
    $\Klesh[]=\Cry$. That is, if $\bmu\in\Parts$ then $\UnD_\F^\bmu\ne0$ if
    and only if~$\bmu\in \Cry$.
  \end{Corollary}

  \begin{proof}If $\F$ is a field of characteristic zero then
    $\set{\bmu+qL_\A(\Lambda)|\bmu\in\Klesh[]}$
    is a basis of $L_\A(\Lambda)/qL_\A(\Lambda)$ by
    \autoref{P:Positivity} and \autoref{C:GradedDecomp}. This basis of
    $L_\A(\Lambda)/qL_\A(\Lambda)$ is exactly the crystal basis of
    $L(\Lambda)$ by \autoref{C:GradedDecomp}, so $\Klesh=\Cry$ in
    characteristic zero.  If $\F$ is a field of positive characteristic
    then a straightforward modular reduction argument shows that
    $\UnD_\F^\bmu\ne0$ only if $\UnD^\bmu_\C\ne0$, for $\bmu\in\Parts$
    (for example, see~\autoref{S:Adjustments} below). So,
    $\Klesh[]\subseteq \Cry$.  By \autoref{P:Fock}, the number of
    irreducible $\H$-modules depends only on~$e$, and not
    on the field~$\F$, so $\Klesh[]=\Cry$ as required.
  \end{proof}

  The idea in \autoref{P:Fock} that over any field the natural bases
  $\{[D^\bmu]\}$ and $\{[Y^\bmu]\}$ of~$\Rep[\A]$ and~$\Proj[\A]$, respectively,
  are distinguished bases of~$L(\Lambda)$ goes back to at least Lascoux, Leclerc
  and Thibon~\cite{LLT}. This was generalized to higher levels by
  Ariki~\cite{Ariki:can} and it played an important role in the classification
  of the irreducible $\H$-modules~\cite{Ariki:class,AM:simples} and in
  Grojnowski and Vazirani's
  work~\cite{Klesh:book,Groj:control,GrojVaz,Vazirani:Parametrisation}. The role
  of the crystal graphs in the representation theory of~$\R$ is explored further
  in \cite{KangPark:IrredKLRTypeA,LaudaVazirani:CatCrystals}.

  \subsection{Homogeneous Garnir relations}
  By \autoref{T:Zfree}, $\R$ is a graded cellular algebra and, as a consequence,
  that there exist graded lifts of the Specht modules for arbitrary $\Lambda\in
  P^+$. However, at this point we cannot really compute inside the graded Specht
  modules because we do not know how to write basis elements indexed by
  non-standard tableaux in terms of standard ones. This section shows how to do
  this. First, some combinatorics.

  Fix a multipartition~$\blam$ and a node $A=(l,r,c)\in\blam$. A (row)
  \textbf{Garnir node} of~$\blam$ is any node $A=(l,r,c)$ such that
  $(l,r+1,c)\in\blam$.  The \textbf{$(e,A)$-Garnir belt} is the set of nodes
  \begin{align*}
    \Belt=&\set{(l,r,k)\in\blam|k\ge c\text{ and }
            e\ceil{k-c+1}/e\leq\lambda^{(l)}_r-c+1}\\
      &\quad\cup\quad
         \set{(l,r+1,k)\in\blam|k\le c\text{ and }c\geq e\ceil{c-k+1}/e}.
  \end{align*}

  Let $b_A=\#\Belt/e$ and write $b_A=a_A+c_A$ where $ea_A$ is the number of nodes
  in~$\Belt$ in row $(l,r)$. Let $\D_A$ be the set of minimal length right coset
  representatives of~$\Sym_{a_A}\times\Sym_{c_A}$ in~$\Sym_{b_A}$; see, for
  example, \cite[Proposition~3.3]{M:Ulect}. When $e=\infty$ these definitions
  should be interpreted as $\Belt=\emptyset$, $b_A=0=a_A=c_A$ and $\D_A=1$.

  Suppose $A$ is a Garnir node of $\blam$. The rows of $\blam$ are indexed by
  pairs $(l,r)$, corresponding to row~$r$ in $\lambda^{(l)}$ where $1\le l\le\ell$
  and $r\ge1$. Order the row indices lexicographically. Let $\t_A$ be the
  $\blam$-tableau that agrees with~$\tlam$ for all numbers
  $k<\tlam(A)=\tlam(l,r,c)$ and $k>\tlam(l,r+1,c)$ and where the remaining entries
  in~rows $(l,r)$ and $(l,r+1)$ are filled in increasing order from left to
  right first along the nodes in row~$(l,r+1)$ that are in the first
  $c$~columns but not in~$\Belt$, then along the nodes in row~$(l,r)$ of $\Belt$
  followed by the nodes in row~$(l,r+1)$ of~$\Belt$, and then along the
  remaining nodes in row~$(l,r)$.

  \begin{Example}  As Garnir belts are contained in consecutive rows of the same
    component, the general case can be understood by looking at a two-rowed
    partition (of level one), so we consider the case $e=3$, $\lambda=(14,6)$ and
    $A=(1,1,4)$. Then
    $$\t_A=
    \begin{tikzpicture}[baseline=-3.6mm,scale=0.45,draw/.append style={thick,black}]
      \draw[preaction={fill=blue!30},pattern=bricks,pattern color=white]%
                  (6.5,0.5)rectangle ++(3,-1);
      \draw[preaction={fill=blue!30},pattern=bricks,pattern color=white]%
                  (1.5,-0.5)rectangle ++(3,-1);
      \draw[preaction={fill=red!30},pattern=bricks,pattern color=white]%
                  (3.5,0.5)rectangle ++(3,-1);
      \draw[preaction={fill=red!30},pattern=bricks,pattern color=white]%
                  (9.5,0.5)rectangle ++(3,-1);
      \newcount\col
      \foreach\Row/\row in {{1,2,3,5,6,7,8,9,10,11,12,13,17,18}/0,%
                            {4,14,15,16,19,20}/-1} {
         \col=1
         \foreach\k in \Row {
           \draw[thin](\the\col,\row)+(-.5,-.5)rectangle++(.5,.5);
           \draw[thin](\the\col,\row)node{\k};
           \global\advance\col by 1
          }
       }
       \draw[blue,ultra thick](6.5,0.5)rectangle ++(3,-1);
       \draw[blue,ultra thick](1.5,-0.5)rectangle ++(3,-1);
       \draw[red,ultra thick](9.5,0.5)rectangle ++(3,-1);
       \draw[red,ultra thick](3.5,0.5)rectangle ++(3,-1);
    \end{tikzpicture}$$
    The lines in~$\t_A$ show how the $(3,A)$-Garnir belt decomposes into a
    disjoint union of ``$e$-bricks''.  In general, $b_A$ is equal to the
    number of $e$-bricks in the Garnir belt and~$a_A$ is the number of $e$-bricks
    in its first row. In this case,~$b_A=4$ and~$a_A=3$.  Therefore,
    $\D_A=\{1,s_3,s_3s_2,s_3s_2s_1\}$.
  \end{Example}

  Let $k_A=\t_A(A)$ be the number occupying~$A$ in~$\t_A$ and define
  $$w^A_r=\prod_{a=k_A+e(r-1)}^{k_A+re-1}(a,a+e),$$
  for $1\le r<b_A$.  The subgroup $\<w^A_r\mid 1\le r<b_A\>$
  of~$\Sym_n$ is isomorphic to~$\Sym_{b_A}$ via the map $w^A_r\mapsto s_r$,
  for $1\le r<b_A$. Set $\bi^A=\bi^{\t_A}$ and $\tau^A_r=e(\bi^A)(\psi_{w^A_r}+1)$,
  for $1\le r<b_A$. If $d\in\D_A$ choose a reduced
  expression $d=s_{r_1}\dots s_{r_k}$ for $d$ and define
  $$\tau^A_d = \tau^A_{r_1}\dots\tau^A_{r_k}\in e(\bi^A)\R.$$
  The elements $\tau^A_d$ of $\R$ seem to be very special and deserving
  of further study.  They are homogeneous elements in~$\R$ of degree
  zero that are independent of all choices of reduced expressions.
  Moreover, by \cite[Theorem~4.13]{KMR:UniversalSpecht}, the elements
  $\set{\tau^A_r|1\le r<b_A}$ satisfy the braid relations when they act
  on~$S^\blam_\Zcal$ and they generate a copy of~$\Sym_{b_A}$
  inside~$\End_{\Zcal}(S^\blam_\Zcal)$!

  \begin{Theorem}[\protect{%
    Kleshchev, Mathas and Ram~\cite[Theorem~6.23]{KMR:UniversalSpecht}}]
    \label{T:Garnir}
    Suppose that $\blam\in\Parts$ and that $\Zcal$ is an integral domain. The
    graded Specht module $S_\Zcal^\blam$ of~$\R(\Zcal)$ is isomorphic to the graded
    $\R$-module generated by a homogeneous element $v_{\tlam}$ of
    degree~$\deg\tlam$ subject to the relations:
    \begin{enumerate}
      \item $v_{\tlam} e(\bi)=\delta_{\bi\bi^\blam}v_{\tlam}$.
      \item $v_{\tlam} y_s=0$, for $1\le s\le n$.
      \item $v_{\tlam} \psi_r=0$ whenever $s_r\in\Sym_\blam$,
      for $1\le r<n$.
      \item $\sum_{d\in\D_A}v_{\tlam}\psi_{\t_A}\tau^A_d=0$, for all Garnir nodes
      $A\in\blam$.
    \end{enumerate}
  \end{Theorem}

  Relations (a)--(c) already appear in \cite{BKW:GradedSpecht} and, in terms of
  the cellular basis machinery, they are a consequence of
  \autoref{P:PsiProperties}.

  The relations in part~(d) are the \textbf{homogeneous Garnir relations}.
  These relations are a homogeneous form of the well-known Garnir relations of
  the symmetric group~\cite[Theorem~7.2]{James}.  There is an analogous
  description of the dual Specht modules $S_\blam$ in terms of \textit{column
  Garnir relations}~\cite[\S7]{KMR:UniversalSpecht}. Using Dyck tilings,
  Fayers~\cite{Fayers:DyckTilings} has shown how to write the
  homogeneous Garnir relations in terms of $\psi$-basis of the Specht
  module.

  The most difficult part of the proof of \autoref{T:Garnir} is showing
  that the~$\tau^A_d$ satisfy the braid relations.  This is proved using
  the Khovanov-Lauda diagram calculus that was briefly mentioned in
  \autoref{S:CycQuiverHecke}.  Like \autoref{T:Zfree} this result holds
  over an arbitrary ring. To prove that the graded module defined by
  \autoref{T:Garnir} has the correct rank the construction of the graded
  Specht module~$S^\blam$ over a field in \autoref{T:PsiBasis},
  from~\cite{HuMathas:GradedCellular,BKW:GradedSpecht}, is used.

  One of the main points of \autoref{T:Garnir} is that it makes it possible to
  calculate in the graded Specht module over any ring.  Prior to
  \autoref{T:Garnir} the only way to compute inside the graded Specht modules
  was, in effect, to use the isomorphism $\R\bijection\H$ of \autoref{T:BKiso}
  to work in the ungraded setting then use the inverse isomorphism
  $\H\bijection\R$ to get back to the graded setting.  This made it difficult to
  keep track of, and to exploit, the grading on~$S^\blam$ --- and it was only
  possible to work with Specht modules defined over a field.

  \autoref{T:Garnir} also gives the relations for $S^\blam$ as an $\Rn$-module. From
  this perspective \autoref{T:Garnir} can be used to give another
  construction of the graded Specht modules. For $\alpha,\beta\in Q^+$ let
  $\Rn[\alpha,\beta]=\Rn[\alpha]\otimes\Rn[\beta]$. \autoref{D:QuiverRelations}
  implies that there is a non-unital embedding
  $\Rn[\alpha,\beta]\hookrightarrow\Rn[\alpha+\beta]$ that maps
  $e(\bi)\otimes e(\bj)$ to $e(\bi\vee\bj)$, where $\bi\vee\bj$ is the
  concatenation of~$\bi$ and $\bj$. Under this embedding
  the identity element of~$\Rn[\alpha,\beta]$ maps to
  $$e_{\alpha,\beta}=\sum_{\bi\in I^\alpha,\, \bj\in I^\beta}e(\bi\vee\bj).$$
  \autoref{D:QuiverRelations} implies that $\Rn[\alpha+\beta]$ is free as an
  $\Rn[\alpha,\beta]$-module, so the functor
  $$\Ind_{\alpha,\beta}^{\alpha+\beta}(M\boxtimes N)
  =(M\boxtimes N)e_{\alpha,\beta}\otimes_{\Rn[\alpha,\beta]}\Rn[\alpha+\beta]$$
  is a left adjoint to the natural restriction map. Iterating this
  construction, given $\beta_1,\dots,\beta_\ell\in Q^+$ and $\Rn[\beta_k]$
  modules $M_k$, for $1\le k\le\ell$, define
  $$M_1\circ\dots\circ M_\ell =
  \Ind_{\beta_1,\dots,\beta_\ell}^{\beta_1+\dots+\beta_\ell}
     (M_1\boxtimes\dots\boxtimes M_\ell).$$

  The definition of the graded Specht modules by generators and relations in
  \autoref{T:Garnir} makes the following result almost obvious. This description
  of the Specht modules is part of the folklore of these algebras with several
  authors~\cite{Vazirani:Parametrisation,BK:DegenAriki} using it as the
  definition of Specht modules.

  \begin{Corollary}[\protect{%
    Kleshchev, Mathas and Ram~\cite[Theorem~8.2]{KMR:UniversalSpecht}}]
    \label{TSpechtHigherLevel}
    Suppose that $\lambda^{(k)}\in\Pcal_{1,\beta_k}$, for $\beta_k\in Q^+$
    and $1\le k\le\ell$, so that $\blam\in\Parts[\beta]$, where
    $\beta=\beta_1+\dots+\beta_\ell$. Then
    there is an isomorphism of graded $\R$-modules $($and graded
    $\Rn$-modules$)$,
    $$S^\blam\<\deg\t^{\lambda^{(1)}}+\dots+\deg\t^{\lambda^{(\ell)}}\>
         \cong(S^{\lambda^{(1)}}\circ\dots
             \circ S^{\lambda^{(\ell)}})\<\deg\t^\lambda\>,
    $$
    where on the right hand side $S^{\lambda^{(k)}}$ is considered as an
    $\Rn[\beta_k]$-module, for $1\le k\le\ell$.
  \end{Corollary}

  A second application of \autoref{T:Garnir} is a generalization of James'
  famous result \cite[Theorem 8.15]{James} for symmetric groups that describes
  what happens to the Specht modules when they are tensored with the sign
  representation. First some notation.

  Following~\cite[\S3.3]{KMR:UniversalSpecht}, for $\bi\in I^n$ let
  $-\bi=(-i_1,\dots-i_n)\in I^n$. Recalling the multicharge $\charge$
  from \autoref{S:Quivers}, set $\charge'=(-\kappa_\ell,\dots,-\kappa_1)$ and
  let $\Lambda'=\Lambda(\charge')\in P^+$. Similarly, if
  $\beta=\sum_ia_i\alpha_i\in Q^+$ let $\beta'=\sum_{i\in I}a_i\alpha_{-i}$.
  Inspecting \autoref{D:CycQuiverHecke}, there is a unique isomorphism of graded
  algebras
  \begin{equation}\label{epsilon}
    \sgn\map{\R[\beta]}\Rp;\quad e(\bi)\mapsto e(-\bi),
           \quad y_r\mapsto -y_r, \quad\text{and}\quad \psi_s\mapsto-\psi_s,
  \end{equation}
  for all admissible $r$ and $s$ and $\bi\in I^\beta$. The involution $\sgn$
  induces an equivalence of categories
  $\rep(\Rp[\beta'])\longrightarrow\rep(\R[\beta])$ that sends an $\Rp$-module
  $M$ to the $\R[\beta]$-module $M^\sgn$, where the $\R[\beta]$-action is
  twisted by~$\sgn$.

  \begin{Corollary}[\protect{%
    Kleshchev, Mathas and Ram~\cite[Theorem~8.5]{KMR:UniversalSpecht}}]
    \label{C:SignDualSpecht}
    Suppose that $\bmu\in\Parts[\beta]$, for $\beta\in Q^+$. Then
    $$S^\bmu\cong(S_{\bmu'})^\sgn\quad\text{ and }\quad S_\bmu\cong(S^{\bmu'})^\sgn$$
    as $\R[\beta]$-modules.
  \end{Corollary}

  In \cite{KMR:UniversalSpecht} this is proved by checking the relations in
  \autoref{T:Garnir}. As noted in
  \cite[Proposition~3.26]{HuMathas:QuiverSchurI}, this can be proved more
  transparently by noting that, up to sign, the involution $\sgn$ maps the
  $\psi$-basis of~$\R$ to the $\psi'$-basis of~$\Rp$. Some care must be taken
  with the notation here. For example, if $\bmu\in\Parts[\beta]$ then
  $\bmu'\in\Parts[\beta']$. See \cite[\S3.7]{HuMathas:QuiverSchurI} for more
  details.

  We give an application of these results to the graded decomposition numbers.
  First, by \autoref{C:SimpleClass} if $\bmu\in\Klesh$ there exists $\bi\in I^n$
  and a sequence of multipartitions $\bmu_0=\zero,\bmu_1,\dots,\bmu_n=\bmu$ in
  $\Klesh[]$ such that $\bmu_{k+1}$ is obtained from~$\bmu_k$ by adding a good
  $i_k$-node, for $0\le k<n$. It follows from the modular branching rules
  \cite[Theorem 4.12]{BK:GradedDecomp}, and properties of crystal graphs, that
  there exists a unique sequence of multipartitions
  $\Mull(\bmu_0)=\zero, \Mull(\bmu_1), \dots, \Mull(\bmu_n)=\Mull(\bmu)$ such that
  $\Mull(\bmu_{k+1})$ is obtained from~$\Mull(\bmu_k)$ by adding a good
  $-i_k$-node and $\Mull(\bmu_{k+1})\in\mathcal{K}^{\Lambda'}_{k+1}$, for $1\le
  k\le n$. The \textbf{Mullineux conjugate} of $\bmu$ is the multipartition
  $\Mull(\bmu)$. Thus, $D^{\Mull(\bmu)}$ is a non-zero irreducible
  $\Rp$-module.  We emphasize that the $\Rp$-module $D^{\Mull(\bmu)}$ is defined
  using the $\psi$-basis of $\Rp$ and hence the crystal theory used
  in~\autoref{S:Categorification}, with respect to the multicharge $\charge'$.

  \begin{Theorem}\label{T:Mullineux}
    Suppose that $\bmu\in\Klesh[\beta]$, for $\beta\in Q^+$. Then
    $$(D^{\Mull(\bmu)})^\sgn\cong D^{\bmu}$$
    as $\R[\beta]$-modules.
  \end{Theorem}

  \begin{proof} As $\sgn$ is an equivalence of categories,
    $(D^{\Mull(\bmu)})^\sgn\cong D^\bnu\<d\>$ for some $\bnu\in\Klesh[\beta]$ and
    $d\in\Z$ by \autoref{C:GradedSimples}. Since $\sgn$ is homogeneous,
    by \autoref{T:GradedSimples}(a),
    $$\gdim (D^{\Mull(\bmu)})^\sgn=\gdim D^{\Mull(\bmu)}
                                  =\overline{\gdim D^{\Mull(\bmu)}}
                                  =\overline{\gdim (D^{\Mull(\bmu)})^\sgn},$$
    so that $d=0$ and $(D^{\Mull(\bmu)})^\sgn\cong D^\bnu$. To show that
    $\bnu=\bmu$ it is now enough to work in the ungraded setting. Therefore, we
    can either use the modular branching rules
    of~\cite{Ariki:branching,Groj:control}, or
    their graded counterparts from \cite[Theorem~4.12]{BK:GradedDecomp},
    together with what is by now a standard argument due to
    Kleshchev~\cite[Theorem~4.7]{Klesh:III}, to show that $\bnu=\bmu$.
  \end{proof}

  The $\sgn$ map induces an equivalence
  $\rep(\Rp)\longrightarrow\rep(\R[\beta])$. As~$\sgn$ is an involution, we also
  write $\sgn\map{\rep(\R[\beta])}\rep(\Rp)$ for the inverse equivalence. The
  last two results can now be written as
  $(S^\blam)^\sgn\cong S_{\blam'}$ and $(D^\bmu)^\sgn\cong D^{\Mull(\bmu)}$ as
  $\Rp$-modules, for $\blam\in\Parts[\beta]$ and $\bmu\in\Klesh[\beta]$.

  \begin{Corollary}\label{C:GradedDecompDegrees}
    Suppose that $F$ is a field and that $\blam\in\Parts[\beta]$ and
    $\bmu\in\Klesh[\beta]$. Then $d_{\bmu\bmu}(q)=1$,
    $d_{\Mull(\bmu)'\bmu}(q)=q^{\defect{\bmu}}$ and
    $d_{\blam\bmu}(q)\ne0$ only if
    $\Mull(\bmu)'\gedom\blam\gedom\bmu$. Moreover, if $F=\C$ then
    $0<\deg d^\C_{\blam\bmu}(q)<\defect\bmu$ whenever
    $\Mull(\bmu)'\gdom\blam\gdom\bmu$.
  \end{Corollary}

  \begin{proof}Suppose that $\blam\in\Parts[\beta]$ and $\bmu\in\Klesh[\beta]$.
    Then
    $$
      [S^\blam: D^\bmu]_q= [ (S^\blam)^\sgn: (D^\bmu)^\sgn ]_q
           = [ S_{\blam'} : D^{m(\bmu)} ]_q,
    $$
    by  \autoref{C:SignDualSpecht} and \autoref{T:Mullineux},
    respectively.  Therefore, using \autoref{C:DualSpechts} and
    \autoref{T:GradedSimples}(a),
    $$  [S^\blam: D^\bmu]_q
          = q^{\defect\blam}[(S^{\blam'})^\circledast:D^{m(\bmu)} ]_q,
           = q^{\defect\bmu}  \overline{ [ S^{\blam'} : D^{m(\bmu)} ]_q }.
    $$
    By \autoref{T:GradedSimples}(c), if $\btau\in\Klesh[\beta]$ and
    $\bsig\in\Parts[\beta]$ then $d_{\btau\btau}(1)=1$ and
    $d_{\bsig\btau}(q)\ne0$ only if $\bsig\gedom\btau$.  Therefore,
    $d_{\Mull(\bmu)'\bmu}(q)
        =q^{\defect\bmu}\,\overline{d_{\Mull(\bmu)\Mull(\bmu)}(q)}
        =q^{\defect\bmu}$
    and $d_{\blam\bmu}(q)\ne0$ only if $\Mull(\bmu)'\gedom\blam\gedom\bmu$.
    The argument so far is valid over any field. Now suppose that $F=\C$. Then
    $d_{\blam\bmu}(q)\in\delta_{\blam\bmu}+q\N[q]$, by \autoref{C:GradedDecomp},
    so the remaining statement about the degrees of the graded
    decomposition numbers follows.
  \end{proof}

  \autoref{C:GradedDecompDegrees} is the easy half of a conjecture of
  Fayers~\cite{Fayers:LLT}, which he was interested in because it leads
  to a faster algorithm for computing the graded decomposition numbers
  of~$\H$. At the level of canonical bases the last two results
  correspond to the fact shifting by the defect transforms an upper
  crystal base into a lower crystal
  base~\cite[Lemma~2.4.1]{Kashiwara:CrystalBases}. See also
  \cite[Remark~3.19]{BK:GradedDecomp}.


  \subsection{Graded adjustment matrices}\label{S:Adjustments}
  All of the results in this section have their origin in the work of
  James~\cite{James10} and Geck~\cite{Geck:Brauertrees} on
  \textit{adjustment matrices}. Brundan and Kleshchev have given two different
  approaches to graded decomposition matrices in~\cite[\S6]{BK:GradedKL} and
  \cite[\S5.6]{BK:GradedDecomp}.  In this section we give third cellular algebra
  approach. Even though our definitions and proofs are different, it is easy to
  see that everything in this section is equivalent to definitions or theorems of
  Brundan and Kleshchev --- or to graded analogues of results of James and Geck.

  Before we introduce the adjustment matrices, let $\A[I^n]$ be the free
  $\A$-module generated by $I^n$. The \textbf{$q$-character} of a
  finite dimensional  $\Rn$-module~$M$ is
  $$\Ch M=\sum_{\bi\in I^n} \gdim M_\bi\cdot \bi \in\A[I^n],$$
  where $M_\bi=Me(\bi)$, for $\bi\in I^n$. For example,
  $\Ch S^\blam=\sum_{\t\in\Std(\blam)}q^{\deg(\t)}\cdot\bi^\t$.

  \begin{Theorem}[\protect{\cite[Theorem~3.17]{KhovLaud:diagI}}]
    \label{T:InjectiveCh}
    Suppose that $\Zcal$ is a field. Then the map
    $$\Ch\map{[\rep(\Rn)]}\A[I^n]; [M]\mapsto\Ch M$$
    is injective.
  \end{Theorem}

  As every $\R$-module can be considered as an $\Rn$-module by inflation, it
  follows that the restriction of $\Ch$ to $[\rep(\R)]$ is still injective.
  Extend the map $\circledast$ to~$\A[I^n]$ by defining $\(\sum_\bi
  f_\bi(q)\cdot\bi\)^\circledast=\sum_\bi\overline{f_\bi(q)}\cdot\bi$.  Then
  $(\Ch[M])^\circledast=\Ch[M^\circledast]$, for all $M\in\rep(\R)$.

  This section compares representations of cyclotomic KLR algebras over
  different fields. Write $S^\blam_\Zcal$ and $D^\bmu_\Zcal$ to emphasize
  that these modules are $\R(\Zcal)$-modules, for $\blam\in\Parts$ and
  $\bmu\in\Klesh$. If $\Zcal=F$ is a field, and $K$ is an extension of~$F$, then
  $D^\bmu_K\cong D^\bmu_F\otimes_F K$ since $D^\bmu_F$ is absolutely irreducible
  by \autoref{T:GradedSimples}. Therefore $\Ch D^\bmu_F$ depends only on~$\bmu$
  and the characteristic of~$F$.

  By \autoref{T:Zfree}, or by
  \autoref{T:Garnir}, the graded Specht module $S^\bmu_\Z$ is defined over~$\Z$
  and $S^\bmu_\Zcal\cong S^\bmu_\Z\otimes_\Z\Zcal$ for any commutative ring~$\Zcal$.
  The graded Specht module $S^\bmu_\Z$ has basis
  $\set{\psi_\t|\t\in\Std(\bmu)}$ and  it comes equipped with a
  $\Z$-valued bilinear form $\<\ ,\ \>$ that is determined by
  \begin{equation}\label{E:form}
    \<\psi_\s,\psi_\t\>\psi_{\tlam}=\psi_\s\psi_{\t\tlam}
           =\psi_{\s}\psi_{d(\t)}^\star y^\bmu e(\ilam).
  \end{equation}
  Following \autoref{E:radical}, define the radical of~$S^\bmu_\Z$ to be
  $$\rad S^\bmu_\Z=\set{x\in S^\bmu_\Z|\<x,y\>=0\text{ for all }
          y\in S^\bmu_\Z}.$$
  In fact, by \autoref{E:form},
  $\rad S^\bmu_\Z=\set{x\in S^\bmu_\Z|xa=0\text{ for all
  }a\in(\R)^{\gedom\bmu}}$.

  \begin{Definition}
    Suppose that $\bmu\in\Parts$. Let $D^\bmu_\Z=S^\bmu_\Z/\rad S^\bmu_\Z$.
  \end{Definition}

  By definition, $\rad S^\bmu_\Z$ is a graded submodule of $S^\bmu_\Z$, so
  $D^\bmu_\Z$ is a graded $\R(\Z)$-module. Hence, $D^\bmu_\Z\otimes_\Z\Zcal$ is
  a graded $\R(\Zcal)$-module for any ring~$\Zcal$.

  The following result should be compared with \cite[Theorem~6.5]{BK:GradedKL}.

  \begin{Theorem}\label{T:DLattice}
    Suppose that $\bmu\in\Parts$. Then $\rad S^\bmu_\Z$ is a $\Z$-lattice in
    $\rad S^\bmu_\Q$ and $D^\bmu_\Z$ is a $\Z$-lattice in $D^\bmu_\Q$.
    Consequently, $D^\bmu_\Q=D^\bmu_\Z\otimes_\Z\Q$ and
    $\Ch D^\bmu_\Z=\Ch D^\bmu_\Q$.
  \end{Theorem}

  \begin{proof}
    Let $G^\bmu_\Z=(\<\psi_\s,\psi_\t\>)$ be the Gram matrix of~$S^\bmu_\Z$.
    As $\Z$ is a principal ideal domain, by the Smith normal form
    there exists a pair of bases $\{a_r\}$ and $\{b_s\}$ of~$S^\bmu_\Z$ such
    that $(\<a_r,b_s\>)=\mathop{\rm diag}(d_1,d_2,\dots,d_z)$ for some
    non-negative integers such that $d_1|d_2|\dots|d_z$, where $d_r=0$ only if
    $d_s=0$ for all $s\ge r$.  That is, $d_1,\dots,d_z$ are the elementary
    divisors of~$G^\bmu_\Z$. As the form is homogeneous, we may
    assume that the bases $\{a_r\}$ and $\{b_s\}$ are homogeneous with $\deg
    a_r=\deg\t_r=-\deg b_r$, for some ordering $\Std(\bmu)=\{\t_1,\dots,\t_z\}$.
    Moreover, in view of \autoref{P:PsiProperties}(a), we can also assume that
    $a_re(\bi)=\delta_{\bi^{\t_r},\bi}a_r$ and $b_s
    e(\bi)=\delta_{\bi^{\t_s},\bi}b_s$, for $1\le r,s\le z$ and $\bi\in I^n$.
    Comparing with the definitions above, it follows that $\set{a_r|d_r=0}$ is a
    basis of $\rad S^\bmu_\Z$ and that $\set{a_r+\rad S^\bmu_\Z|d_r\ne0}$ is a
    basis of $D^\bmu_\Z$. All of our claims now follow.
  \end{proof}

  For an arbitrary field $F$, it is usually not the case that $D^\bmu_F$ is
  isomorphic to $D^\bmu_\Z\otimes_\Z F$ as an $\R(F)$-module. Indeed, if~$F$ is a
  field of characteristic $p>0$ then the argument of \autoref{T:DLattice} shows
  that
  $$\dim_F D^\bmu_F=\set{1\le r\le z|d_r\not\equiv0\pmod p}
                   \le \text{rank}_\Z\  D^\bmu_\Z=\dim_\Q D^\bmu_\Q,$$
  with equality if and only if all of the non-zero elementary divisors of
  $G^\bmu_\Z$ are coprime to~$p$.

  \begin{Definition}[cf. \protect{%
    Brundan and Kleshchev~\cite[\S5.6]{BK:GradedDecomp}}]\label{D:Adjust}
    Suppose that $F$ is a field. For $\blam,\bmu\in\Klesh$ define Laurent
    polynomials $a^F_{\blam\bmu}(q)\in\N[q,q^{-1}]$ by
    $$a^F_{\blam\bmu}(q)
         =\sum_{d\in\Z}\,[D^\blam_\Z\otimes_\Z F: D^\bmu_F\<d\>]\, q^d.$$
    The matrix $\adj=\(a^F_{\blam\bmu}(q)\)$ is the \textbf{graded adjustment
    matrix} of~$\R(F)$.
  \end{Definition}

  Recall that $d_{\blam\bmu}(q)$ is a graded decomposition number of~$\R$.
  If we want to emphasize the field~$F$ then we write
  $d^F_{\blam\bmu}(q)=[S^\blam_F:D^\bmu_F]_q$ and
  $\dec^F=\(d^F_{\blam\bmu}(q)\)$. Note that $e$ is always fixed.

  \begin{Theorem}[\protect{cf. %
    Brundan and Kleshchev~\cite[Theorem~5.17]{BK:GradedDecomp}}]
    Suppose that $F$ is a field. Then:
    \label{T:Adjustment}
    \begin{enumerate}
      \item If $\blam,\bmu\in\Klesh$ then
      $a^F_{\blam\blam}(1)=1$ and $a^F_{\blam\bmu}(q)\ne0$ only if
      $\blam\gedom\bmu$. Moreover,
          $\overline{a^F_{\blam\bmu}(q)}=a^F_{\blam\bmu}(q)$.
      \item We have, $\dec^F=\dec^\Q\circ\adj$. That is, if $\blam\in\Parts$
      and $\bmu\in\Klesh$ then
      $$[S^\blam_F:D^\bmu_F]_q=d^F_{\blam\bmu}(q)
               =\sum_{\bnu\in\Klesh}d^\Q_{\blam\bnu}(q)a^F_{\bnu\bmu}(q).$$
    \end{enumerate}
  \end{Theorem}

  \begin{proof}By construction, every composition factor of $D^\blam_\Z\otimes
    F$ is a composition factor of $S^\blam_F$, so the first two properties of
    the Laurent polynomials $a^F_{\blam\bmu}(q)$ follow from
    \autoref{T:GradedSimples}.  By \autoref{T:DLattice}, the adjustment matrix
    induces a well-defined map of Grothendieck groups
    $\adj\map{[\rep(\R(\Q))]}[\rep(\R(F))]$ given by
    $$\adj\([D^\blam_\Q]\)=[D^\blam_\Z\otimes F]
                       =\sum_{\bmu\in\Klesh}a^F_{\blam\bmu}(q)[D^\bmu_F].$$
   Taking $q$-characters,
   $\Ch D^\blam_\Q=\sum_\bmu a^F_{\blam\bmu}(q)\Ch D^\bmu_F$. Applying
   $\circledast$ to both sides gives
   $\Ch D^\blam_\Q=\sum_\bmu \overline{a^F_{\blam\bmu}(q)}\Ch D^\bmu_F$. Therefore,
   $\overline{a^F_{\blam\bmu}(q)}=a^F_{\blam\bmu}(q)$ by \autoref{T:InjectiveCh},
   completing the proof of part~(a). For~(b), since
   $S^\blam_F\cong S^\blam_\Z\otimes_\Z F$,
   $$[S^\blam_F]=\adj\([S^\blam_\Q]\)
        =\adj\Big(\!\!\sum_{\bnu\in\Klesh} d^\Q_{\blam\bnu}(q)[D^\bnu_\Q]\Big)
        =\sum_{\bnu\in\Klesh}\sum_{\bmu\in\Klesh}
        d^\Q_{\blam\bnu}(q)a^F_{\bnu\bmu}(q)[D^\bmu_F].$$
    Comparing the coefficient of $[D^\bmu_F]$ on both sides completes the proof.
  \end{proof}

  \autoref{C:GradedDecomp} determines the graded decomposition numbers of the
  cyclotomic Hecke algebras in characteristic zero. There are several different
  algorithms for computing the graded decomposition numbers in characteristic
  zero~\cite{KleshchevNash:LLT,LLT,Fayers:LLT,HuMathas:QuiverSchurI,Uglov,GW}.
  To determine the graded decomposition numbers in positive characteristic it is
  enough to compute the adjustment matrices of \autoref{T:Adjustment}. The
  simplest case will be when $a^F_{\blam\bmu}(q)=\delta_{\blam\bmu}$, for all
  $\blam,\bmu\in\Klesh$.  Unfortunately, we currently have no idea when this
  happens. Two failed conjectures for when $\adj$ is the identity matrix are
  discussed in \autoref{Ex:KleshchevRam} and \autoref{Ex:JamesConjecture} below.

  We now compute the integral Gram matrices $G_\Z^\blam=\(\<\psi_\s,\psi_\t\>\)$
  and some adjustment matrix entries in several examples.

  \begin{Example}[Semisimple algebras]\label{Ex:SSSpechts}
    Suppose that $e>n$ and that $(\Lambda,\alpha_{i,n})\le1$, for all $i\in I$. Let
    $\blam\in\Parts$ and $\s,\t\in\Std(\blam)$. Then
    $\<\psi_\s,\psi_\t\>=\delta_{\s\t}$ because $\bi^\s=\bi^\t$ if and only if
    $\s=\t$ by \autoref{L:UniqueResidues}. Hence, $G_\Z^\blam$ is the identity
    matrix for all $\blam\in\Parts$.
  \end{Example}

  \begin{Example}[Nil-Hecke algebras]\label{Ex:NilHeckeSpechts}
    Suppose that $\Lambda=n\Lambda_i$ and $\beta=n\alpha_i$, for some $i\in I$.
    Let $\blam=(1|1|\dots|1)\in\Parts$, as in \autoref{S:NilHecke}, and suppose
    $\s,\t\in\Std(\blam)$ then
    $\<\psi_\s,\psi_\t\>\psi_{\tlam}
    =\psi_{\s}\psi_{d(\t)}^\star y_1^{n-1}y_2^{n-2}\dots y_{n-1},$
    by \autoref{E:form} and \autoref{Ex:NilHeckeylam}. By
    \autoref{P:SchubertPolynomials}, $\psi_\s\psi_{d(\t)}^\star=\psi_\u$,
    where $\u=\s d(\t)^{-1}$, if $\ell(d(\u))=\ell(d(\s))+\ell(d(\t))$ and
    otherwise $\psi_\s\psi_{d(\t)}^\star=0$. On the other hand, by the last
    paragraph of the proof of \autoref{P:SchubertPolynomials}, or simply by
    counting degrees, $\psi_\u y_1^{n-1}y_2^{n-2}\dots y_{n-1}=0$ if
    $\u\ne\t_\blam$ and
    $\psi_{\t_\blam}y_1^{n-1}y_2^{n-2}\dots y_{n-1}=(-1)^{n(n-2)/2}\psi_{\tlam}$.
    Hence, $\<\psi_\s,\psi_\t\>=\delta_{\s\t'}$, where
    $\t'=\t_\blam d'(\t)$ is the tableau that is conjugate to~$\t$. Hence,
    $G_\Z^\blam$ is $(-1)^{n(n-2)/2}$ times the anti-diagonal identity matrix.
    Consequently, $D^\blam_\Z=S^\blam_\Z$ and $S^\blam_F$ is irreducible
    for any field~$F$.
  \end{Example}

  \begin{Example}
    Suppose $e=2$, $\Lambda=\Lambda_0$ and $\lambda=(2,2,1)$. Then
    $\Std(\lambda)$ contains the five tableaux:
    $$\begin{array}{*6c}
             & \t_1=\t^\lambda & \t_2&\t_3&\t_4&\t_5\\\toprule
             \t&  \tableau({1,2},{3,4},{5}) &
                \tableau({1,3},{2,4},{5}) &
                \tableau({1,3},{2,5},{4}) &
                \tableau({1,2},{3,5},{4}) &
                \tableau({1,4},{2,5},{3}) \\
          d(\t)  & 1 & s_2 & s_2s_4 & s_4 &s_2s_4s_3\\
          \deg\t & 2&0&-2&0&0\\
          \bi^\t & 01100 & 01100 & 01100 & 01100 & 01010
      \end{array}
    $$
    We want to compute the Gram matrix $G^\blam_\Z=\(\<\psi_\s,\psi_\t\>\)$.
    Now $\<\psi_\t,\psi_\t\>\ne0$ only if~$\bi^\s=\bi^\t$, by
    \autoref{P:PsiProperties}(a), and only if $\deg\s+\deg\t=0$, since the
    bilinear form is homogeneous of degree zero. Hence, the only
    possible non-zero inner products are
    $$\<\psi_{\t_1},\psi_{\t_3}\>
           =\<\psi_{\t^\lambda},\psi_{\t^\lambda}\psi_2\psi_4\>
           =\<\psi_{\t^\lambda}\psi_4,\psi_{\t^\lambda}\psi_2\>
           =\<\psi_{\t_4},\psi_{\t_2}\>,$$
   together with $\<\psi_{\t_2},\psi_{\t_2}\>$,
   $\<\psi_{\t_4},\psi_{\t_4}\>$ and
   $\<\psi_{\t_5},\psi_{\t_5}\>$. If $a\in\{2,4\}$ then
   $$\<\psi_{\t^\lambda}\psi_a,\psi_{\t^\lambda}\psi_a\>
         =\<\psi_{\t^\lambda}\psi_a^2,\psi_{\t^\lambda}\>
         =\pm\<\psi_{\t^\lambda}(y_a-y_{a+1}),\psi_{\t^\lambda}\>
         =0,$$
    because $\psi_{\t^\lambda}y_r=0$ by \autoref{E:quadratic}, for $1\le r\le 5$.
    To compute the remaining inner products we have to go back to the definition
    of the bilinear form~\autoref{E:form}. By \autoref{D:psis},
    $y^\lambda=y_2y_4$ so
    $$
      \<\psi_{\t_1},\psi_{\t_3}\>\psi_{\t^\lambda}
            =\psi_{\t^\lambda}\psi_2\psi_4y_2y_4
            =\psi_{\t^\lambda}\psi_2y_2\psi_4y_4
            =\psi_{\t^\lambda}(y_3\psi_2+1)(y_5\psi_4+1)
            =\psi_{\t^\lambda},
    $$
    by \autoref{P:PsiProperties}(c). Hence,
    $\<\psi_{\t_1},\psi_{\t_3}\>=1=\<\psi_{\t_2},\psi_{\t_4}\>$. Finally,
    $$
      \<\psi_{\t_5},\psi_{\t_5}\>\psi_{\tlam}
         =\psi_{\tlam}\psi_2\psi_4\psi_3^2\psi_2\psi_4y_2y_4
         =\psi_{\tlam}\psi_2\psi_4(2y_3y_4-y_3^2-y_4^2)
               \psi_2\psi_4y_2y_4.
    $$
    where the second equality uses \autoref{E:quadratic}.
    Now $v_{\tlam}\psi_2y_3=v_{\tlam}(y_2\psi_1+1)=v_{\tlam}$ and, similarly,
    $v_{\tlam} y_4\psi_4=-v_{\tlam}$. Consequently $v_{\tlam} \psi_2\psi_4 y_a^2=0$,
    for $a=3,4$, so it follows that $\psi_{\tlam}\psi_2\psi_4\psi_3^2=-2\psi_{\tlam}$
    and hence that $\<\psi_{\t_5},\psi_{\t_5}\>=-2$. Therefore, the Gram matrix
    of $S^{(2,2,1)}$ is
    $$G_\Z^\lambda=\left(\begin{smallmatrix}
                            0 & 1 & 0 & 0 & 0\\
                            1 & 0 & 0 & 0 & 0\\
                            0 & 0 & 0 & 1 & 0\\
                            0 & 0 & 1 & 0 & 0\\
                            0 & 0 & 0 & 0 & -2
                          \end{smallmatrix}\right)
    $$
    Consequently, the elementary divisors of $G_\Z^\lambda$ are $1,1,1,1,2$.
    Therefore, if~$e=2$  and $\Zcal=\Q$ (so $v=-1$), then
    $S^\lambda_\Q=D^\lambda_\Q$ is irreducible, as is easily checked
    using \autoref{C:IrredSpechts}. Now suppose that $\Zcal=\Fp[2]$ (so
    $v=1$), so that $\H\cong\Fp[2]\Sym_5$. Then the calculation of
    $G_\Z^\lambda$ shows that the Specht module $S^\lambda$ is reducible
    with $\dim_{\Fp[2]} D^\lambda_{\Fp[2]}=4<5=\dim_\Q D^\lambda_\Q$.  It
    follows that if $e=p=2$ then $D^{(1^5)}$ is also a composition
    factor of $S^\blam$, so $a^{\Fp[2]}_{(2,2,1),(1^5)}(q)=1$.
  \end{Example}

  \begin{Example}\label{Ex:KleshchevRam}
    Kleshchev and Ram~\cite[Conjecture~7.3]{KleshchevRam:IrredKLR} made a
    conjecture that, in type~$A$, is equivalent to saying that the adjustment
    matrices $\adj$ of the (cyclotomic) KLR algebras are trivial when $e=\infty$.
    Williamson~\cite{Williamson:KLRcounterexample} has given an example that
    shows that, in general, this is not true.  Williamson's example comes from
    geometry~\cite{KashiwaraSaito:GeometricCrystals},  however, when it is
    translated into the language that we are using it corresponds to a
    statement about the simple module $D^\bmu$, for $\bmu=(2|2|1|1|3|3|2|2)$,
    for the cyclotomic quiver Hecke algebra $\R[16]$ with $e=\infty$ and
    $\Lambda=2\Lambda_1+2\Lambda_2+2\Lambda_3+2\Lambda_4$.  Fix the multicharge
    $\charge=(4, 4, 3, 3, 2, 2, 1, 1)$ and set
    $\bi=(4,5,3,4,2,3,4,5,2,3,1,2,3,4,1,2)$. So $y^\bmu=y_1y_9y_{15}y_{19}$.
    There are $5$ standard $\bmu$-tableaux of degree zero with residue
    sequence~$\bi$, namely:
    $$\begin{array}{c@{\,}*7{c@{\,|\,}}c@{\,}cc}
        \multicolumn{10}{c}{\t} & \ell(d(\t))\\\toprule
          \big( &\ootab(4,8)& \ootab(1,2)& \ootab(13)& \ootab(3)& \ootab(5,6,7)&
               \ootab(9,10,14)& \ootab(15,16)& \ootab(11,12)& \big) & 23 \\[1.5mm]
          \big( &\ootab(1,2)& \ootab(4,8)& \ootab(13)& \ootab(3)& \ootab(9,10,14)&
               \ootab(5,6,7)& \ootab(15,16)& \ootab(11,12)& \big) & 28 \\[1.5mm]
          \big( &\ootab(4,8)& \ootab(1,2)& \ootab(13)& \ootab(3)& \ootab(9,10,14)&
               \ootab(5,6,7)& \ootab(11,12)& \ootab(15,16)& \big) & 28 \\[1.5mm]
          \big( &\ootab(4,8)& \ootab(1,2)& \ootab(10)& \ootab(3)& \ootab(9,13,14)&
               \ootab(5,6,7)& \ootab(15,16)& \ootab(11,12)& \big) & 31 \\[1.5mm]
          \big( &\ootab(4,8)& \ootab(1,2)& \ootab(3)& \ootab(13)& \ootab(9,10,14)&
               \ootab(5,6,7)& \ootab(15,16)& \ootab(11,12)& \big) & 31
     \end{array}$$
    The negative oGram matrix for $\bi$-weight space $S^\bmu e(\bi)$
    of~$S^\bmu$ is
    $$\left(\begin{smallmatrix}
              0 & 0 & 1 & 1 & 0 \\
              0 & 0 & 1 & 1 & 0 \\
              1 & 1 & 0 & 1 & 1 \\
              1 & 1 & 1 & 0 & 1 \\
              0 & 0 & 1 & 1 & 0
            \end{smallmatrix}\right).
    $$
    Calculating this matrix is non-trivial because the lengths of the
    permutations $d(\t)$ are reasonably large. This matrix was computed using
    the author's implementation of the graded Specht modules in
    \textsf{Sage}~\cite{SageCombinat}.  Brundan, Kleshchev and McNamara
    \cite[Example~2.16]{BrundanKleshchevMcNamara:FiniteKLR} obtain exactly the
    same matrix, up to a permutation of the rows and columns, as part of the
    Gram matrix for the homogeneous bilinear form of the corresponding proper
    standard module for~$\Rn$.

    The elementary divisors of this matrix are $1, 1, 2, 0, 0$, so the dimension
    of~$D^\bmu e(\bi)$ is~$2$ in characteristic 2 and~$3$ in all other
    characteristics.  Consequently, the dimension of $D^\bmu$, and hence the
    adjustment matrix $\adj$ for $\R[16](F)$, depends on the characteristic
    of~$F$. That $\dim D\bmu_F$ depends on~$F$ was first proved by
    Williamson who computed a one dimensional \textit{intersection form}
    coming from geometry.
  \end{Example}

  \begin{Example}\label{Ex:JamesConjecture}
    Consider the case when $\Lambda=\Lambda_0$, so that $\H$ is the
    Iwahori-Hecke algebra of the symmetric group. The \textit{James
    conjecture}~\cite[\S4]{James10} says that if $F$ is a field of
    characteristic~$p>0$ and $\lambda,\mu\in\Parts$ then
    $a_{\lambda\mu}(q)=\delta_{\lambda\mu}$ if $ep>n$. A natural strengthening
    of this conjecture is that the adjustment matrix of $\R[\beta]$ is trivial
    whenever $\defect\beta<p$. For the symmetric groups, the condition
    $\defect\beta<p$ exactly corresponds to the case when the defect group of
    the block~$\R[\beta]$ is abelian.

    The James conjecture is known to be true for blocks of weight at
    most~$4$~\cite{James10,Richards,Fayers:weight3,Fayers:weight4}. Moreover,
    for every defect $w\ge0$ there exists a \textit{Rouquier block} of
    defect~$w$ for which the James conjecture holds~\cite{JLM}. Starting from
    the Rouquier blocks, there was some hope that the derived equivalences of
    Chuang and Rouquier~\cite{ChuangRouq:sl2} could be used to prove the
    James conjecture for all blocks.

    Notwithstanding all of the evidence in favour of the James conjecture, it
    turns out that the conjecture is wrong!  Again,
    Williamson~\cite[\S6]{Williamson:JamesLusztig} has cruelly (but ultimately
    kindly) produced counter-examples to the James conjecture. At the same time
    he also found counter-examples to the \textit{Lusztig
    conjecture}~\cite{Lusztig:conj} for $\text{SL}_n$. These examples rely upon
    Williamson's recent work with Elias that gives generators and relations for
    the category of Soergel bimodules~\cite{EliasWilliamson:Schubert}. As of
    writing, the smallest known counter-example to the James conjecture occurs in
    a block of defect $561$ in~$\Fp[839]\Sym_{467874}$. Williamson has not
    revealed which Specht modules his counter-examples appear in, so the size of Gram
    matrix that needs to be computed in order to verify this example is not
    known.  The Gram matrices of the Specht modules will be
    significantly larger, and harder to compute, than the one
    dimensional intersection form that Williamson reduces to (using a
    chain of deep results in geometric representation theory), and then
    calculates, using elementary techniques (and a computer).
  \end{Example}

  Williamson's counter-examples to the James and Lusztig conjectures
  suggest that there is no block theoretic criterion for the adjustment
  matrix of a block to be trivial, except asymptotically where the
  Lusztig conjecture is known to hold~\cite{AJS}. With hindsight,
  perhaps this is not so surprising because the condition given in
  \autoref{C:IrredSpechts} for a Specht module to be irreducible is
  rarely a block invariant. The failure of the James and Lusztig
  conjectures suggests that we should, instead, look for necessary and
  sufficient conditions for the $\R(F)$-modules $D^\bmu_\Z\otimes F$ to
  be irreducible, for $\bmu\in\Klesh$. Some steps towards such a
  criterion are made in \autoref{Conjecture} below.

  \medskip
  Brundan and Kleshchev \cite[\S5.6]{BK:GradedDecomp} remarked that
  $a^F_{\blam\bmu}(q)\in\N$ in all of the examples that they had computed.
  They asked whether this might always be the case. The next examples show
  that, in general, $a^F_{\lambda\mu}(q)\notin\N$.

  \begin{Example}[\protect{Evseev~\cite[Corollary~5]{Evseev:BadAdjustment}}]
    \label{Ex:Evseev}
    Suppose that $e=2$, $\Lambda=\Lambda_0$ and let $\lambda=(3,2^2,1^2)$ and
    $\mu=(1^9)$. Take $F=\Fp[2]$ to be a field of characteristic~$2$ and let
    $\adj[{\Fp[2]}]=(a^{\Fp[2]}_{\lambda\bmu}(q))$ be the corresponding
    adjustment matrix.

    As part of a general argument Evseev shows that
    $a^{\Fp[2]}_{\lambda\mu}(q)\notin\N$. In fact, it is not hard to see directly
    that $a^{\Fp[2]}_{\lambda\mu}(q)=q+q^{-1}$.
    Comparing the decomposition matrix for $\Fp[2]\Sym_9$ given by
    James~\cite{James} with the graded decomposition matrices when $e=2$ given
    in~\cite{M:Ulect}, shows that $d^\Q_{\lambda\mu}=0$,
    $d^{\Fp[2]}_{\lambda\mu}=2$, and that $a^{\Fp[2]}_{\lambda\mu}(1)=2$. Now
    $D_{\Fp[2]}^\mu=D_{\Fp[2]}^\mu e(\bi^\mu)$ is one dimensional, so any
    composition factor of $S_{\Fp[2]}^\lambda$ that is isomorphic to
    $D_{\Fp[2]}^\mu\<d\>$, for some $d\in\Z$, must be contained in
    $S_{\Fp[2]}^\lambda e(\bi^\mu)$. There are exactly
    six standard $\lambda$-tableau with residue sequence $\bi^\mu$, namely:
    $$\begin{array}{*7c}
      \deg\t&1&1&1&1&1&-1\\
      \t &
      \tableau({1,6,9},{2,7},{3,8},{4},{5})&
      \tableau({1,4,9},{2,5},{3,8},{6},{7})&
      \tableau({1,4,5},{2,7},{3,8},{6},{9})&
      \tableau({1,2,3},{4,7},{5,8},{6},{9})&
      \tableau({1,4,7},{2,5},{3,6},{8},{9})&
      \tableau({1,4,7},{2,5},{3,8},{6},{9})
    \end{array}$$%
    As $D^\mu$ is one dimensional, and concentrated in degree zero, it follows
    that $a^{\Fp[2]}_{\lambda\mu}=d^{\Fp[2]}_{\lambda\mu}(q)=q+q^{-1}$. We can see
    a shadow of the adjustment matrix entry in the Gram matrix of
    $S^\blam_\Z e(\bi^\mu)$, that is equal to
    {\scriptsize
      $$\left( \begin{smallmatrix}
                  0 & 0 & 0 & 0 & 0 & 0 \\
                  0 & 0 & 0 & 0 & 0 & 0 \\
                  0 & 0 & 0 & 0 & 0 & 4 \\
                  0 & 0 & 0 & 0 & 0 & -2 \\
                  0 & 0 & 0 & 0 & 0 & 2 \\
                  0 & 0 & 4 & -2 & 2 & 0
                \end{smallmatrix} \right)
      $$%
    }%
    The elementary divisors of this matrix are $2, 2, 0, 0, 0, 0$, with the
    $2$'s in degrees~$\pm1$. Therefore, the graded dimension of
    $D_{\Fp[2]}^\lambda e(\bi^\mu)$ decreases by~$q+q^{-1}$ in
    characteristic~$2$.
  \end{Example}

  \begin{Example}
    Motivated by the runner removable theorems
    of~\cite{JM:equating,ChuangMiyachi:runnerMorita} and \autoref{Ex:Evseev},
    take $e=3$, $F=\Fp[2]$, $\lambda=(3,2^4,1^3)$ and $\mu=(1^{14})$. (The
    partitions $\lambda$ and~$\mu$ are obtained from the corresponding
    partitions in \autoref{Ex:Evseev} by conjugating, adding an empty runner,
    and then conjugating again.) Again, we work over $\Fp[2]$ and consider the
    corresponding adjustment matrices.

    Calculating with \textsc{Specht}~\cite{SPECHT} we find that
    $d^\Q_{\lambda\mu}=0$ and that $d^{\Fp[2]}_{\lambda\mu}=2$.  Once again, it
    turns out that there are six $\lambda$-tableaux with $3$-residue
    sequence $\bi^\mu$, with five of these having degree~$1$ and one having
    degree~$-1$. (Moreover, the Gram matrix of $S^\lambda e(\bi^\mu)$ is
    the same as the Gram matrix given in \autoref{Ex:Evseev}.) Hence, as in
    \autoref{Ex:Evseev},
    $a^{\Fp[2]}_{\lambda\mu}(q)=q+q^{-1}=d^{\Fp[2]}_{\lambda\mu}(q)$.

    As the runner removable theorems compare blocks for different~$e$ over the
    same field we cannot expect to find an example of a non-polynomial
    adjustment matrix entry in odd characteristic in this way. Nonetheless, it
    seems fairly certain that non-polynomial adjustment matrix entries
    exist for all~$e$ and all~$p>0$.

    Evseev~\cite[Corollary~5]{Evseev:BadAdjustment} gives three other examples
    of adjustment matrix entries that are equal to $q+q^{-1}$ when $e=p=2$. All
    of them have similar analogues when $e=3$ and $p=2$.
    Finally, if we try adding further empty runners to the partitions $\lambda$
    and~$\mu$, so that $e\ge4$, then the corresponding adjustment matrix entry is
    zero (all of these partitions have weight~$4$).
  \end{Example}

\section{Seminormal bases and the KLR grading}
  In this final section we link the KLR grading on~$\R$ with the semisimple
  representation theory of~$\H$ using the seminormal bases. We start by showing
  that by combining information from all of the KLR gradings for different
  cyclic quivers leads to an integral formula for the Gram determinants of the
  ungraded Specht modules.

  \subsection{Gram determinants and graded dimensions}\label{S:GramDetsII}
  In \autoref{T:SchaperDet} we gave a ``rational'' formula for the Gram
  determinant of the ungraded Specht modules $\UnS^\blam$, for
  $\blam\in\Parts$. We now give an integral formula for these determinants and
  give both a combinatorial and a representation theoretic interpretation of
  this formula.

  Suppose that the Hecke parameter~$v$ from \autoref{D:HeckeAlgebras} is an
  indeterminate over~$\Q$ and consider an integral cyclotomic Hecke algebra $\H$
  over the field $\Zcal=\Q(v)$ where $\Lambda\in P^+$ such that $e>n$
  and $(\Lambda,\alpha_{i,n})\le 1$, for all $i\in I$. Then $\H$ is semisimple by
  \autoref{C:IntegralSSimple}.

  \begin{Definition}
    Suppose that $\blam\in\Parts$. For $e\ge2$ and
    $\bi\in I^n_e$ define
    $$\deg_{e,\bi}(\blam) = \sum_{\t\in\Std_\bi(\blam)}\deg_e\t,$$
    where $\Std_\bi(\blam)=\set{\t\in\Std(\blam)|\bi^\t=\bi}$.
    Set $\deg_e(\blam)=\sum_{\bi\in I^n_e}\deg_{e,\bi}(\blam)$.
    For a prime integer~$p>0$ set $\Deg_p(\blam)=\sum_{k\ge1}\deg_{p^k}(\blam)$.
  \end{Definition}

  By definition, $\deg_e(\blam),\Deg_p(\blam)\in\Z$. For $e>0$ let $\Phi_e(x)\in\Z[x]$
  be the $e$th cyclotomic polynomial in the indeterminate~$x$.

  \begin{Theorem}[\protect{Hu-Mathas~\cite[Theorem~C]{HuMathas:SeminormalQuiver}}]
    \label{T:IntegralDet}
    Suppose that $\Lambda\in P^+$, $e>n$ and that $(\Lambda,\alpha_{i,n})\le1$,
    for all $i\in I$. Let $\blam\in\Parts$. Then
    $$\det\UnG^\blam=\prod_{e>1}\Phi_e(v^2)^{\deg_e(\blam)}.$$
    Consequently, if $v=1$ then
    $\displaystyle\det\UnG^\blam=\prod_{p\text{ prime}}p^{\Deg_p(\blam)}.$
  \end{Theorem}

  Proving this result is not hard: it amounts to interpreting \autoref{D:SNCS}
  in light of the KLR degree functions on~$\Std(\blam)$. There is a power of~$v$
  in the statement of this result in~\cite{HuMathas:SeminormalQuiver}. This is
  not needed here because we have renormalised the quadratic relations in the
  Hecke algebra given in \autoref{D:HeckeAlgebras}.

  The Murphy basis is defined over $\Z[v,v^{-1}]$. Therefore,
  $\det\UnG^\blam\in\Z[v,v^{-1}]$ and \autoref{T:IntegralDet} implies that
  $\deg_e(\blam)\ge0$ for all $\blam\in\Parts$ and $e\ge2$. In fact,
  \cite[Theorem~3.24]{HuMathas:SeminormalQuiver} gives an analogue of
  \autoref{T:IntegralDet} for the determinant of the Gram matrix restricted to
  $\UnS^\blam e(\bi)$, suitably interpreted, and the following is true:

  \begin{Corollary}[\protect{\cite[Corollary~3.25]{HuMathas:SeminormalQuiver}}]
    \leavevmode\newline
    Suppose that $e\ge2$, $\blam\in\Parts$ and
    $\bi\in I_e^n$. Then $\deg_{e,\bi}(\blam)\ge0$.
  \end{Corollary}

  The definition of the integers $\deg_{e,\bi}(\blam)$ is purely
  combinatorial, so it should be possible to give a combinatorial proof
  of this result perhaps using \autoref{T:Induction}. We think, however,
  that this is probably difficult.

  Fix an integer $e\ge2$ and a dominant weight
  $\Lambda\in P^+$ and consider the Hecke algebra $\H$ over a field~$F$. If
  $\blam\in\Parts$ then, by definition,
  $$\Ch S^\blam = \sum_{\bmu\in\Klesh}d_{\blam\bmu}(q)\Ch D^\bmu\in\A[I^n].$$
  Let $\partial\map{\A[I^n]}\Z[I^n]$ be the linear map given by
  $\partial(f(q)\cdot\bi)=f'(1)\bi$, where $f'(1)$ is the derivative of $f(q)\in\A$
  evaluated at $q=1$. Then $\partial\Ch
  S^\blam=\sum_\bi\deg_{e,\bi}(\blam)\cdot\bi$.
  The KLR idempotents are orthogonal, so
  $\gdim D^\bmu_\bi=\overline{\gdim D^\bmu_\bi}$ since
  $(D^\bmu)^\circledast\cong D^\bmu$. Therefore, $\partial\Ch D^\bmu=0$.
  Hence, applying $\partial$ to the formula for $\Ch S^\blam$ shows that
  \begin{equation}\label{E:ChSmu}
    \sum_{\bi\in I^n}\deg_{e,\bi}(\blam)\cdot\bi =\partial\Ch S^\blam
         =\sum_{\bi\in I^n}\sum_{\bmu\in\Klesh}
               d_{\blam\bmu}'(1)\dim D_\bi^\bmu \cdot\bi.
  \end{equation}
  Consequently, $\deg_{e,\bi}(\blam)=\sum_\bmu d'_{\blam\bmu}(1)\dim D^\bmu_\bi$.
  So far we have worked over an arbitrary field. If $F=\C$ then
  $d_{\blam\bmu}(q)\in\N[q]$, by \autoref{P:CanonicalBasis}, so that
  $d'_{\blam\bmu}(1)\ge0$. Therefore, $\deg_{e,\bi}(\blam)\ge0$ as claimed.
  (In fact, by \autoref{T:Adjustment}, the right-hand
  side of \autoref{E:ChSmu} is independent of~$F$, as it must be.)

  \autoref{T:Jantzen} shows that taking the $\p$-adic valuation of the Gram
  determinant of $\UnS^\blam$ leads to the Jantzen sum formula for~$\UnS^\blam$.
  Therefore, \autoref{E:ChSmu} suggests that
  \begin{equation}\label{E:GradingFiltration}
    \sum_{k>0}[J_k(\UnS^\blam_\C)]
          =\sum_{\bmu\gdom\blam} d'_{\blam\bmu}(1)[\UnD_\C^\bmu],
  \end{equation}
  where we use the notation of \autoref{T:Jantzen}. That is,
  \autoref{T:IntegralDet} corresponds to writing the Jantzen sum formula as a
  non-negative linear combination of simple modules. In fact, we have not done
  enough to prove \autoref{E:GradingFiltration}. (One way to do this would be
  to establish analogous statements for the Gram determinants of the Weyl modules of
  the cyclotomic Schur algebras~\cite{DJM:cyc}.) Nonetheless,
  \autoref{E:GradingFiltration} is true, being proved by
  Ryom-Hansen~\cite[Theorem~1]{RyomHansen:Schaper} in level one and by
  Yvonne~\cite[Theorem~2.11]{Yvonne:Conjecture} in general.

  A better interpretation of \autoref{E:ChSmu} is in terms of
  \textit{grading filtrations}~\cite[\S2.4]{BGS:Koszul}. Let
  $\dot\R=\ZHom_{\R}(Y,Y)$, where $Y=\bigoplus_{\bmu\in\Klesh}Y^\bmu$ is
  a progenerator for $\R$. Then $\dot\R$ is a graded basic algebra for
  $\R$ and the functor
  $$\Fun\map{\rep(\R)}\rep(\dot\R); M\mapsto \ZHom_{\R}(Y,M),
          \qquad\text{for }M\in\rep(\R),$$
  is a graded Morita equivalence; see, for example,
  \cite[\S2.3-2.4]{HuMathas:QuiverSchurI}. Recall that
  $\car=\(c_{\blam\bmu}(q)\)=\dec^T\circ\dec$ is the Cartan matrix of~$\R$.
  By \autoref{C:CartanSymmetric},
  $c_{\blam\bmu}(q)=\gdim\ZHom_{\R}(Y^\blam,Y^\bmu)$ so that
  $$\gdim\dot\R=\sum_{\blam,\bmu\in\Klesh}c_{\blam\bmu}(q)\in\N[q,q^{-1}].$$

  For the rest of this section assume that $F=\C$. Then
  $c_{\blam\bmu}(q)\in\N[q]$ by \autoref{C:GradedDecomp}. Therefore,
  $\gdim\dot\R\in\N[q]$ so that $\dot\R$ is a positively graded algebra.  Let
  $\dot M=\bigoplus_{d=a}^z\dot M_d$ be a graded $\dot\R$-module. The \textbf{grading
  filtration} of $\dot M$ is the filtration
  $\dot M=G_a(\dot M)\supseteq G_{a+1}(\dot M)\supseteq\dots\supseteq
                G_z(\dot M) \supset0,$
  where
  $$G_d(\dot M)=\bigoplus_{k\ge d}\dot M_k.$$
  Then $G_r(\dot M)$ is a graded $\dot\R$-module precisely because~$\dot\R$ is
  \textit{positively} graded. The grading filtration of an $\Rn$-module $M$ is
  the filtration given by $G_r(M)=\Fun^{-1}(G_r(\Fun(M)))$, for $r\in\Z$. By
  \autoref{C:GradedDecompDegrees}, $S^\blam=G_0(S^\blam)$ and $G_r(S^\blam)=0$
  for $r>\defect\blam$.

  For $\blam\in\Parts$ and $\bmu\in\Klesh$ write
  $d_{\blam\bmu}(q)=\sum_{r\ge0}d_{\blam\bmu}^{(r)}\, q^r$, for
  $d_{\blam\bmu}^{(r)}\in\N$.

  \begin{Lemma}\label{L:GradingFiltration}
    Suppose that $F=\C$ and $\blam\in\Parts$. If $0\le r\le\defect\blam$ then
    $$G_r(S^\blam)/G_{r+1}(S^\blam)
      \cong\bigoplus_{\bmu\in\Klesh}\(D^\bmu\<r\>\)^{\oplus d^{(r)}_{\blam\bmu}}.$$
  \end{Lemma}

  \begin{proof}This is an immediate consequence of the definition of the grading
    filtration and \autoref{C:GradedDecomp}.
  \end{proof}

  Comparing this with \autoref{E:GradingFiltration} suggests that
  $J_r(\UnS^\blam)=G_r(S^\blam)$, for $r\ge0$. Of course, there is no reason to
  expect that $J_r(\UnS^\blam)$ is a graded submodule of~$S^\blam$. Nonetheless,
  establishing a conjecture of Rouquier~\cite[(16)]{LLT}, Shan has proved the
  following when $\Lambda$ is a weight of level~$1$.

  \begin{Theorem}[\protect{Shan~\cite[Theorem~0.1]{Shan:JantzenFiltrations}}]
    \label{T:Shan}
    Suppose that  $F$ is a field of characteristic zero, $\Lambda=\Lambda_0$,
    and that $\blam\in\Parts$. Then $J_r(\UnS^\blam)=G_r(S^\blam)$ is a graded
    submodule of~$S^\blam$ and
    $[J_r(\UnS^\blam)/J_{r+1}(\UnS^\blam):D^\mu\<s\>]
              =\delta_{rs}d_{\blam\bmu}^{(r)}$,
    for all $\bmu\in\Klesh$ and $r\ge0$.
  \end{Theorem}

  Shan actually proves that the Jantzen, radical and grading filtrations of
  graded Weyl modules coincide for the Dipper-James $v$-Schur
  algebras~\cite{DJ:Schur}. This implies the result above because the Schur
  functor maps Jantzen filtrations of Weyl modules to Jantzen filtrations of
  Specht modules. There is a catch, however, because Shan remarks that it is
  unclear how her geometrically defined grading relates to the grading on the
  $v$-Schur algebra given by Ariki~\cite{Ariki:GradedSchur} and hence to the KLR
  grading on~$\R$. As we now sketch, \autoref{T:Shan} can be deduced from Shan's
  result using recent work.

  Since Shan's paper cyclotomic quiver Schur algebras have been
  introduced for arbitrary dominant weights
  \cite{HuMathas:QuiverSchurI,StroppelWebster:QuiverSchur,Ariki:GradedSchur},
  thus giving a grading on all of the cyclotomic Schur algebras
  introduced by Dipper, James and the author~\cite{DJM:cyc}.  The key
  point, which is non-trivial, is that the module categories of the
  cyclotomic quiver Schur algebras are Koszul.  When $e=\infty$ this is
  proved in \cite{HuMathas:QuiverSchurI} by reducing to parabolic
  category~$\O$ for the general linear groups, which is known to be
  Koszul by~\cite{BGS:Koszul,Backelin:Koszul}.  Using similar ideas,
  Maksimau~\cite{Maksimau:QuiverSchur} proves that Stroppel and
  Webster's cyclotomic quiver Schur algebras are Koszul for
  arbitrary~$e$ by using \cite{RouquierShanVaragnoloVasserot} to reduce
  to affine parabolic category~$\O$.

  As the module categories of the cyclotomic quiver Schur are Koszul, an
  elementary argument \cite[Proposition~2.4.1]{BGS:Koszul} shows that the
  radical and grading filtrations of the graded Weyl modules of these algebras
  coincide. By definition, the analogue of \autoref{L:GradingFiltration}
  describes the graded composition factors of the grading (=radical) filtrations
  of the graded Weyl modules --- compare with
  \cite[Corollary~7.24]{HuMathas:QuiverSchurI} when $e=\infty$ and
  \cite[Theorem~1.1]{Maksimau:QuiverSchur} in general.  The graded Schur
  functors of \cite{HuMathas:QuiverSchurI,Maksimau:QuiverSchur} send graded
  Weyl modules to graded Specht modules, graded simple modules to graded simple
  $\R$-modules (or zero), grading filtrations to grading filtrations and Jantzen
  filtrations to Jantzen filtrations. Combining these facts with Shan's
  work~\cite{Shan:JantzenFiltrations} implies \autoref{T:Shan} when
  $\Lambda=\Lambda_0$. We note that the $v$-Schur algebras were first shown to
  be Koszul by Shan, Varagnolo and
  Vasserot~\cite{ShanVaragnoloVasserot:KoszulDualityKacMoody}. It is also
  possible to match up Shan's grading on the $v$-Schur algebras with the
  gradings of~\cite{Ariki:GradedSchur,StroppelWebster:QuiverSchur} using the
  uniqueness of Koszul gradings~\cite[Proposition~2.5.1]{BGS:Koszul}.  As these
  papers use different conventions, it is necessary to work with the graded
  Ringel dual.

  The obstacle to extending \autoref{T:Shan} to arbitrary weights
  $\Lambda\in P^+$ is in showing that the Jantzen and radical
  (=grading) filtrations of the graded Weyl modules of the cyclotomic quiver
  Schur algebras coincide. As the cyclotomic quiver Schur algebras are Koszul it
  is possible that this is straightforward. It seems to the author, however,
  that it is necessary to generalize Shan's
  arguments~\cite{Shan:JantzenFiltrations} to realize the Jantzen filtration
  geometrically using the language of~\cite{RouquierShanVaragnoloVasserot}.


  \subsection{A deformation of the KLR grading}\label{S:OKLR}
  Following~\cite{HuMathas:SeminormalQuiver}, especially the appendix, we now
  sketch how to use the seminormal basis to prove that $\R\cong\H$ over a field
  (\autoref{T:BKiso}). The aim in doing this is not so much to give a new proof
  of the graded isomorphism theorem.  Rather, we want to build a bridge between
  the KLR algebras and the well-understood semisimple representation theory of
  the cyclotomic Hecke algebras. In \autoref{S:NewBasis} we cross this bridge to
  construct a new graded cellular basis $\{B_{\s\t}\}$ of~$\H$ that is
  independent of the choices of reduced expressions that are necessary in
  \autoref{T:PsiBasis}.

  Throughout this section we consider a cyclotomic Hecke algebra $\H$
  defined over a field~$F$ that has Hecke parameter $v\in F^\times$ of
  quantum characteristic~$e\ge2$. As in \autoref{S:Quivers}, $\Lambda\in
  P^+$ is determined by a multicharge $\charge\in\Z^\ell$. We set up a
  modular system for studying~$\H=\H(F)$.

  Let $x$ be an indeterminate over $F$ and let $\O=F[x]_{(x)}$ be the
  localization of~$F[x]$ at the principal ideal generated by~$x$. Let
  $K=F(x)$ be the field of fractions of~$\O$.  Let $\HO$ be the
  cyclotomic Hecke algebra with Hecke parameter $t=x+v$, a unit
  in~$\O$, and cyclotomic parameters $Q_l=x^l+[\kappa_l]_t$,
  for $1\le l\le\ell$.  Then $\HK=\HO\otimes_\O K$ is a split semisimple algebra by
  \autoref{T:SSKLRBasis}.  Moreover, by definition,
  $\H=\H(F)\cong\HO\otimes_\O F$, where we consider~$F$ as an $\O$-module by
  letting $x$ act as multiplication by~$0$.

  As the algebra $\HK$ is semisimple, it has a seminormal basis $\{f_{\s\t}\}$ in
  the sense of \autoref{D:SeminormalBasis}. With our choice of parameters, the
  content functions from \autoref{E:content} become
  $$c^\Zcal_r(\s)=t^{2(c-b)}x^l+[\kappa_l+c-b]_t=t^{2(c-b)}x^l+[\c^\Z_r(\s)]_t$$
  if $\s(l,b,c)=r$, for $1\le k\le n$. Then,
  $L_r f_{\s\t}=c^\Zcal_r(\s)f_{\s\t}$, for $(\s,\t)\in\Std^2(\Parts)$. By
  \autoref{C:SeminormalClass}, the basis $\{f_{\s\t}\}$ determines a seminormal
  coefficient system $\balpha=\set{\alpha_r(\t)|\t\in\Std(\Parts)\text{ and
  }1\le r<n}$ and a set of scalars $\set{\gamma_\t|\t\in\Std(\Parts)}$.

  For $\bi\in I^n$ let $\Std(\bi)=\set{\s\in\Std(\Parts)|\bi^\s=\bi}$ be the set
  of standard tableaux with residue sequence~$\bi$. Define
  \begin{equation}\label{E:ResideIdempotent}
    \fo=\sum_{\t\in\Std(\bi)} F_\t.
  \end{equation}
  By definition, $\fo\in\HK$ but, in fact, $\fo\in\HO$. This idempotent lifting
  result dates back to Murphy~\cite{M:Nak} for the symmetric groups. For higher
  levels it was first proved in~\cite{M:seminormal}. In
  \cite{HuMathas:SeminormalQuiver} it is proved for a more general class of
  rings~$\O$.

  \begin{Lemma}[\protect{\cite[Lemma~4.4]{HuMathas:SeminormalQuiver}}]
    \label{L:ResidueIdempotents}
    Suppose that $\bi\in I^n$. Then $\fo\in\HO$.
  \end{Lemma}

  We will see that $\fo\otimes_\O1_F$ is the KLR idempotent
  $e(\bi)$, for $\bi\in I^n$. Notice that $1=\sum_\bi\fo$ and, further, that
  $\fo\fo[\bj]=\delta_{\bi\bj}\fo$, for $\bi, \bj\in I^n$,
  by \autoref{T:SeminormalBasis}.

  As detailed after \autoref{T:BKiso}, Brundan and Kleshchev construct their
  isomorphisms $\R\bijection\H$ using certain rational functions $P_r(\bi)$ and
  $Q_r(\bi)$ in $F[y_1,\dots,y_n]$. The advantage of working with seminormal
  forms is that, at least intuitively, these rational functions ``converge'' and
  can be replaced with ``nicer'' polynomials. The main tool for doing this is
  the following result, generalizing \autoref{L:ResidueIdempotents}.

  Let $M_r=1-t^{-1}L_r+tL_{r+1}$ , for $1\le r<n$. Then
  $M_rf_{\s\t}=M_r^\Zcal(\s)f_{\s\t}$, where
  $M^\Zcal_r(\s)=1-t^{-1}c^\Zcal_r(\s)+tc^\Zcal_{r+1}(\s)$. The constant
  term of $M^\Zcal_r(\s)$ is equal to
  $v^{2c^\Z_r(\s)-1}[1-c^\Z_r(\s)+c^\Z_{r+1}(\s)]_v\ne0$. Consequently, $M_r$
  acts invertibly on $f_{\s\t}$ whenever $\s\in\Std(\bi)$ and
  $1-i_r+i_{r+1}\ne0$ in $I=\Z/e\Z$.  This
  observation is part of the proof of part~(a) of the next result.
  Similarly, set $\rho^\Zcal_r(\s)=c^\Zcal_r(\s)-c^\Zcal_{r+1}(\s)$.
  Then $\rho^\Zcal_r(\s)$ is invertible in~$\O$ if $i_r\ne i_{r+1}$.

  \begin{Corollary}[\protect{%
       Hu-Mathas \cite[Corollary~4.6]{HuMathas:SeminormalQuiver}}]
    \label{C:Inverting}
    Suppose that $1\le r<n$ and $\bi\in I^n$.
    \begin{enumerate}
      \item If $i_r\ne i_{r+1}+1$ then
        $\dfrac1{M_r}\fo=\Sum_{\s\in\Std(\bi)}\frac1{M^\Zcal_r(\s)}F_\s \in\HO$.
      \item If $i_r\ne i_{r+1}$ then
        $\dfrac1{L_r-L_{r+1}}\fo=\Sum_{\s\in\Std(\bi)}\frac1{\rho^\Zcal_r(\s)}F_\s
                       \in\HO$.
  \end{enumerate}
  \end{Corollary}

  The invertibility of $M_r\fo$, when $i_r\ne i_{r+1}+1$, allows us to define
  analogues of the KLR generators of $\R$ in $\HO$. The invertibility of
  $(L_r-L_{r+1})\fo$ is needed to show that these new elements generate~$\HO$.

  Define an embedding $I\hookrightarrow \Z; i\mapsto\i$ by letting $\i$ be the
  smallest non-negative integer such that $i=\i+e\Z$, for $i\in I$.

  \begin{Definition}\label{D:KLRLift}
    Suppose that $1\le r<n$. Define elements
    $\psi_r^\O=\sum_{\bi\in I^n}\psi_r^\O \fo$ in~$\HO$ by
    $$\psi_r^\O \fo=\begin{cases}
      (T_r+t^{-1})\frac{t^{2\i_r}}{M_r}\fo,&\text{if }i_r=i_{r+1},\\
      (T_rL_r-L_rT_r)t^{-2\i_r}\fo, &\text{if }i_r=i_{r+1}+1,\\
      (T_rL_r-L_rT_r)\frac1{M_r}\fo, &\text{otherwise}.
    \end{cases}$$
    If $1\le r\le n$ then define
    $\yo_r=\sum_{\bi\in I^n}t^{-2\i_r-1}(L_r-[\i_r])\fo$.
  \end{Definition}

  We now describe an $\O$-deformation of cyclotomic KLR algebra~$\R$. This is a
  special case of one of the main results of~\cite{HuMathas:SeminormalQuiver},
  which allows greater flexibility in the choice of the ring~$\O$.

  \begin{Theorem}[\protect{Hu-Mathas~\cite[Theorem~A]{HuMathas:SeminormalQuiver}}]
    \label{T:KLRDeformation}
    As an $\O$-algebra, the algebra $\HO$ is generated by the elements
    $$\set{\fo|\bi\in I^n}\cup\set{\psio_r|1\le r<n}\cup\set{\yo_r|1\le r\le n}$$
    subject only to the following relations:
    \bgroup
      \setlength{\abovedisplayskip}{1pt}
      \setlength{\belowdisplayskip}{1pt}
      \begin{equation*}
      \prod_{\substack{1\le l\le\ell\\\kappa_i\equiv i_1\pmod e}}
                (\yo_1-x^l-[\kappa_l-i_1])\fo=0,
      \end{equation*}
      \begin{equation*}
        \fo \fo[\bj] = \delta_{\bi\bj} \fo, \qquad
           {\textstyle\sum_{\bi \in I^n}} \fo= 1, \qquad
            \yo_r \fo = \fo \yo_r,
      \end{equation*}
      \begin{align*}
        \psio_r \fo&= \fo[s_r{\cdot}\bi] \psio_r,
        &\yo_r \yo_s &= \yo_s \yo_r,\\
        \psio_r \yo_{r+1} \fo&=(\yo_r\psio_r+\delta_{i_ri_{r+1}})\fo,&
        \yo_{r+1}\psio_r\fo&=(\psio_r \yo_r+\delta_{i_ri_{r+1}})\fo,\\
        \psio_r \yo_s  &= \yo_s \psio_r,&&\text{if }s \neq r,r+1,\\
        \psio_r \psio_s &= \psio_s \psio_r,&&\text{if }|r-s|>1,
    \end{align*}
    \begin{align*}
      (\psio_r)^2\fo &=\begin{cases}
        (\dyo[\eps]_r-\yo_{r+1})(\dyo[\spe]_{r+1}-\yo_r)\fo,&
                   \text{if }i_r\leftrightarrows i_{r+1},\\
        (\dyo[\eps]_r-\yo_{r+1})\fo,& \text{if }i_r\rightarrow i_{r+1},\\
        (\dyo[\spe]_{r+1}-\yo_r)\fo,& \text{if }i_r\leftarrow i_{r+1},\\
           0,&\text{if }i_r=i_{r+1},\\
           \fo,&\text{otherwise,}
    \end{cases}
    \end{align*}
    and
    $\big(\psio_r\psio_{r+1}\psio_r-\psio_{r+1}\psio_r\psio_{r+1}\big)\fo$ is
    equal to
    \begin{align*}
    \begin{cases}
        (\dyo[\eps]_r+\dyo[\eps]_{r+2}-
        \dyo[\eps]_{r+1}-\dyo[\spe]_{r+1})\fo,
                &\text{if }i_{r+2}=i_r\rightleftarrows i_{r+1},\\
        -t^{\eps}\fo,&\text{if }i_{r+2}=i_r\rightarrow i_{r+1},\\
        \fo,&\text{if }i_{r+2}=i_r\leftarrow i_{r+1},\\
        0,&\text{otherwise},
      \end{cases}
    \end{align*}
  \egroup
  where $\rho_r(\bi)=\i_r-\i_{r+1}$ and $\dyo[d]_r=t^{2d}\yo_r+t^{-1}[d]$,
  for $d\in\Z$.
  \end{Theorem}

  The statement of \autoref{T:KLRDeformation} is slightly different to
  \cite[Theorem~A]{HuMathas:SeminormalQuiver} because we are using a
  different choice of modular system $(K,\O,F)$ and because
  \autoref{D:HeckeAlgebras} renormalises the quadratic relations for the
  generators~$T_r$ of~$\HO$, for $1\le r<n$.


  The strategy behind the proof of \autoref{T:KLRDeformation} is quite simple:
  we compute the action of the elements defined in \autoref{D:KLRLift} on the
  seminormal basis and use this to verify that they satisfy the relations in the
  theorem. To bound the rank of the algebra defined by the presentation in
  \autoref{T:KLRDeformation} we essentially count dimensions. By specializing
  $x=0$, we obtain \autoref{T:BKiso} as a corollary of
  \autoref{T:KLRDeformation}.

  To give a flavour of the type of calculations that were used to verify that
  the elements in \autoref{D:KLRLift} satisfy the relations in
  \autoref{T:KLRDeformation}, for $\s\in\Std(\bi)$ and $1\le r<n$ define
  \begin{equation}\label{E:beta}
  \beta_r(\s)= \begin{cases}
    \dfrac{\alpha_r(\s)t^{2\i_r}}{M_r^\Zcal(\s)}, &\text{if }i_r=i_{r+1},\\[1pt]
    \alpha_r(\s)\rho^\Zcal_r(\s)t^{-2\i_r}, &\text{if }i_r=i_{r+1}+1,\\
    \dfrac{\alpha_r(\s)\rho^\Zcal_r(\s)}{M^\Zcal_r(\s)},
        &\text{otherwise},
  \end{cases}
  \end{equation}
  Then \autoref{T:SeminormalBasis} easily yields the following.

  \begin{Lemma}\label{L:PsiExpansion}
    Suppose that $1\le r<n$ and that $(\s,\t)\in\SStd(\Parts)$. Set
    $\bi=\bi^\s$, $\bj=\bi^\t$, $\u=\s(r,r+1)$ and $\v=\t(r,r+1)$. Then
    $$
    \psio_rf_{\s\t} = \beta_r(\s)f_{\u\t}
        -\delta_{i_ri_{r+1}}\frac1{\rho^\Zcal_r(\s)}f_{\s\t}.
    $$
    Moreover, if $\s(l,b,c)=r$ then
    $$y^{\<d\>}_rf_{\s\t}=t^{-1}\(t^{2(c-b+d-i_r)}x^l+[c^\Z_k(\s)+d-\i_r]\)f_{\s\t},
    $$
    for $1\le r\le n$ and $d\in\Z$.
  \end{Lemma}

  Armed with \autoref{L:PsiExpansion}, and \autoref{D:SNCS}, it is an easy
  exercise to verify that all of the relations in \autoref{T:KLRDeformation}
  hold in~$\HO$. For the quadratic relations, \autoref{L:PsiExpansion} implies
  that $(\psio_r)^2f_{\s\t}=0$ if $\s\in\Std(\bi)$ and $i_r=i_{r+1}$ whereas if
  $i_r\ne i_{r+1}$ then $(\psio_r)^2f_{\s\t}=\beta_r(\s)\beta_r(\u)f_{\s\t}$,
  where $\u=\s(r,r+1)$.  The quadratic relations in \autoref{T:KLRDeformation}
  now follow using \autoref{E:beta} and \autoref{L:PsiExpansion}. For example,
  suppose that $i_r\rightarrow i_{r+1}$ and $\s\in\Std(\bi)$. Pick nodes
  $(l,b,c)$ and $(l',b',c')$ such that  $\s(l,b,c)=r$ and $\s(l',b',c')=r+1$. Then,
  using \autoref{L:PsiExpansion} and \autoref{D:SNCS},
  $$ (\psio_r)^2f_{\s\t}=t^{-2\i_{r+1}}\beta_r(\s)\beta_r(\u)f_{\s\t}
           =t^{-2\i_{r+1}}M_r^\Zcal(\u)f_{\s\t}
  $$
  On the other hand, by \autoref{L:PsiExpansion},
  $(y^{\<1+\rho_r(\bi)\>}_r-\yo_{r+1})$ acts on $f_{\s\t}$ as multiplication by
  the same scalar. It follows that
  \begin{align*}
    (\psio_r)^2\fo
         &=(\psio_r)^2\sum_{\s\in\Std(\bi)}\frac1{\gamma_\s}f_{\s\s}
          =(y^{\<1+\rho_r(\bi)\>}_r-\yo_{r+1})
             \sum_{\s\in\Std(\bi)}\frac1{\gamma_\s}f_{\s\s}\\
         &=(y^{\<1+\rho_r(\bi)\>}_r-\yo_{r+1})\fo
  \end{align*}
  when $i_r\rightarrow i_{r+1}$.  These calculations are perhaps not very
  pretty, but nor are they are difficult. As indicated by
  \autoref{R:QuiverHecke}, the quadratic relations appear in, and simplify, the
  proof of the deformed braid relations.

\subsection{A distinguished homogeneous basis}\label{S:NewBasis}
  One of the advantages of \autoref{T:KLRDeformation} is that it allows us to
  transplant questions about the KLR algebra~$\R$ into the language of
  seminormal bases. \autoref{D:SNCS} defines $*$-seminormal bases, which
  provide a good framework for studying the semisimple cyclotomic Hecke
  algebras. The algebra~$\H$ comes with two cellular algebra automorphisms,~$*$
  and~$\star$, where~$\star$ is the unique anti-isomorphism fixing the
  homogeneous generators of \autoref{D:CycQuiverHecke} and~$*$ is the
  unique anti-isomorphism fixing the generators of
  \autoref{D:HeckeAlgebras}. In general, these automorphisms are
  different.

  \begin{Definition}[\protect{Hu-Mathas~\cite[\S5]{HuMathas:SeminormalQuiver}}]
    \label{D:StarSNCS}
    A \textbf{$\star$-seminormal coefficient system} is a collection of scalars
    $$\bbeta=\set{\beta_r(\t)|\t\in\Std(\Parts)\text{ and }1\le r\le n}$$
    such that $\beta_r(\t)=0$ if $\v=\t(r,r+1)$ is not standard, if
    $\v\in\Std(\Parts)$ then
    $\beta_r(\v)\beta_r(\t)$ is given by the product of the particular
    $\beta$-coefficients in~$\autoref{E:beta}$, and
    $\beta_r(\t)\beta_{r+1}(\t s_r)\beta_r(\t s_rs_{r+1})
          =\beta_{r+1}(\t)\beta_r(\t s_{r+1})\beta_{r+1}(\t s_{r+1}s_r)$,
          and if~\mbox{$|r-r'|>1$}~then
    $\beta_r(\t)\beta_{r'}(\t s_r)=\beta_{r'}(\t)\beta_{r}(\t s_{r'})$
    for $1\le r,r'<n$.
  \end{Definition}

  Exactly as in \autoref{C:SeminormalClass}, a $\star$-seminormal
  coefficient system determines a $\star$-seminormal basis
  $\{f_{\s\t}\}$ that, similar to \autoref{D:SeminormalBasis}, consists
  of elements $f_{\s\t}\in H_{\s\t}$ such that
  $f_{\s\t}^\star=f_{\t\s}$, for $(\s,\t)\in\Std^2(\Parts)$.  The left
  (and right) the action of $\psio_r$ on $f_{\s\t}$ is exactly as in
  \autoref{L:PsiExpansion} but where the coefficients come from an
  arbitrary $\star$-seminormal coefficient system~$\bbeta$.

  \autoref{D:StarSNCS} gives extra flexibility in choosing a
  $\star$-seminormal basis. By \cite[(5.8)]{HuMathas:SeminormalQuiver} there
  exists a $\star$-seminormal basis $\{f_{\s\t}\}$ such that the $\psi$-basis of
  \autoref{T:PsiBasis} lifts to a $\psio$-basis $\{\psio_{\s\t}\}$ with the
  property that
  \begin{equation}\label{E:PsiFTriangular}
    \psio_{\s\t}=f_{\s\t}+\sum_{(\u,\v)\Gdom(\s,\t)} r_{\u\v}f_{\u\v},
  \end{equation}
  for some $r_{\u\v}\in K$. In this way we recover \autoref{T:PsiBasis} and with
  quicker proof than the original arguments
  in~\cite{HuMathas:GradedCellular}. Perhaps most significantly, by working
  with~$\HO$ we can improve upon the $\psi$-basis.

  \begin{Theorem}[\protect{Hu-Mathas
      \cite[Theorem~6.2, Corollary~6.3]{HuMathas:SeminormalQuiver}}]
    \label{T:KLbasis}
    Suppose that $(\s,\t)\in\SStd(\Parts)$. There exists a unique
    element $B^\O_{\s\t}\in\HO$ such that
    $$B^\O_{\s\t}=f_{\s\t}
          +\sum_{\substack{(\u,\v)\in\SStd(\Parts)\\(\u,\v)\Gdom(\s,\t)}}
             p_{\u\v}^{\s\t}(x^{-1})f_{\u\v},$$
    where $p_{\u\v}^{\s\t}(x)\in xK[x]$. Moreover,
    $\set{B^\O_{\s\t}|(\s,\t)\in\SStd(\Parts)}$
    is a cellular basis of~$\HO$.
  \end{Theorem}

  The existence and uniqueness of this basis essentially come down to Gaussian
  elimination, although for technical reasons it is necessary to work over the
  $x\O$-adic completion of~$\O$. Proving that $\{B^\O_{\s\t}\}$ is
  a cellular basis is more involved and, ultimately, this relies on the
  uniqueness properties of the $B^\O$-basis elements.

  As the $B^\O$-basis is determined by a $\star$-seminormal basis, the
  basis $\{B^\O_{\s\t}\}$ behaves well with respect to the KLR grading
  on~$\H$. The main justification for using this seminormal basis as a
  proxy for choosing a ``nice'' basis for~$\H$, apart from the fact that
  it works, is that \autoref{T:SSKLRBasis} shows that the natural
  homogeneous basis of the semisimple cyclotomic quiver Hecke algebras
  is a $\star$-seminormal basis.

  In characteristic zero the non-zero polynomials $p_{\u\v}^{\s\t}(x)$ satisfy
  \begin{equation}\label{E:Ppolys}
    0<\deg p^{\s\t}_{\u\v}(x)\le\tfrac12(\deg\u-\deg\s+\deg\v-\deg\t),
  \end{equation}
  whenever $(\u,\v)\Gdom(\s,\t)$ by
  \cite[Proposition~6.4]{HuMathas:SeminormalQuiver}. Moreover, if $\s,\t,\u,\v$
  are all standard tableaux of the same shape then
  $p_{\u\v}^{\s\t}(x)=p_\u^\s(x)p_\v^\t(x)$, where $0<\deg
  p^\s_\u(x)\le\tfrac12(\deg\u-\deg\s)$ and $0<\deg
  p^\t_\v\le\frac12(\deg\v-\deg\t)$, whenever $\u\gdom\s$ and $\v\gdom\t$,
  respectively.

  As the basis $\{B^\O_{\s\t}\}$ is defined over~$\O$ we can reduce
  modulo the ideal~$x\O$ to obtain a basis
  $\{B^\O_{\s\t}\otimes_\O1_K\}$ of~$\H=\H(K)$. This basis is hard to
  compute and we do not know if the elements of
  $\{B^\O_{\s\t}\otimes_\O1_K\}$ are homogeneous in general.
  Nonetheless, it is possible to construct a homogeneous basis
  $\{B_{\s\t}\}$ of $\H$ from $\{B^\O_{\s\t}\}$. If~$\blam\in\Parts$
  then define $B_{\tlam\tlam}$ to be the homogeneous component
  of~$B^\O_{\tlam\tlam}\otimes1_K$ of degree $2\deg\tlam$. More
  generally, for $\s,\t\in\Std(\blam)$ we define $B_{\s\t}=D_\s^\star
  B_{\tlam\tlam}D_\t$, where $D_\s,D_\t\in\H$ are certain homogeneous
  elements in~$\H$. In characteristic zero, $B_{\s\t}$ is the
  homogeneous component of $B^\O_{\s\t}\otimes1_K$ of degree
  $\deg\s+\deg\t$, and all other components are of larger degree. For
  any field, by \autoref{E:PsiFTriangular} and \autoref{T:KLbasis},
  \begin{equation}\label{E:BTriangularity}
    B_{\s\t}=\psi_{\s\t}+\sum_{(\u,\v)\Gdom(\s,\t)}a_{\u\v}\psi_{\u\v},
  \end{equation} for some $a_{\u\v}\in K$ that are non-zero only if
  $\bi^\u=\bi^\s$, $\bi^\v=\bi^\t$ and $\deg\u+\deg\v=\deg\s+\deg\t$.
  Therefore, the $B$-basis resolves the ambiguities of
  \autoref{P:PsiProperties}(b). More importantly, we have the following.

  \begin{Theorem}[\protect{Hu-Mathas~\cite[Theorem~6.9]{HuMathas:SeminormalQuiver}}]
    \label{T:BstBasis}
    Suppose that $K$ is a field. Then $\set{B_{\s\t}|(\s,\t)\in\Std^2(\Parts)}$
    is a graded cellular basis of~$\R$ with weight poset $(\Parts,\gedom)$,
    cellular algebra automorphism~$\star$ and with
    $\deg B_{\s\t}=\deg\s+\deg\t$, for $(\s,\t)\in\Std^2(\Parts)$.
    Moreover, if $(\s,\t)\in\Std^2(\Parts)$ then $B_{\s\t}+\Hlam$ depends only
    on~$\s$ and $\t$ and not on the choice of reduced expressions for the
    permutations $d(\s),d(\t)\in\Sym_n$.
  \end{Theorem}

  By construction, the basis $\{B_{\s\t}\}$ depends on the field~$F$.
  If~$F$ is a field of positive characteristic then $B_{\s\t}$ depends upon the
  choice of the elements~$D_\s$ and $D_\t$, which are uniquely determined
  modulo the ideal~$\Hlam$.

  \subsection{A simple conjecture}\label{S:Conjecture}
  The construction of the basis $\{B^\O_{\s\t}\}$ of $\HO$ in
  \autoref{T:KLbasis}, together with the degree constraints on the
  polynomials $p^{\s\t}_{\u\v}(x)$ in \autoref{E:Ppolys}, is reminiscent
  of the Kazhdan-Lusztig basis~\cite{KL}. There is no known analogue of
  the Kazhdan-Lusztig bar involution in this setting. On the other hand,
  we do require that the basis elements $B_{\s\t}$ are homogeneous,
  which might be an appropriate substitute for being bar invariant in
  the graded setting.  Partly motivated by this analogy with the
  Kazhdan-Lusztig basis, we now define analogues of \textit{cell
  representations} for the $B$-basis.

  The basis $\{B_{\s\t}\}$ of \autoref{T:BstBasis} is a graded cellular
  basis so it defines a new homogeneous basis
  $\set{B_\t|\t\in\Std(\blam)}$ of the graded Specht module~$S^\blam$.
  Let the pre-order $\succeq_B$ on $\Std(\blam)$ be the transitive
  closure of the relation $\dot\succeq_B$ where $\t\dot\succeq_B\v$ if
  there exists $a\in\R$ such that $B_\t a=\sum_\s r_\s B_\s$ with
  $r_\v\ne0$.  (So $\succeq_B$ is reflexive and transitive but not
  anti-symmetric.) Let $\sim_B$ be the equivalence relation on
  $\Std(\blam)$ determined by~$\succeq_B$ so that $\t\sim_B\v$ if and
  only if $\t\succeq_B\v\succeq_B\t$.  For example,
  $\tlam\succeq_B\t\succeq_B\t_\blam$, for all $\t\in\Std(\blam)$.

  Let $\Std[\blam]$ be the set of $\sim_B$-equivalence
  classes in~$\Std(\blam)$. The set $\Std[\blam]$ is partially ordered by
  $\succeq_B$, where $\T\succeq_B\V$ if $\t\succeq_B\v$ for some $\t\in\T$
  and $\v\in\V$. Write $\T\succeq_B\v$ if $\t\succeq_B\v$ for some $\t\in\T$
  and $\T\succ_B\v$ if $\T\succeq_B\v$ and $\v\notin\T$.
  Define $S^\blam_{\T\succeq}$ to be the vector subspace of~$S^\blam$ with basis
  $\set{B_\v|\T\succeq_B\v}$. Similarly, let $S^\blam_{\T\succ}$
  be the vector space with basis $\set{B_\v|\T\succ_B\v}$. The definition
  of~$\succeq_B$ ensures that $S^\blam_{\T\succeq}$ and
  $S^\blam_{\T\succ}$ are both graded $\H$-submodules of $S^\blam$ and that
  $S^\blam_{\T\succ}\subsetneq S^\blam_{\T\succeq}$.
  Therefore, $S^\blam_\T=S^\blam_{\T\succeq}/S^\blam_{\T\succ}$ is a
  graded $\H$-module. By choosing any total order on~$\Std[\blam]$ that
  extends the partial order~$\succeq_B$, it is easy to see that
  $S^\blam$ has a filtration with subquotients being precisely the
  modules $S^\blam_\T$, for $\T\in\Std[\blam]$.

  For $\blam\in\Parts$ let $\T^\blam=\set{\t\in\Std(\blam)|\t\sim_B\tlam}$.  In
  view of \autoref{E:form}, if~$\s,\t\in\Std(\blam)$ and $\<B_\s,B_\t\>\ne0$
  then $\s\sim_B\tlam\sim_B\t$ so that $\s,\t\in T^\blam$. Therefore, $\dim
  D^\blam\le|\T^\blam|$. Of course, if $\blam\notin\Klesh$ then this bound is not
  sharp because $D^\blam=0$  whereas $|\T^\blam|\ge1$.

  \begin{Conjecture}\label{Conjecture}
    Suppose that $F$ is a field of characteristic zero and that
    $\blam\in\Parts$. Then $S^\blam_\T$ is an irreducible $\H$-module,
    for all $\T\in\Std[\blam]$.
  \end{Conjecture}

  As discussed in \cite[\S3.3]{HuMathas:SeminormalQuiver}, and is
  implicit in \autoref{E:ChSmu}, by fixing a composition series for
  $S^\blam$ and using a Gaussian elimination argument, it is possible to
  construct a basis $\{C_\t\}$ of~$S^\bmu$ such that (1) each module in
  the composition series has a basis contained in~$\{C_\t\}$, and (2),
  if $\t\in\Std(\blam)$ then $C_\t=\psi_\t$ plus a linear combination of
  ``higher terms'' with respect to some total order on~$\Std(\blam)$.
  This defines a partition of $\Std(\blam)=T_1\sqcup\dots\sqcup T_z$
  (disjoint union), where the tableaux in the set $T_k$ are in bijection
  with a basis of the $k$th composition factor. Therefore, there exists
  an equivalence relation on~$\Std(\blam)$, together with an associated
  composition series, such that the analogue of \autoref{Conjecture}
  holds for this equivalence relation. Our conjecture attempts to make
  this equivalence relation on $\Std(\blam)$ explicit and canonical.

  If $\mathcal{T}\subseteq\Std(\blam)$ define its character to be
  $\ch\mathcal{T}=\sum_{\t\in\mathcal{T}}q^{\deg\t}\cdot\bi^\t\in\A[I^n]$.
  The point of this definition is that $\ch\mathcal{T}$ is a purely
  combinatorial invariant of~$\mathcal{T}$.  As two examples,
  $\Ch S^\blam=\ch\Std(\blam)$ and $\Ch S^\blam_\T=\ch\T$.

  \begin{Proposition}\label{P:conjecture}
    Suppose that \autoref{Conjecture} holds when $F=\C$.
    \begin{enumerate}
      \item Suppose that $\bmu\in\Klesh$. Then
      $D_\C^\bmu\cong S^\bmu_{\T^\bmu}$ and
      $\Ch D_\C^\bmu=\ch\T^\bmu$.
      \item If $\blam\in\Parts$ and $\T\in\Std[\blam]$ then there is a
      unique pair $(\bnu_\T,d_\T)$ in~$\Klesh\times\N$ such that
      $\ch\T=q^{d_\T}\Ch D_\C^{\bnu_\T}=q^{d_\T}\ch\T^{\bnu_\T}$. Moreover,
      $$d_{\blam\bmu}(q)
             =\sum_{\substack{\T\in\Std[\blam]\\\bnu_\T=\bmu}}q^{d_\T}.$$
    \end{enumerate}
  \end{Proposition}

  \begin{proof}By \autoref{C:GradedSimples}, $D_\C^\bmu\ne0$ since $\bmu\in\Klesh$.
    The irreducible module~$D_\C^\bmu$ is generated by
    $B_{\tmu}+\rad S_\C^\bmu=\psi_{\tmu}+\rad S_\C^\bmu$,
    so~$D_\C^\bmu\cong S^\bmu_{T^\bmu}$ since both modules are irreducible by
    \autoref{Conjecture}. Hence, (a) follows.

    For part~(b), $S^\blam_\T\cong D_\C^\bnu\<d\>$, for some
    $\bnu\in\Klesh$ and $d\in\Z$, because $S^\blam_\T$ is irreducible
    by \autoref{Conjecture}. Therefore,
    $\Ch S^\blam_\T=q^d\Ch D_\C^\bnu$. The uniqueness of
    $(\bnu_\T,d_\T)=(\bnu,d)\in\Klesh\times\Z$ now follows from
    \autoref{T:InjectiveCh}. Moreover, $d\ge0$ by
    \autoref{C:GradedDecomp}. As every composition factor
    of~$S_\C^\blam$ is isomorphic to $S^\blam_\T$, for some
    $\T\in\Std[\blam]$, the formula for $d_{\blam\bmu}(q)$ is now
    immediate.
  \end{proof}

  \autoref{P:conjecture} shows that \autoref{Conjecture} encodes closed formulas
  for the characters and graded dimensions of the irreducible $\H$-modules and
  for the graded decomposition numbers of~$\H$.  For this result to be
  useful we need to first verify \autoref{Conjecture} and then to
  explicitly determine the equivalence relation~$\sim_B$. Our last
  result is a step in this direction.

  \begin{Lemma}
    Suppose that $\s,\t\in\Std(\blam)$ and that
    $\t=\s(r,r+1)$ such that $\bi^\s_{r+1}\ne\bi^\s_r\pm1$, where $1\le r<n$
    and $\blam\in\Parts$. Then~$\s\sim_B\t$.
  \end{Lemma}

  \begin{proof}
    By assumption, either $\s\gdom\t$ or $\t\gdom\s$.  Without loss of
    generality we assume that $\s\gdom\t$.  It follows from
    \autoref{E:BTriangularity}, and \autoref{T:Garnir}, that
    $$B_\s\psi_r=\psi_\t+\sum_\u a_\u\psi_\u=B_\t+\sum_\u b_\u B_\u,$$
    where $a_\u,b_\u\in F$ are non-zero only if
    $\ell(d(\u))<\ell(d(\s))$. Therefore, $\s\succeq_B\t$.  If
    $\bi^\s_{r+1}\ne\bi^\t_r$ then $e(\bi^\s)\psi_r^2=e(\bi^\s)$ by
    \autoref{E:quadratic}, so $\s\sim_B\t$. Now
    consider the more interesting case when $\bi^\s_{r+1}=\bi^\s_r$ or,
    equivalently, $\bi^\s_r=\bi^\t_r$. Using \autoref{E:ypsi},
    $$B_\t y_{r+1} = \(B_\s\psi_r-\sum_\u b_\u B_\u\)y_{r+1}
                   =B_\s(y_r\psi_r+1)-\sum_\u b_\u B_\u y_{r+1}.$$
    In view of \autoref{P:PsiProperties}(c), $B_\s$ appears on the right-hand side
    with coefficient $1$. Hence, $\t\succeq_B\s$ implying that $\s\sim_B\t$ as
    claimed.
  \end{proof}

  Finally, we remark that it is easy to check that \autoref{Conjecture} is true
  in the trivial cases considered in \autoref{Ex:SSSpechts} and
  \autoref{Ex:NilHeckeSpechts}. With considerably more effort, using
  \cite[Lemma~9.7]{BrundanStroppel:KhovanovIII} and results of
  \cite[Appendix]{HuMathas:QuiverSchurI}, it is possible to verify the
  conjecture when $\Lambda\in P^+$ is a weight of level~$2$ and $e>n$. In all of
  these cases, the conjecture can be checked because $B_{\s\t}=\psi_{\s\t}$, for
  all $(\s,\t)\in\Std^2(\Parts)$.

  The $B$-basis, and hence \autoref{Conjecture} and all of the results
  in this section (except that in positive characteristic we can only
  say that $d_\T\in\Z$ in \autoref{P:conjecture}, rather than
  $d_\T\in\N$), make sense over any field. We restrict our conjecture to
  fields of characteristic zero because it would be foolhardy to venture
  into the realms of positive characteristic without strong evidence.
  This said, whether or not our conjecture for the $B$-basis is true, we
  are convinced that, in \textit{all characteristics}, there
  exists a ``canonical'' graded cellular basis~$\{C_{\s\t}\}$ of~$\R$
  such that the analogous version of \autoref{Conjecture} holds for
  the~$\sim_C$ equivalence classes.

  To put it another way, the results of
  \cite[\S3.3]{HuMathas:SeminormalQuiver} show that the KLR-tableau
  combinatorics is rich enough to give closed combinatorial formulas for
  both the graded decomposition numbers and the graded dimensions of the
  irreducible representations of~$\H$. We believe that over any field
  the graded Specht modules have a distinguished homogeneous basis that
  ``canonically'' determines these combinatorial formulas.

\let\u=\uold
\sloppy 

\bibliographystyle{andrew}

\end{document}